\documentclass{amsart}

\usepackage{amsmath}
\usepackage{amssymb}
\usepackage{amsthm}
\usepackage{mathrsfs}
\def\mlprime{\/{\mathsurround=0pt$'$}}

\usepackage{fixmath}
\usepackage[pdftex,colorlinks]{hyperref}

\theoremstyle{plain}
\newtheorem{theorem}{Theorem}[section]
\newtheorem{lemma}[theorem]{Lemma}

\newtheorem{proposition}[theorem]{Proposition}

\newtheorem{corollary}[theorem]{Corollary}

\theoremstyle{definition}

\newtheorem{definition}[theorem]{Definition}

\theoremstyle{remark}
\newtheorem{remark}[theorem]{Remark}

\DeclareMathOperator{\Int}{Int}
\DeclareMathOperator{\Lie}{Lie}

\newcommand{\Vect}[1]{\mathbold{#1}}
\providecommand{\abs}[1]{\lvert#1\rvert}
\providecommand{\norm}[1]{\lVert#1\rVert}

\author{Tal Poznansky}
\title[Linear groups with simple reduced $C^{*}$-algebras]{Characterization of linear groups whose reduced $C^{*}$-algebras are simple}
\begin{document}
\address{Centro di Ricerca Matematica ``Ennio De Giorgi''\\
Scuola Normale Superiore\\
56126 Pisa\\
Italia}
\thanks{The author would like to thank the following
individuals for enlightening conversations: E.~Breuillard, Ts.~Gelander, G.~A.~Margulis,
Sh.~Mozes, G.~A.~So{\u\i}fer, and A.~{\.Z}uk.}
\begin{abstract}
The reduced $C^{*}$-algebra of a countable linear group $\Gamma$ is shown to be simple
if and only if $\Gamma$ has no nontrivial normal amenable subgroups. Moreover, these
conditions are shown to be equivalent to the uniqueness of tracial state on the
aforementioned $C^{*}$-algebra.
\end{abstract}
\maketitle
\tableofcontents
\section{Introduction}
It is of some interest to understand the structure of the reduced $C^{*}$-algebra of a
countable group $\Gamma$. Quite a natural question is: Under which
conditions is the reduced $C^{*}$-algebra of $\Gamma$ simple? This question was
posed, for example, in~\cite{MR87b:22007}. (For a more complete history of this
question, please see~\cite{1123.22004} and the references contained therein.)
Along the same lines, one
might wonder under what circumstances the reduced $C^{*}$-algebra of $\Gamma$ has several
nonproportional traces. The purpose of this note is to resolve these questions in the
case where $\Gamma$ is linear.

This note has as its purpose the demonstration of the following result.
\begin{theorem}\label{thm:main}
Let $\Gamma$ be a countable linear group. The following are equivalent:
\begin{itemize}
\item[(i)] The reduced $C^{*}$-algebra of $\Gamma$ is simple.
\item[(ii)] The reduced $C^{*}$-algebra of $\Gamma$ has a unique trace, up to
normalization.
\item[(iii)] $\Gamma$ has no nontrivial normal amenable subgroups.
\end{itemize}
\end{theorem}

Here is an outline of the structure of the paper: In \S\ref{sec:BCD}, we review the
criteria of Bekka, Cowling, and de la Harpe for the simplicity of the reduced
$C^{*}$-algebra. We then present a modification of their criteria, which we will
eventually verify. In \S\ref{sec:pingpong}, we present a modified ping-pong lemma, which
is better adapted to our purposes: Namely, it will produce ping-pong partners for given
elements, rather than for some powers thereof. In \S\S\ref{sec:tits}
and~\ref{sec:rank1facts}, we collect some lemmas of Tits, and some facts about exterior
powers of representations of rank-one groups, respectively.

In \S\ref{sec:attractorfamilies}, we introduce a tool from dynamics, which will
eventually enable us to make sense of `perturbations' in the context of a Zariski-dense
subgroup. Namely, this tool exploits compactness, and gives a priori, though not
effective, bounds on powers of conjugating elements in certain algebraic expressions.

In \S\ref{sec:CellularDecomposition}, we use the Bruhat decomposition to study the orbit
structure on products of Grassmann varieties. In \S\ref{sec:Transversality}, we apply the results of
\S\ref{sec:CellularDecomposition} to a special case arising from the
study of dynamics of actions on vector spaces over local fields. Roughly speaking, we
explore the subject of transversality for the characteristic subspaces of semisimple
elements of a Zariski-dense subgroup of a simple group defined over a local field.

Eventually, we will encode certain algebraic properties in the dynamical language of
proximality. For example, the `very proximality' of a semisimple element will allow us show
that it lies in a free subgroup (cf.\ the proof of Theorem~\ref{thm:zariski-densePPPs}).
Similarly, we can encode higher-order algebraic information about subgroups of linear
groups using elements proximal with respect to several representations simultaneously.
We tell this story in \S\ref{sec:SimProx}.

In \S\ref{sec:StrongTransversality}, we revisit the subject of transversality for
characteristic subspaces of 
nontorsion, semisimple elements. Using special facts about simply laced groups, we obtain
very precise control over the characteristic subspaces of conjugate semisimple elements
in a Zariski-dense subgroup in any sufficiently large, irreducible representation on a
vector space over a local field. The proof proceeds by reduction modulo a local field of
a group scheme over $\mathbb{Z}$. Hence \S\S\ref{sec:weightspacesforSLgroups} and~\ref{sec:ExteriorPowers}
take place in the setting of complex reductive groups. In \S\ref{sec:ReductionModk}, we
reap the rewards for the case of a general local field. (For arguments in a similar
spirit, cf.~\cite{MR0372054} and~\cite{MR0277536}.)

In \S\ref{sec:ContractiveDynamics}, we consider the problem of finding, for a fixed nontorsion element of a
linear group, a projective representation over a local field with the
property that the fixed element has very contractive dynamics. In \S\ref{sec:Pi-im}, we consider
the case of a semisimple element. In particular, we can write down a very satisfying criterion
for the existence of an irreducible representation of $G(k)$ satisfying the
constraint that a prescribed semisimple element act very proximally. (See Proposition~\ref{prop:key},
below.) In \S\ref{sec:quasiProxUnipotent}, the case of a
nontorsion unipotent.

In the remaining section, \S\ref{sec:proof}, we prove Theorem~\ref{thm:main}, and
offer several related formulations, which might prove of some utility.

\section{Preliminaries and recollections}
\subsection{Notations}
In this paper, all fields mentioned will be made to lie in a universal field $\Omega$---that
is, an algebraically closed field of infinite transcendence degree over its prime field.
All fields under discussion will be infinite.
By the phrase ``a vector space'' we will always mean a vector space over $\Omega$ unless
we specify otherwise, writing for instance ``a vector space over $k$,'' where $k$ is
another field. 
Unless specified otherwise, by the phrase ``an algebraic group,'' we shall intend an affine
algebraic
group defined over $\Omega$, which we implicitly identify with the set of its points over
$\Omega$. Otherwise, we will use the terminology, ``an algebraic group defined over
$k$,'' or ``a $k$-group.'' We will call a $k$-group $G$ absolutely almost simple (resp.\ 
almost $k$-simple)
if it has no proper normal, connected, algebraic subgroups (resp.\  $k$-subgroups).

Let $k$ be a field. A vector space $E$ is said to have a $k$-structure if we have
implicitly associated to $E$ a $k$-submodule $E_{k}$ of $E$ with the property that $E \cong
E_{k}\otimes_{k} \Omega$. If $V$ is an algebraic variety defined over $k$, we shall
denote by $V(k)$ the set of its $k$-points. We will habitually use the following
observations: If $E$ is a vector space with a $k$-structure $E_{k}$, then we can identify
 $\mathrm{GL}(E)(k)$ with $\mathrm{GL}(E_{k})$.  If we denote by $\mathrm{Gr}_{m}(E)$
(resp.\  $P(E)$) the Grassmann variety of $m$-dimensional subspaces of $E$ (resp.\  the
projective space of $E$),
then we can identify $\mathrm{Gr}_{m}(E)(k)$ with $\mathrm{Gr}_{m}(E_{k})$ (resp.\ 
$P(E)(k)$ with $P(E_{k})$).

Let $G$ be a $k$-group, $E$ a vector space endowed with a $k$-structure, $\rho: G \to
\mathrm{GL}(E)$ a $k$-rational representation (i.e., a homomorphism which is a
$k$-morphism). Then $\rho$ is said to be absolutely irreducible if $\rho(G)$ leaves
invariant no proper, nontrivial subspace of $E$. It is said to be irreducible over
$k$ if $\rho(G)$ leaves invariant no proper, nontrivial subspace of $E$ which is defined
over $k$. In other words, the representation $\rho_{k}: G(k) \to \mathrm{GL}(E_{k})$ is
irreducible in the usual sense.

Occasionally, we shall work with local fields $k$. In this case, there might be some
ambiguity between the Zariski topology on a $k$-variety and the topology coming from the
locally compact field. We shall endeavor to be explicit, but we usually mean the
Zariski topology. Therefore, terms like $k$-dense and $k$-open always refer to the
Zariski topology, even when $k$ is local. `Connected' will always mean `connected in the
Zariski topology.' On the other hand, the terms `interior' and `compact' will always
refer to those notions relative to the topology induced from the absolute value on a
normed field. We denote by $\Int X$ the interior of a set $X$.

If we have a field extension, $k \subset\ell$ and a $k$-variety $V$, then confusion might
arise about whether a subset of $X$ is open (or closed) in the Zariski topology coming
from $k$, or from $\ell$ (or from $\Omega$, for that matter). Fortunately, this is not the
case, since the $k$-topology on $V(k)$ coincides with the restriction to
$V(k)$ of the $\ell$-topology on $V(\ell)$.

If $g$ and $h$ are elements of any group
whatsoever, we introduce the notation $\mbox{}^{g}h = ghg^{-1}$ and $h^{g} = g^{-1}hg$.
If $g$ is an element of an algebraic group, we denote by $g_{s}$ and $g_{u}$ its
semisimple and unipotent Jordan components. An algebraic $k$-group is said to be
unipotent if all of its elements are unipotent.
We recall that a connected algebraic $k$-group is called semisimple if it contains no
nontrivial, connected, solvable algebraic subgroup (equivalently, $k$-subgroup). It is
called reductive if it contains no nontrivial, connected, unipotent algebraic subgroup
(equivalently, $k$-subgroup). We refer to~\cite{MR0277536} for the representation theory
of semisimple algebraic groups, and to~\cite{MR1890629} for generalities about root
systems.

If $d$ is a distance on the $k$-points of an affine $k$-variety which is comparable to the
affine distance coming from the absolute value on $k$, then we call $d$ admissible. If
$X$ and $Y$ are subsets of affine $k$-space, and $f: X \to Y$ a function, then we denote
by $\norm{f}$ the supremum $$\sup_{\substack{p, q \in X\\
p\neq q}}\frac{d(f(p), f(q))}{d(p, q)}.$$

Let now $E$ denote a finite-dimensional vector space over a local field $k$, $P$ its
projective space.
If $g$ is a linear automorphism of $E$ (and $\dim E > 1$),
we denote by $A(g) \subset P$ (resp.\  $A'(g)$) the projectivization of the sum of the
eigenspaces of $g$ corresponding to eigenvalues of maximal norm (resp., not of maximal
norm). Write $\mathrm{Cr}(g) = A'(g)\cup A'(g^{-1})$. We note that if $g, h \in \mathrm{GL}(E)$,
then $A(\mbox{}^{g}h) = g
\cdot A(h)$, $A'(\mbox{}^{g}h) = g \cdot A'(h)$; and, for any $z \in \mathbb{N}$,
$A(g^{z}) = A(g)$, $A'(g^{z}) = A'(g)$. Also, if $g_{s} \in \mathrm{GL}(E)$, then
$A(g) = A(g_{s})$, $A'(g) = A'(g_{s})$.
\begin{definition}
We say that $g \in \mathrm{GL}(E)$ is \emph{proximal} if $A(g)$ is a singleton set. We say
that $g$ is \emph{very proximal} if both $g$ and $g^{-1}$ are proximal
\end{definition}
\noindent If $\rho$ is a representation of a group $\Gamma$ on a finite-dimensional
vector space over a local field $k$, and $S$ is a subset of $\Gamma$, then we write
\begin{gather*}
\Omega_{+}(\rho, S) = \{s \in S\mid\rho(s)\text{ is proximal}\}\\
\Omega_{0}(\rho, S) = \{s \in S\mid\rho(s)\text{ is very proximal}\}.
\end{gather*}
We note that if $\rho$ is a $k$-rational representation
of a $k$-group $G$ on a finite-dimensional
vector space (over $\Omega$, that is), then we can regard the restriction of $\rho$
to $G(k)$ as a representation on a finite-dimensional $k$-space. Thus, if $S\subset
G(k)$, then we will use the notation $\Omega_{+}(\rho, S)$, $\Omega_{0}(\rho, S)$
as described above. 

Finally, suppose $E$ is a vector space endowed with a $k$-structure, $P$ its
projectivization. If $h\in
\mathrm{GL}(E)(k)$, it makes sense to talk about $A(h) \subset P(k)$.
It also makes sense, if $g \in \mathrm{GL}(E)$, to talk about $g \cdot A(h) \subset P$.
We caution that this is not the same thing as $A(\mbox{}^{g}h)$, which need not be defined.

\subsection{On the criteria of Bekka, Cowling, and de la Harpe}\label{sec:BCD}
In the paper~\cite{MR96a:22020}, Bekka, Cowling, and de la Harpe provide several
sufficient conditions on a countable
group $\Gamma$ for the reduced $C^{*}$-algebra $C_{r}^{*}(\Gamma)$ of $\Gamma$ to be
simple, and to have a unique normalized trace.

For a countable group $\Gamma$, denote by $\lambda_{\Gamma}: \Gamma \to U(\ell^{2}(\Gamma))$
the left regular representation of $\Gamma$ on the Hilbert space $\ell^{2}(\Gamma)$.
For a square-summable sequence $a = \langle a_{j} \rangle$ of complex numbers, denote by
$\norm{a}_{2}$ the $\ell^{2}$-norm $(\sum_{j = 1}^{\infty} \abs{a_{j}}^{2})^{1/2}$.
The following appears as Lemma~2.1 in~\cite{MR96a:22020}.
\begin{lemma}\label{lem:BCDPana}
If for every finite set $F$ of nonidentity elements
of $\Gamma$, there exist an element $g \in \Gamma$ and a number $C$ such that
\begin{equation}\label{eq:ana}
\left\lVert\sum_{j = 1}^{\infty} a_{j} \lambda_{\Gamma}(g^{-j}hg^{j})\right\rVert
\leqslant C
\norm{a}_{2}
\end{equation}
for each $h \in F$ and every square-summable sequence $a$, then the
reduced $C^{*}$-algebra $C_{r}^{*}(\Gamma)$ is simple, and has a unique trace up to
normalization.
\end{lemma}
In turn the other sufficient conditions of the paper~\cite{MR96a:22020} entail the one
above. The following condition---and the proof of its sufficiency---is a slight
modification of some of those other sufficient conditions.

Given a representation $\rho$ of $\Gamma$ on a vector space over a local field, and a
positive number $c$, denote
by $A^{c}(\rho(g))$ the projectivization of the sum of those generalized eigenspaces of
$\rho(g)$ corresponding to eigenvalues of norm $c$. When $\rho$ is clear from the
context, we shall omit it from our notation.
\begin{lemma}\label{lem:modifiedBCD}
Suppose that for every finite subset $F \subset \Gamma$, there
exists a decomposition $F = F_{p} \cup F_{r}$; a
representation $\rho_{h}$ of $\Gamma$ on a finite-dimensional vector space over a local
field for each element $h \in F_{r}$; and an element $g \in \Gamma$ of infinite order,
such that
\begin{itemize}
\item[(i)] for each $h \in F_{p}$, the subgroup $\langle g, h \rangle$ generated by $g$ and
$h$ is canonically isomorphic to the free product $\langle g \rangle * \langle h
\rangle$; and
\item[(ii)] for each $h \in F_{r}$ and every pair of positive numbers $c$ and $d$, we have
$h \cdot A^{c}(\rho_{h}(g)) \cap A^{d}(\rho_{h}(g)) = \varnothing$.
\end{itemize}
Then $C_{r}^{*}(\Gamma)$ is simple, and has a unique trace up to normalization.
\end{lemma}
\begin{proof}
First, fix $h \in F_{p}$. In the proof of Lemma~2.2 of~\cite{MR96a:22020}, Bekka,
Cowling, and de la Harpe show that there exists $C$ such that
equation~\eqref{eq:ana} holds for $g$, our fixed $h$, and every square-summable sequence
$a$.

Next, applying the proofs of Lemmas~2.3 and~2.4 of~\cite{MR96a:22020} to the set $F_{r}$,
we arrive at the same conclusion for each $h \in F_{r}$. Since $F$ is a finite set, we
can take $C$ large enough to work for all $h \in F$ simultaneously. The result now
follows from Lemma~\ref{lem:BCDPana}.
\end{proof}
\subsection{The ping-pong lemma}\label{sec:pingpong}
The following formulation of the ping-pong lemma is better adapted for our purposes than
the usual formulation. The usual formulation in practice requires one to raise given group
elements to some powers. But the statement of our problem necessitates finding a
ping-pong partner for a \emph{fixed} element.

\begin{lemma}[Ping-pong]\label{lem:pingpong}
Let $L$ be a group generated by a subgroup $K$ and an element $h$,
with $K$ of cardinality exceeding $2$. Assume given a subset $U$ of an $L$-space $X$,
satisfying the following conditions:
\begin{itemize}
\item[(i)] $h\cdot U \neq U$; and
\item[(ii)] For all integers $j$ such that  $h^{j} \neq 1$, and all $g \neq 1$ in $K$, we
have $gh^{j}\cdot U \subset U$.
\end{itemize}
Then $L$ is isomorphic to the free product $K * \langle h\rangle$.
\end{lemma}
\begin{proof}
If $\widetilde{K} \cong K$ and $\langle \widetilde{h}\rangle \cong \langle h\rangle$, then there
is a canonical epimorphism $\eta$ from $\widetilde{K} * \langle \widetilde{h} \rangle$ onto $L$.
Explicitly, $\eta$ sends $\widetilde{g}$ to $g$ for $\widetilde{g} \in \widetilde{K}$,
and sends $\widetilde{h}$ to $h$. To show that the kernel of this epimorphism is trivial, it
is necessary and sufficient to show that any nontrivial reduced word in the free product
$\widetilde{K} * \langle \widetilde{h} \rangle$ remains nontrivial when evaluated in $L$.

Let $w$ be such a word.
Our hypothesis on the order of $K$ guarantees that we can conjugate $w$ to a reduced word
$w'$
which begins \emph{and} ends with nontrivial elements from $\widetilde{K}$. Since $\eta(w) =
1$ if and only if $\eta(w') = 1$, we may assume $w$ has the form
$$w = \widetilde{g_{\ell}}\widetilde{h}^{j_{\ell}}\cdots
\widetilde{h}^{j_{2}}\widetilde{g_{1}}\widetilde{h}^{j_{1}}\widetilde{g_{0}}$$
where $h^{j_{i}} \neq 1$, all $i = 1,\ldots,\ell$, and $g_{i} \in K \setminus \{1\}$,
all $i = 0,\ldots,\ell$.

Consider the union $(h\cdot U) \cup (h^{-1}\cdot U)$. We prove by induction on $m$ that the
subword $w_{m} = \widetilde{g_{m}}\widetilde{h}^{j_{m}}\cdots
\widetilde{h}^{j_{2}}\widetilde{g_{1}}\widetilde{h}^{j_{1}}\widetilde{g_{0}}$ acts (via $\eta$) on this
union by sending it into $U$.

First
of all, $g_{0}\cdot ((h\cdot U) \cup (h^{-1}\cdot U))$ is contained in $U$ by (ii). Next
\begin{equation*}
\begin{split}
\eta(w_{m})\cdot ((h\cdot U) \cup (h^{-1}\cdot U)) & =
 (g_{m}h^{j_{m}})\eta(w_{m-1})\\
&\quad\cdot ((h\cdot U) \cup (h^{-1}\cdot U))\\
& \subset (g_{m}h^{j_{m}})\cdot U\\
& \subset U
\end{split}
\end{equation*}
by the induction hypothesis and another application of (ii).

We conclude that $\eta(w)\cdot ((h\cdot U) \cup (h^{-1}\cdot U)) \subset U$. But the
hypothesis that $h\cdot U
\neq U$ is equivalent to $(h\cdot U) \cup (h^{-1}\cdot U) \not\subset
U$. It follows that $w \notin \ker \eta$.
\end{proof}

In our applications, the subgroup $K$ will always be infinite cyclic.

\subsection{Some fundamental lemmas of Tits}\label{sec:tits}
In this section, we simply collect for reference in the sequel several powerful
observations
of Tits regarding the relationship between the norm of a projective transformation and
proximality. Let $E$ be a finite-dimensional
vector space over a local field $k$, and $P$ its projectivization. All topological
terminology in this section (e.g., interior, compact) will refer to the topology
induced from the locally compact topology of $k$. Given a projective
transformation $b$ on $P$ and a subset $K \subset P$, we understand 
$\norm{p|_{K}}$ to be the norm relative to an admissible distance on $P$.

\begin{lemma}[Lemma 3.5 of~\cite{MR44:4105}]\label{lem:3.5}
A projective transformation has finite norm.
\end{lemma}

\begin{lemma}[Lemma 3.8 of~\cite{MR44:4105}]\label{lem:3.8}
Let $g \in \mathrm{GL}(E)$, let $K \subset P$ be a compact set, and let $q >0$.
\begin{itemize}
\item[(i)] Suppose that $A(g)$ is a singleton and that $K \cap A'(g) = \varnothing$. Then
there exists an integer $N$ such that $\norm{\widehat{g}^{z}|_{K}} < q$ for all $z >
N$; and for every neighborhood $U$ of $A(g)$, there exists an integer $N'$ such that
$\widehat{g}^{z} K \subset U$ for all $z > N'$.
\item[(ii)] Assume that, for some $m \in \mathbb{N}$, we have $\widehat{g}^{m} K \subset
\Int K$ and $\norm{\widehat{g}^{m}|_{K}} < 1$. Then $g$ is proximal and $A(g)
\subset \Int K$.
\end{itemize}
\end{lemma}

For linear subspaces $V$ and $W$ of $P$, let us denote by $V\vee W$ the join of $V$ and
$W$. If $V \vee W = E$ and $V\cap W= 0$,
then we can define a map $\mathrm{proj}(V, W):P\setminus V
\to W$ by sending $p \in P$ to $(\{p\}\vee V)\cap W$.

\begin{lemma}[Lemma 3.9 of~\cite{MR44:4105}]\label{lem:3.9}
Let $g \in \mathrm{GL}(E)$ be semisimple. Let $K$ be a compact subset of $P\setminus A'(g)$.
Set $\pi = \mathrm{proj}(A'(g), A(g))$. Let $U$ be a neighborhood of $\pi(K)$ in $P$.
\begin{itemize}
\item[(i)] The set $\{\norm{\widehat{g}^{z}|_{K}} \mid z \in \mathbb{N}\}$
is bounded.
\item[(ii)] Suppose $N \subset \mathbb{N}$ is any infinite set such that, for any pair
$\lambda, \mu$ of eigenvalues of $g$ whose absolute value is maximum, we have
$(\lambda^{-1}\mu)^{z} \to 1$ as $z$ tends to infinity along $N$. Then $\widehat{g}^{z}K
\subset U$ for almost all $z \in N$.
\end{itemize}
\end{lemma}

\subsection{Attractor families}\label{sec:attractorfamilies}
When working with a given projective transformation, it sometimes helps to conjugate by
a projective transformation with more favorable dynamical properties. If we conjugate
the former element
by a large enough power of the latter, the dynamics of the resulting transformation
resemble the latter more than the former. The content of this section is a result which
allows us bound this ``large enough'' at the expense of introducing non-constructiveness.

Let $k$ be an infinite field, with local
extensions $k_{1}, \ldots,
k_{r}$. Let $G$ be a connected $k$-group. Let $\rho_{1}, \ldots, \rho_{r}$ be absolutely
irreducible
representations of $G$ on vector spaces $E_{1}, \ldots, E_{r}$, respectively, and
suppose $\rho_{i}$ is $k_{i}$-rational, $i = 1, \ldots, r$. Denote by $P_{i}$ the
projective space of $E_{i}$.

In this section, we shall be concerned with the compact Hausdorff topological
space $X = P_{1}(k_{1})\times\cdots\times P_{r}(k_{r})$. For brevity, we shall write $\phi_{i}
= \rho_{i}(\phi)$ when $\phi \in G(k)$.
We wish to associate to an
element $\phi \in \bigcap_{i=1}^{r} \Omega_{+}(\rho_{i}, G(k))$ an
open subset of $X$. For such an element $\phi$, set
$$\mathcal{O}_{\phi} = (P_{1}(k_{1})\setminus A'(\phi_{1}))\times\cdots\times
(P_{r}(k_{r}) \setminus A'(\phi_{r})).$$

\begin{definition}
We shall say a subset $\Phi$ of $\bigcap_{i=1}^{r} \Omega_{+}(\rho_{i}, G(k))$ is an
\emph{attractor family} if
$\{\mathcal{O}_{\phi}\mid \phi \in \Phi\}$ covers $X$.
\end{definition}
\noindent We remark that the notion of an attractor family is defined relative to a given
family of representations. This will play an important role later, but we suppress
mention of the family of representations when no ambiguity can arise.

We also remark that a conjugate of an attractor family is again an attractor family.

Suppose $S \subset G(k)$ is a Zariski-dense subgroup of $G$.

\begin{lemma}\label{simulunnec}
\begin{itemize}
\item[(i)] If  $\mathcal{C}$ is the
conjugacy class in $S$ of some $\phi \in
\bigcap_{i
=1}^{r}\Omega_{+}(\rho_{i}, S)$, then $\mathcal{C}$ contains a \emph{finite} attractor
family.
\item[(ii)] Given an attractor family $\Phi$, and a neighborhood $U \subset X$ of
$$\bigcup_{\phi \in \Phi} A(\phi_{1})\times\cdots\times A(\phi_{r}),$$
then there exists a number
$N'$ such that for any $p \in X$ there exists $\phi \in \Phi$ such
that $(\widehat{\rho_{1}} \times\cdots\times \widehat{\rho_{r}})(\phi^{z}) p \subset U$ for
all $z > N'$.
\end{itemize}
\end{lemma}
\begin{proof}
Let $\phi$ be as in part (i) of the lemma. Suppose, contrary to fact that the
collection
$\{\mathcal{O}_{\phi^{x}}\mid x \in S\}$ does not cover $X$; namely, that there exists
$p \in X\setminus \bigcup_{x \in S}\mathcal{O}_{\phi^{x}}$.
By definition, $p = (p_{1},\ldots, p_{r}) \in X\setminus \bigcup_{x \in S}\mathcal{O}_{\phi^{x}}$ if and only if for every $x
\in S$ there exists
an index $i$ such that $x\cdot p_{i} \subset A'(\phi_{i})$. Let $N(p_{i}, A'(\phi_{i}))$ denote
the 
subvariety $\lbrace x \in G \mid x\cdot p_{i} \subset A'(\phi_{i})\rbrace$. Then $S \subset \bigcup_{i = 1}^{r} N(p_{i},
A'(\phi_{i}))$.

Since $S$ is Zariski-dense, we have also $G \subset \bigcup_{i = 1}^{r} N(p_{i},
A'(\phi_{i}))$. But $G$ is irreducible as an algebraic variety, whence there exists an
index $i$ such that $G = N(p_{i}, A'(\phi_{i}))$. In particular, we have also
$S \subset N(p_{i}, A'(\phi_{i}))$. In other words, $p_{i}$ lies in the subspace
$$\bigcap_{x\in S} A'(\rho_{i}(\phi^{x})).$$ But this intersection is a proper and
$S$-invariant subspace and hence also $G$-invariant. It follows from irreducibility of
$\rho_{i}$ that this intersection is empty, providing our contradiction. Therefore,
$\{\mathcal{O}_{\phi^{x}}\mid x \in S\}$ covers $X$. Part (i) of the lemma now
follows from compactness of $X$.

Suppose now that $\Phi$ and $U$ are as in part (ii) of the lemma. By compactness, it
is necessary and sufficient to prove the statement under the additional hypothesis that
$\Phi$ is finite.
We may assume that $U = \bigcup_{\phi \in \Phi} U_{\phi}$ is a union of ``cubes''
$$U_{\phi} = U(\phi, 1)\times\cdots\times U(\phi, r)$$
with $U(\phi, i)\subset P_{i}(k_{i})$ is a neighborhood of $A(\phi_{i})$, $i = 1,\ldots, r$.
Because $X$ is compact
Hausdorff, it follows that we can take neighborhoods $V_{\phi}$ of $X\setminus\mathcal{O}_{\phi}$
with the
property that $\bigcap_{\phi \in \Phi} V_{\phi}$ is empty. By compactness, we can cover
$X\setminus V_{\phi}$ with a finite number of cubes, each of which lies in the
complement of $\mathcal{O}_{\phi}$. By construction of $\mathcal{O}_{\phi}$, the closure of
the $i$th ``side'' of each of these cubes lies outside $A'(\phi_{i})$.
Part (ii) now follows from Lemma~\ref{lem:3.8}.(i).
\end{proof}

If $E$ is a vector space, then we denote by $E^{\vee}$ the dual space of $E$. We 
begin with the useful observation that if $n = \dim E$ is finite, then the projective space
of $E^{\vee}$ is naturally isomorphic to the Grassmann variety $\mathrm{Gr}_{n-1}(E)$
of
hyperspaces in $E$. Explicitly, this isomorphism is given by sending $t \in P(E^{\vee})$ to
$\ker \check{t}$, where $\check{t} \in E^{\vee}\setminus \{0\}$ is a
representative of the line $t$. Moreover, if $E$ has a
$k$-structure, then this isomorphism is defined over $k$; and if $k$ is local, then the
restriction to the $k$-points is a homeomorphism.

If $\rho: G \to \mathrm{GL}(E)$ is a representation of a group $G$ on a space $E$, let
us denote by $\rho^{\vee}$ the dual representation $G \to \mathrm{GL}(E^{\vee})$. We
return to the notation introduced at the beginning of this section. Let $\rho_{1},\ldots,
\rho_{r}, P_{1},\ldots, P_{r} $ be as above.
\begin{corollary}\label{attractfamily}
Let $\Phi \subset
\bigcap_{i = 1}^{r}\Omega_{0}(\rho_{i}, G(k))$
be an attractor
family for the the family of representations
$\{\rho_{1},\ldots, \rho_{r},\rho_{1}^{\vee},\ldots, \rho_{r}^{\vee}\}$.
Assume given subsets $U_{i}, K_{i} \subset P_{i}(k_{i})$ for $i = 1,\ldots, r$,
where $U_{i}$ is a neighborhood of $\bigcup_{\phi \in \Phi} A(\rho_{i}(\phi))$
and $K_{i}$ is compact with
$K_{i} \cap \bigcup_{\phi \in \Phi} A'(\rho_{i}(\phi^{-1})) = \varnothing$,
all $i = 1,\ldots, r$.

Then there exists a number $N'$ such that for any $\gamma_{i} \in
\Omega_{0}(\rho_{i}, G(k))$ there exists a $\phi \in \Phi$ such that
\begin{gather*}
\phi^{z}\cdot A(\gamma_{i}) \subset U_{i}\\
\phi^{z}\cdot A'(\gamma_{i}^{-1}) \cap K_{i} = \varnothing
\end{gather*}
all $i = 1,\ldots, r$ and $z > N'$.
\end{corollary}
\begin{proof}
The corollary follows from two observations: First, for any representation $\rho$ on a
vector space $E$ over a local field, $\rho(\gamma^{-1})$ is proximal implies that
$\rho^{\vee}(\gamma)$ is proximal. Moreover, under the correspondence $P(E^{\vee})
\cong \mathrm{Gr}_{\dim E -1}(E)$ described in the
previous paragraph, the image in $P(E)$ of the ``hyperspace'' $A(\rho^{\vee}(\gamma))$
is exactly $A'(\rho(\gamma^{-1}))$.

The second observation is that if $K\subset P(E)$ and $w \cap K = \varnothing$, then the
condition $v\cap K = \varnothing$ defines a neighborhood of $w$ in the locally compact
topology on $\mathrm{Gr}_{\dim E
-1}(E)$. The result now follows from Lemma~\ref{simulunnec}.(ii).
\end{proof}

\subsection{Cellular decomposition}\label{sec:CellularDecomposition}
In this section, we take $G$ to be a connected, $k$-split, almost $k$-simple
$k$-group. We begin by recalling some standard facts about parabolic subgroups and the Weyl group.
Fix a maximal $k$-split torus $T$ and a Borel subgroup $B$ defined
over $k$ and containing $T$. We denote by $\Pi$ the set of simple roots relative to the
notion of positivity induced on $\Delta(T, G)$ from $B$. For each simple root
$\alpha \in \Pi$ we have an associated simple reflection $w_{\alpha} \in W =W(T,G)$.
The set $\{w_{\alpha} \mid \alpha \in \Pi\}$ is a minimal generating set for $W$. Relative
to this set of generators, there is a unique element $w_{0} \in W$ whose length in the
word metric $\ell_{\Pi}$ is maximal. In particular, $\ell_{\Pi}(w_{0}) = \abs{\Delta^{+}}$.

We call a $k$-subgroup $P$ of $G$ parabolic if it contains a Borel subgroup defined over
$k$, and \emph{standard} parabolic if it contains our fixed Borel $B$. The notion of
``standard'' therefore is not invariant, but it will be clear from the context which Borel
is the standard one. In any case, every standard parabolic $P$ is of the form
\begin{equation}
P_{I} = \bigcup_{w \in W_{I}} BwB\label{doublecosetdecomp}
\end{equation}
where $I$ is a subset of $\Pi$ and $W_{I}$ the subgroup of $W$ generated by the simple
reflections associated to roots in $I$. From the characterization of $\{w_{\alpha}\}$ as
a \emph{minimal} generating set, it follows that $I \subsetneq \Pi$ implies $W_{I} < W$
and $P_{I} < G$. Also, it follows that every $k$-subgroup containing $B$ is $k$-closed.

We recall the fact that the unipotent radical $B_{u}$ of $B$ is isomorphic as a $k$-variety
to a `pointed' affine space of dimension $\abs{\Delta^{+}}$. Explicitly, the map
$$\Vect{t} = (t_{\beta})_{\beta \in \Delta^{+}} \mapsto \prod_{\beta \in \Delta^{+}}
u_{\beta}(t_{\beta})$$
is an isomorphism of algebraic $k$-varieties, where the product is taken in any fixed
order on $\Delta^{+}$.

Let $U^{-}$ be the unipotent radical of an `opposite' Borel subgroup, $B^{\dot{w}_{0}}$. Having fixed a total ordering on $\Delta^{+}$, we denote by $\pi$ the projection of $U^{-}B \cong U^{-} \times B$ onto the pointed affine space $A^{\abs{\Delta^{+}}}$ of dimension $\abs{\Delta^{+}}$ given by
$$\prod_{\beta < 0} u_{\beta}(t_{\beta})b \mapsto \Vect{t} = (t_{\beta})_{\beta < 0}.$$
Here the product on the left-hand side is taken with respect to the total ordering on $\Delta$. In particular, the leftmost term comes from $U_{-\lambda}$, where $\lambda$ is the highest
root.

\begin{proposition}\label{prop:conjtobigcell}
If $h \in G(k)$ is not central in an absolutely almost simple $k$-subgroup of full rank in
$G$, then $h$ is conjugate to an element of $Bw_{0}B$ by an element of $G(k)$.
\end{proposition}
\begin{proof}
For $w \in W$, let us denote by $X(w)$ the closure of the Bruhat cell $BwB$. We can identify the tangent space to
$X(w)$ at the identity with a certain subspace of the Lie algebra of $G$. For a root $\alpha \in \Delta$, let us denote by $L_{\alpha}$ the
corresponding root subspace of $\Lie G$, and by $X_{\alpha}
\in L_{\alpha}$ the normalized derivation corresponding to conjugation by the root subgroup $U_{\alpha}$. Since $X(w)$ is invariant under the conjugation action
of $T$, it follows that this subspace is spanned by the root subspaces it contains. Notably, if $w \neq w_{0}$, then $T_{1}
X(w) \subset \Lie B \oplus \bigoplus_{-\lambda < \nu < 0} L_{\nu}$, where $\lambda$ is the highest
root of $\Delta$ with respect to the partial ordering induced from
$B$. Let $\{\Vect{e}_{\nu}(a) \mid \nu < 0\}$ be the `standard' basis of the tangent space to $A^{\abs{\Delta^{+}}}$ at $a$.

If $s \in
T$, then the pair $(X(w), s)$ is invariant under the conjugation action of
$T$. It follows that if $w \neq w_{0}$, then $d_{s}\pi$ maps the tangent space $T_{s} X(w)$ into the span of those vectors $\Vect{e}_{\nu}(\Vect{0})$ with $-\lambda < \nu < 0$.

For an element $g \in G$, let us denote by $\phi_{g}$ the morphism $G \to A^{\abs{\Delta^{+}}}$ given by $x \mapsto \pi(\mbox{}^{x}g)$. We shall be interested in the differential $d_{1}\phi_{g}$ of $\phi$ at the identity. If $s \in T \setminus \ker \lambda$, then $d_{1}\phi_{s}(X_{-\lambda})$ is a nonzero multiple of $\Vect{e}_{-\lambda}(\Vect{0})$. But this implies that $d_{s}\pi(T_{s}[s]) = \operatorname{im} d_{1}\phi_{s} \not\subset d_{s}\pi(\bigcup_{w \neq w_{0}} T_{s} X(w))$, whence 
$T_{s}[s] \not\subset \bigcup_{w \neq w_{0}} T_{s} X(w)$. Thus in this case we have that $[s] \not\subset \bigcup_{w \neq w_{0}} X(w)$, namely that $[s]$ meets $Bw_{0}B$.

But the orbit of $-\lambda$ under the Weyl group is the set of all (long) roots in $\Delta$. This implies that $[s]$ meets $Bw_{0}B$ unless $s$ lies in the intersection of the kernels of all (long)
roots. Namely, $[s]$ meets $Bw_{0}B$ unless $Z_{G}(s)$ contains
the connected reductive subgroup $M$ of $G$
generated by those subgroups $U_{\alpha}$ with $\alpha \in \Delta$ (long).
(Put otherwise, $M$ is the smallest absolutely almost simple, closed $k$-subgroup of $G$
containing T.)
Incidentally,
this completes the proof of the proposition in the special case where $h$ is semisimple.

We recall the fact that if $s$ is the semisimple Jordan component of an element $h \in G$, then $[s]$ lies in the closure of $[h]$.
If $w \in W \setminus \{w_{0}\}$, then $\bigcap_{x \in G} \operatorname{Ad} x (T_{1}X(w)) =
0$. Equivalently, the connected component of $\bigcap_{x \in G} \mbox{}^{x} X(w)$ containing the identity is trivial. Combining the preceding remarks leads us to conclude that
$X(w)$ contains no nontrivial unipotent classes. By irreducibility, the same holds for
$\bigcup_{w \neq w_{0}} X(w) = G \setminus Bw_{0}B$. In particular, the proposition holds
in the special case where $h$ is unipotent.

Now let $s$ (resp.\ $u$) be the semisimple (resp.\ unipotent) Jordan component of a fixed but arbitrary element $h \in G$. If the
centralizer $Z_{G}(s)$ of $s$ does not contain a conjugate of $M$, then $B[s]B = G$;
since $[s]$ lies in the closure of $[h]$, it follows that $B[h]B$ is dense.

On the other hand, suppose that $Z_{G}(s)$ contains a conjugate of $M$. Denote by $H$
the connected reductive (in fact simple) subgroup
$Z_{G}(s)^{\circ}$. Upon
conjugation, we may assume
that $B_{H} = B \cap H$ is a Borel subgroup of $H$. If $h$ is
not central in $H$, then $u
\neq 1$. By the preceding remarks,
$B_{H}(h^{H})B_{H} = B_{H}(u^{H})B_{H}$ is dense in $H$. But then $B[h]B$ is dense in
$BHB$. The latter is a closed, irreducible algebraic variety, hence coincides with a Bruhat cell closure
$X(w)$. Since $L_{-\lambda} \subset \Lie M$, it follows that $w = w_{0}$.
\end{proof}

\begin{remark}\label{rem:conjtobigcell}
Suppose $h = su$. As a consequence of the proof of Proposition~\ref{prop:conjtobigcell},
we see that $[h]$ meets $Bw_{0}B$ if either $u \neq 1$, or $s$ is not itself central in
a simple subgroup of full rank in $G$.

Suppose $h = s \in T$ lies in the intersection of the kernels of all
(long) roots of $\Delta(G, T)$. Thus, if $\Delta$ is simply laced, then $s$ is central.
On the other hand, if $\Delta$ is not simply laced, then either $s$ is central; or $s$
lies in the kernel of no short root of $\Delta$, and $Z_{G}(s)^{\circ}$ is
that maximal proper connected simple subgroup of $G$ containing $T$
whose root system coincides with the set of long roots of $\Delta$.

If $\alpha$ is a
short root contained in a $G_{2}$-subsystem (resp.\ $B_{2}$-subsystem) of $\Delta$, then 
there exist a pair of long roots of $\Delta$ whose difference is $3\alpha$ (resp.\
$2\alpha$). It follows that $\alpha(s)$ is a cube root (resp.\ square root) of $1$. Note that any root in a root system of type $C_{3}$ can be written as the sum of two short roots.
If $s$ lies outside the kernels of all short roots of $\Delta$, it follows that $\Delta$
contains no subsystem of type $C_{3}$,
which excludes from consideration the possibility that $\Delta$ is of type $C_{n}$, $n
\geqslant 3$, or of type $F_{4}$. In other words, $B[h]B$ is dense unless 
either:
\begin{itemize}
\item[(i)] $\Delta$ is of type $B_{n}$, $n \geqslant 2$; $k$ is not of characteristic two;
and $\alpha(s) = -1$ for all short roots $\alpha$; or
\item[(ii)] $\Delta$ is of type $G_{2}$; $k$ is not of characteristic three; $\alpha(s) = \epsilon$ for an equilateral
triangle of short roots $\alpha \in \Delta$, where $\epsilon$
is a root of the polynomial $x^{2} + x + 1$; and $\alpha(s) = \epsilon^{2} = \epsilon^{-1}$ for the other three short roots
$\alpha \in \Delta$.
\end{itemize}
\end{remark}

This proposition is useful in analyzing
the diagonal action of $G$ on the variety $G/P \times G/Q$,
where $P$ and $Q$ are standard parabolic subgroups. To be precise, this is the action
$g\cdot (x, y) = (gx, gy)$; or, analogously, $g\cdot (x, y) = (\mbox{}^{g}x,
\mbox{}^{g}y)$, if we prefer to think of points of $G/P$ as subgroups conjugate to $P$.
We write $[(xP, yQ)]$ for the $G$-orbit of a point of $G/P \times G/Q$. 

The fundamental fact about this action is that the $G$-orbits
are in one-to-one correspondence with double cosets $PwQ$, or equivalently, with the set
$W_{P}\backslash W/W_{Q}$ (cf.\ \S3 of~\cite{MR0430094} and \S5 of~\cite{MR0207712}). Here
$W_{P}$ is the subgroup of the Weyl group which indexes
the union in \eqref{doublecosetdecomp}, and $W_{Q}$ is similarly defined. The following
observation will be helpful. For this next lemma, we can relax our standing assumption that
$G$ is almost $k$-simple. It will be sufficient to assume that $G$ is reductive.

\begin{lemma}\label{lem:openorbit}
The orbit $[(P, w_{0}Q)]$ is open and dense in $G/P \times G/Q$.
\end{lemma}
\begin{proof}
The orbit $[(P, w_{0}Q)]$ coincides with the image of $G \times Pw_{0}Q$ under the morphism
$G \times G \to G/P \times G/Q$ sending $(x, y)$ to $(xP, xyQ)$. This morphism is the
composition of the $k$-variety isomorphism $G \times G \to G \times G$ given by $(x, y)
\mapsto (x, xy)$ with the usual quotient morphism $G \times G \to G/P \times G/Q$. The
result follows from the remark that $G \times Pw_{0}Q$ is open.
\end{proof}
\begin{corollary}\label{cor:GkOpenOrbit}
If $k$ is a local field, then $G(k) \cdot (P, w_{0}Q)$ is open and closed in $[(P,
w_{0}Q)](k)$.
\end{corollary}
\begin{proof}
Lemma~\ref{lem:openorbit} implies that the orbit morphism $g \mapsto g \cdot
(P, w_{0}Q)$ is separable. The result follows from paragraph~3.18 of~\cite{MR0316587}.
\end{proof}

\subsection{Transversality}\label{sec:Transversality}
The main aim of this section is to
prove Proposition~\ref{prop:dyslexictits},
which will later enable us to choose elements with favorable properties generically.

Let $G$ be a $k$-group, and let $\rho: G \to
{\mathrm GL}(E)$ be a finite-dimensional representation on a vector space $E$. Suppose
moreover
that $E$ is endowed with a $k$-structure with respect to which $\rho$ is rational. In
this section we also assume that $\rho$ is \emph{irreducible over $k$}---namely, $\rho(G)$
leaves invariant no subspace of $E$ which is defined over $k$.

We shall suppress mention of $\rho$ whenever possible, writing $\rho(g) \Vect{e} = g\cdot
\Vect{e}$ when $g \in G$ and $\Vect{e} \in E$.

Let $V$ and $W$ be subspaces of $E$ defined over $k$ with $V > 0$ and $ W < E$. Let us
recall the following lemma of Tits.

\begin{lemma}[Lemma 3.10 of \cite{MR44:4105}]\label{lem:tits}
The set $\{g \in G \mid g\cdot V \not\subset W\}$
is $k$-open and nonempty.
\end{lemma}

The representation $\rho$
induces a $k$-morphism $G \times G \to \mathrm{Gr}_{\dim V}(E) \times \mathrm{Gr}_{\dim W}(E)$,
where
$\mathrm{Gr}_{m}(E)$ denotes the Grassmann variety of $m$-dimen\-sional subspaces of $E$.
Suppose now that $P$ is the stabilizer in $G$ of $V$ (i.e., $P = \{p \in G \mid
p \cdot V = V\}$), and $Q$ is the stabilizer in $G$ of $W$. If $P$ and $Q$
are parabolic subgroups, then this morphism descends to a $G$-equivariant $k$-morphism
$$G/P \times G/Q \to \mathrm{Gr}_{\dim V}(E) \times \mathrm{Gr}_{\dim W}(E),$$
given explicitly by $(g_{1}P, g_{2}Q) \mapsto (g_{1}\cdot V, g_{2}\cdot W)$. The
morphism is equivariant under the diagonal $G$-action.

Denote by $\mathrm{Inc}$ the
preimage under this morphism of $\{(v,w) \mid v \subset w\}$.
(We intentionally suppress mention of the dimensions $\dim V$ and $\dim W$.)
Note that $\mathrm{Inc}$ is a $k$-closed and $G$-invariant subvariety of $G/P \times
G/Q$.

For the remainder of this section we assume $G$ is $k$-split, $k$-connected, and almost $k$-simple.
Fix a maximal $k$-split torus $T$. Suppose $B$ is a
Borel subgroup of $G$ defined over $k$ and containing $T$. As we discuss above,
$B$ gives rise to a word length on the Weyl group $W =
W(T, G)$. For clarity, we shall refer to the longest word according to $B$ as $w_{B}$
(rather than the more conventional $w_{0}$ used above).

We say that $B$ stabilizes a subset $X \subset E$ if $\rho(b)X = X$ for all $b$ in $B$.

\begin{corollary}\label{cor:opencell}
Take $V$ and $W$ as in the previous lemma, with $V > 0$ and $ W < E$, and take an element $h \in Bw_{B}B$.
Suppose that $B$ stabilizes both $V$ and $W$. Then
$h\cdot V \not\subset W$.
\end{corollary}
\begin{proof}
By hypothesis, $P = \mathrm{Stab}_{G}(V)$ and $Q = \mathrm{Stab}_{G}(W)$ are (standard)
parabolic subgroups whose intersection contains $B$. If $h\cdot V \subset W$, then
$(hP, Q) \in \mathrm{Inc}$. But $\mathrm{Inc}$ is $G$-invariant, and hence contains the
orbit $[(hP, Q)]=[(P, w_{B}Q)]$. Because $\mathrm{Inc}$ is closed and $[(P, w_{B}Q)]$ is
dense in $G/P \times G/Q$, it follows that $G/P \times G/Q  \subset \mathrm{Inc}$.

In particular, we have that $(gP, Q) \in \mathrm{Inc}$ for all $g \in G$. Hence $g\cdot
V \subset W$ for all $g$. The desired conclusion now follows from Lemma \ref{lem:tits}.
\end{proof}

To simplify notation, we write $A(g) = A(\rho(g))$ for $g \in G(k)$, and
similarly for $A'(g)$ and $\mathrm{Cr}(g)$.

\begin{proposition}\label{prop:dyslexictits}
Fix $h\in G(k)$. If $g \in G(k)$ is semisimple and $A(g) \neq E$, then the set
$$U_{h, g} = \{u \in G \mid \mbox{}^{u}h\cdot A(g) \not\subset \mathrm{Cr}(g)\}$$
is $k$-open; and is $k$-dense provided $h$ is not central in an absolutely almost simple
$k$-subgroup of full rank in $G$.
\end{proposition}

Given an element $g \in T(k)$ we
can choose a positive system $\Delta^{+}$ for $\Delta = \Delta(T, G)$ with the
following
property: If $\alpha \in \Delta$ is such that $\abs{\alpha(g)} > 1$, then $\alpha
\in \Delta^{+}$.
If $B$ is the Borel subgroup associated to this choice of positive system, then $B$
stabilizes $A(g)$ as well as $A'(g^{-1})$.
In fact, if $\lambda$ is the highest weight
of $\rho$ (with respect to the notion of positivity
$\Delta^{+}$), then we can write $A(g)$ explicitly as the direct sum of those weight
spaces of the form $E^{\lambda - \mu}$ where $\mu$ is a nonnegative integral sum of
simple roots $\alpha$ with $\abs{\alpha(g)} = 1$. Likewise, we can realize
$A'(g^{-1})$  as the direct sum of the weight spaces \emph{not} of the form $E^{w_{B}\lambda +
\mu}$, where $\mu$ ranges over the same set of sums.

\begin{proof}[Proof of Proposition~\ref{prop:dyslexictits}]
The condition defining $U_{h, g}$ is an open condition. To complete the proof, it
suffices to show that the open sets $U_{1} = \{u \in G \mid \mbox{}^{u}h\cdot A(g) \not\subset
A'(g)\}$ and $U_{2} = \{u \in G \mid \mbox{}^{u}h\cdot A(g) \not\subset
A'(g^{-1})\}$ are nonempty.

Let $T$ be a maximal $k$-split torus containing $g$. By the remarks in the paragraph
following the statement of Proposition
\ref{prop:dyslexictits}, we see that $A(g)$ and $A'(g^{-1})$ are stabilized by a common
$k$-defined Borel subgroup $B$ containing $T$. If $u$ is such that $\mbox{}^{u}h \in
Bw_{B}B$, then it follows from Corollary \ref{cor:opencell} that
$\mbox{}^{u}h \cdot A(g) \not\subset A'(g^{-1})$. It follows from Proposition~\ref{prop:conjtobigcell}
that the set $U_{2}$ of such $u$ is nonempty.

On the other hand it is even easier to see that the set of $u \in G$
such that $\mbox{}^{u}h\cdot A(g) \not\subset A'(g)$ is large. In fact, $U_{1}$
is nonempty simply by virtue of the conjugacy (over $k$) of all Borel
subgroups over $k$. Namely, $h$ can be conjugated into the
stabilizer of $A(g)$, and we certainly have that $A(g) \not\subset A'(g)$.
\end{proof}

In particular, if $h$ is a fixed element of $G(k)$ not central in an absolutely almost simple $k$-subgroup of full rank in $G$, then there is a nonempty
open set of elements $g$ satisfying the condition: If $g \in G(k)$ is proximal, then
$h\cdot A(g) \notin \mathrm{Cr}(g)$. This follows from the dominance of the morphism
$T \times G \to G$ given by $(t, v) \mapsto t^{v}$.

We also note that it is not necessary in Proposition \ref{prop:dyslexictits} to assume that the
element $g$ is semisimple. Jordan decomposition tells us that the Proposition holds without
this assumption
either if the field $k$ is perfect; or if the representation $\rho$ is absolutely
irreducible and $G$ is absolutely almost simple.

We conclude this section with a useful result in a similar spirit. Denote by 
$P$ the projective space of our $G$-module $E$.
The following results immediately from Lemma~\ref{lem:tits}.
\begin{corollary}[Lemma 11 of~\cite{MR613853}]\label{cor:xverse}
If $g_{1}, g_{2} \in G(k)$ are such that $A(g_{1})$ and $A(g_{2})$ are both proper 
subspaces of  $P(k)$,
then the set $$\{x \in G \mid x\cdot A(g_{1}), x \cdot A(g^{-1}_{1}) \not\subset
\mathrm{Cr}(g_{2})\text{ and }
A(g_{2}), A(g^{-1}_{2}) \not\subset x\cdot\mathrm{Cr}(g_{1})\}$$ is $k$-open and
$k$-dense in $G$.
\end{corollary}

\subsection{Unipotent one-parameter subgroups}\label{sec:rank1facts}
Let $G$ be a connected, reductive complex algebraic group.
We fix a maximal torus $T$
and a system $\Pi$ of simple roots for $\Delta = \Delta(T, G)$.

For a root $\alpha \in \Delta$, we denote by $U_{\alpha}$ the corresponding root
subgroup; and by $S_{\alpha}$, the simple rank one
subgroup generated by $U_{\alpha}$ together with $U_{-\alpha}$. Also, if $S \subset
\Pi$, let us denote by $G_{S}$ the subgroup generated by all
$S_{\alpha}$ for all $\alpha \in S$.

Let $E$ be a finite-dimensional module for $G$. Denote by $E^{\mu}$ the weight space
of $E$ corresponding to weight $\mu$.

Suppose that $\mu_{0}$ is a weight of $E$ with the property that $\mu_{0} - \alpha$ is
not a weight of $E$. Consider the
maximal $\alpha$-string through $\mu_{0}$: Namely, the longest possible string
$$\mu_{0}, \mu_{0} + \alpha, \ldots, \mu_{0} + m \alpha$$
of weights of $E$. Since $S_{\alpha}$ stabilizes $\bigoplus_{i = 0}^{m} E^{\mu_{0} +
i \alpha}$, we can infer that if $M$ is any $S_{\alpha}$-irreducible submodule of $E$,
then $M$ is spanned by weight vectors for the action of $G$ on $E$, or equivalently,
with its intersections with the weight spaces $E^{\mu}$.

It follows that for fixed $\alpha \in \Delta$, we can choose a basis $\{\Vect{e}^{\mu}_{i}\}$
of $E$ with the properties that
\begin{itemize}
\item[(i)] For any particular weight $\mu$, the vectors $\Vect{e}^{\mu}_{i}$ for $i = 1, \ldots, \dim E^{\mu}$
form a basis for $E^{\mu}$;
\item[(ii)] If $M$ is any $S_{\alpha}$-submodule of $E$, then $\{\Vect{e}^{\mu}_{i}\}$
contains a basis for $M$; and
\item[(iii)] For all $i = 1, \ldots, \dim E^{\mu}$, either $U_{\alpha}$ fixes
$\Vect{e}^{\mu}_{i}$ or there
exists a $j$, $1 \leqslant j \leqslant \dim E^{\mu + \alpha}$, such that for every $t \in
\mathbb{C}^{\times}$, we have
$$u_{\alpha}(t) \cdot \Vect{e}^{\mu}_{i} = \Vect{e}^{\mu}_{i} + t \Vect{e}^{\mu +
\alpha}_{j} + \text{correction},$$ where the correction vector lies in
$\bigoplus_{i = 2}^{\infty} E^{\mu + i\alpha}$ and is of quadratic magnitude in $t$.
\end{itemize}

Consider the family of nilpotent endomorphisms $u_{\alpha}(t) - I$, where $I$ is the
identity matrix. Given $\Vect{e} \in E^{\mu}$, there is a number $z$ such that the vector-valued
function $(u_{\alpha}(t) - I)^{z} \cdot \Vect{e}$ is not zero identically in $t$, but
this vector lies in the kernel of $u_{\alpha}(t) - I$ for all $t \in \mathbb{C}$. It follows from property~(iii) that
for this $z$,
the vector $(u_{\alpha}(t) - I)^{z} \cdot \Vect{e}$ is
independent of $t$, up to scalar.

For $M$ as in item~(ii), we observe that $M$ is spanned by its intersections
with the weight spaces $E^{\mu}$. Also, if $\dot{w}_{\alpha} \in G$ represents
the simple reflection $w_{\alpha}
\in W$ corresponding to $\alpha$, then $\dot{w}_{\alpha}$ normalizes $S_{\alpha}$. It
follows that $\dot{w}_{\alpha} \cdot (M \cap E^{\mu}) = M
\cap E^{w_{\alpha} \cdot \mu}$.

The set of pairs $(\mu, i)$ indexing the aforementioned basis
$\{ \Vect{e}^{\mu}_{i} \mid 1 \leqslant i \leqslant \dim E^{\mu}\}$ of $E$
possesses additional structure.
By property~(iii) above, we have that either $U_{\alpha}$ fixes
$\Vect{e}^{\mu}_{i}$, or for some $j$ and every $t \in
\mathbb{C}^{\times}$, we have
$u_{\alpha}(t) \cdot \Vect{e}^{\mu}_{i} = \Vect{e}^{\mu}_{i} + tu_{\alpha}(t) \cdot \Vect{e}^{\mu +
\alpha}_{j}$.
In the latter case, let us say that $(\mu, i)$ is \emph{$\alpha$-linked to} $(\mu + \alpha,
j)$. It is true---though no mere consequence of notation---that $(\mu_{1}, i_{1})$ is
$\alpha$-linked to $(\mu_{2}, i_{2})$ if and only if $(\mu_{2}, i_{2})$ is
$(-\alpha)$-linked to $(\mu_{1}, i_{1})$. By an \emph{$\alpha$-chain of $E$} we intend a
sequence of
pairs $\langle(\mu_{j}, i_{j})\rangle_{j = 1, \ldots, \ell}$ such that each
$(\mu_{j}, i_{j})$ is $\alpha$-linked to its successor $(\mu_{j + 1}, i_{j + 1})$.

Let us introduce a bit of notation. Let $\Lambda$ be the weight lattice of $G$.
Let $\pounds: \Lambda \times \mathbb{N} \to \Lambda$ be the projection onto the first
coordinate. We observe that if $I$ is an $\alpha$-chain of $E$, then $\pounds I$ is an
$\alpha$-string of weights of $E$.

If $X \subset E$, we define the set $\pounds_{E} X$ to be the smallest subset $K$ of
$\Lambda$ satisfying
$$X \subset \bigoplus_{\mu \in K} E^{\mu}.$$
For example, for a basis vector $\Vect{e}^{\mu}_{i}$, we have the identity
$\pounds_{E}(S_{\alpha} \cdot \Vect{e}^{\mu}_{i}) = \pounds I$, where $I$ is the
maximal $\alpha$-chain of $E$ containing $(\mu, i)$.
As we observed above, $w_{\alpha}$ acts as an involution on the set $\pounds_{E}(S_{\alpha}
\cdot \Vect{e}^{\mu}_{i})$.

Suppose $M$ is a subspace of a $E$ which is spanned by its intersections with
the weight spaces $E^{\mu}$. Write $M^{\mu} = M \cap E^{\mu}$. By hypothesis, $M =
\bigoplus_{\mu \in \pounds_{E} M} M^{\mu}$.
\begin{lemma}\label{lem:unipotentorbits}
Let $M$ be an $S_{\alpha}$-submodule of $E$.
Suppose $\mu$ and $\nu$ are weights of $\pounds_{E} M$ and satisfying
$$\abs{\langle \mu, \alpha \rangle} \leqslant \langle \nu, \alpha \rangle.$$
Upon restriction to $E^{\mu}$ and projection to $E^{\nu}$, for $t \in \mathbb{C}^{\times}$, $u_{\alpha}(t)$
induces a linear map $M^{\mu} \cap \operatorname{Span} U_{-\alpha}
\cdot M^{\nu} \to M^{\nu}$. This map is injective for all but finitely many values of $t$.

Likewise, the map $M^{\nu} \to M^{\mu}$ induced from $u_{-\alpha}(t)$ is injective for all but finitely many
$t \in \mathbb{C}$.
\end{lemma}
\begin{proof}
Without loss of generality, we assume that
$G$ has rank one, $G = S_{\alpha}$, and that $E$ is irreducible. In this case, the weight
spaces are one-dimensional.

We consider the assertion about the map $E^{\mu} \to E^{\nu}$. Since $w_{\alpha}$
induces an involution on the set of weights of $E$ the hypothesis entails that $E^{\nu}$
is nonzero. The map in question is linear, and hence completely determined by the image
of $\Vect{e}^{\mu}$. Having chosen a basis correctly, this map is given by
$\Vect{e}^{\mu} \mapsto t^{\langle \nu - \mu, \varpi_{\alpha}\rangle} \Vect{e}^{\nu}$,
hence has full rank for almost every $t$.

The second assertion follows by symmetry.
\end{proof}

We conclude with a remark about exterior powers of representations of rank-one groups.
For simplicity we assume that $G$ itself has rank one, $G = S_{\alpha}$. In this case,
the weight lattice is just $\mathbb{Z}\varpi_{\alpha}$, so we shall identify
a weight $\mu = z\varpi_{\alpha}$ with the number $\langle \mu, \alpha \rangle
\norm{\varpi_{\alpha}} =
z / \sqrt{2}$. If $E$ is an irreducible, finite-dimensional representation of $G$, let us
denote by $\Vect{e}^{i}$ a nonzero weight vector of weight $i$.

Suppose now $m \leqslant n = \dim E$ and let $I$ be an $m$-element subset of
$\{\frac{1-n}{\sqrt{2}},
\frac{1-n}{\sqrt{2}} + \sqrt{2}, \ldots, \frac{n - 1}{\sqrt{2}}\}$. Then we define
$$\Vect{e}_{I} = \bigwedge_{i \in I}\Vect{e}^{i}.$$
Denote
by $[I]$ the weight $\sum_{i \in I} i$. We observe that the vectors
$\Vect{e}_{I}$ form a basis of the $m$-fold exterior power $\textstyle{\bigwedge^{\!m}} E$ of
$E$. Moreover, the vector $\Vect{e}_{I}$ is a weight vector of weight $[I]$ for the induced representation of $G$.
\begin{lemma}\label{lem:rank1exteriorpower}
Let $E$ be an irreducible, finite-dimensional module for a rank-one group $G$. Let $\Vect{k}$
be an arbitrary nonzero vector in $E$, and $m < n$. If $E_{0}$ is an
irreducible component of $\textstyle{\bigwedge^{\!m}} E$, then the map $\Vect{x} \mapsto
\Vect{x} \wedge \Vect{k}$ does not annihilate $E_{0}$.
\end{lemma}
\begin{proof}
If $m = n - 1$, then $\textstyle{\bigwedge^{\!m}} E$ is irreducible as a $G$-module, and
the result follows from Lemma~\ref{lem:tits}, for example. Therefore, we assume $m < n
-1$.

Let $\Vect{h} = \sum c_{I} \Vect{e}_{I}$ be a nonzero vector in $E_{0}$.
If $i$ is a weight of $E$, and $I$ is a set of such satisfying $i \notin I$ and $c_{I}
\neq 0$, then $\Vect{e}_{i} \wedge \Vect{h} \neq 0$.

In turn, we can write $\Vect{e}_{i}$ as a sum $\sum c_{t}u_{\alpha}(t) \cdot
\Vect{k}$ or $\sum c_{t}u_{-\alpha}(t) \cdot
\Vect{k}$ over a finite set of times $t$. (Up to exchanging $\alpha$ with $-\alpha$, we
may assume the former.) It follows
that $u_{\alpha}(t) \cdot \Vect{k} \wedge \Vect{h}$ is nonzero for some/almost every $t
\in \mathbb{C}$. The same holds for $\Vect{k} \wedge u_{\alpha}(-t) \cdot \Vect{h}$.
\end{proof}

Combining with Lemma~\ref{lem:tits}, we obtain the following corollary:

\begin{corollary}\label{cor:rank1exteriorpower}
Let $G$, $E$ and $E_{0}$ be as in the Lemma~\ref{lem:rank1exteriorpower}. If $\Vect{k}$
and $\Vect{h}$
are arbitrary nonzero vectors in $E$ and $E_{0}$, respectively, then the set of $g \in G$
satisfying the equation $$\Vect{k} \wedge g \cdot \Vect{h} \neq 0$$ is nonempty and open.
\end{corollary}
\begin{proof}
Lemma~\ref{lem:rank1exteriorpower} implies that the set of $\Vect{x}$ satisfying $\Vect{k} \wedge \Vect{x} =
0$ is a proper subspace.
\end{proof}

\section{Simultaneously proximal elements}\label{sec:SimProx}
We begin this chapter with a lemma due to Margulis and So{\u\i}fer. Let $k_{1},\ldots, k_{r}$ be local fields. 
We denote by $E_{1},\ldots, E_{r}$ vector spaces over the fields $k_{1},\ldots, k_{r}$,
respectively; and by $P_{1},\ldots, P_{r}$ their respective projective spaces.

\begin{lemma}[Compare Lemma~3 of~\cite{MR613853}]\label{lem:rotationofeigs}
Assume given for $i = 1,\ldots, r$ a sequence of semisimple elements $h_{i} \in
\mathrm{GL}(E_{i})$. There exists an infinite subset $N \subset \mathbb{N}$ such
that for any $i$, if $\lambda$ and $\mu$ are eigenvalues of $h_{i}$ whose absolute value is
maximum, we have $(\lambda^{-1}\mu)^{z} \to 1$ as $z$ tends to infinity along $N$.
\end{lemma}
\begin{proof}
We prove by induction on $q$ that there exists a set $N_{q} \subset \mathbb{N}$ with the
properties described in the statement of the lemma satisfied for $h_{i}$, $i = 1,\ldots,
q$. The base case is given in Lemma~3.9.(i) of~\cite{MR44:4105}. 

Assume given $N_{q-1}$ and suppose $\lambda$ and $\mu$ are eigenvalues of $h_{q}$
of absolute value maximal among eigenvalues of $h_{q}$. By local compactness, the unit
ball in $k_{q}$ is compact. Since $\abs{\lambda\mu^{-1}} = 1$, there exists a
sequence $\{z_{m}\}_{m \in \mathbb{N}}$ of elements of $N_{q-1}$ such that $z_{m} \to
\infty$ and $(\lambda^{-1}\mu)^{z_{m}}$ converges as $m \to \infty$. Extracting a
subsequence $\{z_{m_{n}}\}$ such that successive differences $z_{m_{n}} - z_{m_{n-1}}$
tend to infinity in $n$, we see that the set $N_{q} = \{z_{m_{n}} - z_{m_{n-1}} \mid n \in 
\mathbb{N}\}$ satisfies all the conditions set forth in the statement of the lemma.
\end{proof}

\subsection{Density theorems}
Eventually, we will want to choose some elements of $G(k)$ which act very proximally
with respect to several given representations. Speaking informally, we shall call such
elements \emph{simultaneously very proximal} elements, and likewise for simultaneously
proximal elements. In this section, we will show that if we
can find simultaneously proximal elements, then not only are such elements abundant, but
so are simultaneously very proximal elements. In the next section, we will show that
there exist simultaneously proximal elements.

Let $k$ be an infinite field (not
necessarily local); $k_{1},\ldots, k_{r}$, locally compact, valued extensions of $k$; 
$E_{1},\ldots, E_{r}$ vector spaces, with each $E_{i}$ endowed
with a fixed $k_{i}$-structure, $i = 1,\ldots, r$. Denote by $P_{i}$ the projective
space of $E_{i}$. Let
$G$ be a connected, reductive $k$-group, $S$ a Zariski-dense subgroup of $G(k)$, and
assume given a collection $\{\rho_{i}: G \to \mathrm{GL}(E_{i})\}_{i = 1,\ldots, r}$ of 
finite-dimensional representations of $G$. Suppose moreover that
each $\rho_{i}$ is $k_{i}$-rational and irreducible over $k_{i}$.

\begin{proposition}[Compare to Lemma 6 of~\cite{MR613853}]\label{prop:dense}
Suppose that the intersection $\bigcap_{i = 1}^{r} \Omega_{+}(\rho_{i}, S)$ is nonempty.
Then $\bigcap_{i = 1}^{r} \Omega_{+}(\rho_{i}, S)$ is $k$-dense in $G(k)$.
\end{proposition}
\begin{proof}
Let $g$ be an element of the intersection $\bigcap_{i = 1}^{r}\Omega_{+}(\rho_{i}, S)$,
the existence of which
is guaranteed by our
hypothesis. For brevity of notation, we write $g_{i} = \rho_{i}(g) \in
\mathrm{GL}(E_{i})(k_{i})$.
Consider the set $$
U_{i} =
\{x
\in G \mid x\cdot A(g_{i})
\not\subset A'(g_{i})\}.$$
Each $U_{i}$ is a Zariski-open and nonempty subset of $G$.
By connectedness, the intersection
$U = \bigcap_{i=1}^{r} U_{i}$ is Zariski-dense. Since $G$ is reductive, it follows that
$G(k)$---and 
hence also $S$---is Zariski-dense in $G$ (cf.\ Lemma 13.3.9 of~\cite{MR1642713}).
Therefore $S \cap U$ is $k$-dense in $G(k)$.

 Let $h$ be taken from the intersection $S \cap U$. Then for all $i = 1,\ldots, r$ we have
$h\cdot A(g_{i})
\subset P_{i} \setminus A'(g_{i})$. Let $K_{i}$ be a compact subset of $P_{i}(k_{i})$ such
that the singleton $A(g_{i})$ is contained in the interior of $K_{i}$.
(Here we mean `compact' and `interior' in the topology induced from the topology of the
local field $k_{i}$.)

By Lemma~\ref{lem:3.5}, there exists a number $q > 0$ such that $\norm{
\widehat{\rho_{i}}(h)|_{K_{i}}} < q$, all $i = 1,\ldots, r$. Let us write $L_{i} = h\cdot
K_{i}$. Since $L_{i} \subset P_{i} \setminus A'(g_{i})$, by Lemma~\ref{lem:3.8}, it
follows
that there exists an integer $N_{h}$ such that $\norm{\widehat{g_{i}}^{z}|_{L_{i}}} <
1/q$ and $g^{z}\cdot L_{i} \subset \Int K_{i}$, all $i = 1,\ldots, r$ and $z
> N_{h}$. Let $N \subset \mathbb{N}$ be given so that the set $\{g^{z} \mid z \in N\}$
is $k$-connected, and let $N'(h) = \{m \in N \mid m > N_{h}\}$. By
definition of $q$, we have $\norm{\widehat{\rho_{i}}(g^{z}h)|_{K_{i}}} < 1$ and
$g^{z}h \cdot K_{i} \subset \Int K_{i}$, all
$i = 1,\ldots, r$ and $z$ exceeding $N_{h}$.

Therefore Lemma~\ref{lem:3.8}.(ii) implies that $g^{z}h$ lies in
the intersection $\bigcap_{i = 1}^{r} \Omega_{+}(\rho_{i}, S)$ for all $z > N_{h}$.
Denote by $\mathrm{Cl}_{k}(\cdot)$ the $k$-closure operator. 
Since $N\setminus N'(h)$ is finite and the set $\{g^{z}h \mid z \in N\}$
is $k$-connected, it follows that $$g^{z}h \in \mathrm{Cl}_{k}\left(\bigcap_{i = 1}^{r}
\Omega_{+}(\rho_{i}, S)\right)$$
for all $z \in N$. That this is true for every $h$
in the $k$-dense set $S \cap U$ implies that $G \subset \mathrm{Cl}_{k}(\bigcap_{i = 1}^{r} \Omega_{+}(\rho_{i},
S))$. This completes the proof.
\end{proof}

\begin{proposition}[Compare Proposition 3.11 of \cite{MR44:4105}]\label{prop:veryprox}
If $\bigcap_{i = 1}^{r} \Omega_{+}(\rho_{i}, S)$ is $k$-dense in $G(k)$, then so is
$\bigcap_{i = 1}^{r} \Omega_{0}(\rho_{i}, S)$.
\end{proposition}
\begin{proof}
$G$ contains a Zariski-open set of semisimple elements~\cite{MR0180554}. Because $G$ is
reductive, $G(k)$---and hence also $\bigcap_{i = 1}^{r} \Omega_{+}(\rho_{i}, S)$---is Zariski-dense in $G$. Let
$g^{-1}$ be a semisimple element of $\bigcap_{i = 1}^{r} \Omega_{+}(\rho_{i}, S)$. For
convenience we use the notation $g_{i} = \rho_{i}(g)$.

By
Lemma~\ref{lem:tits}, the set
$$\bigcap_{i = 1}^{r}\left(\{x \in G \mid x \cdot A(g_{i}) \not\subset A'(g_{i}^{-1})\} \cap \{x
\in G \mid x^{-1} \cdot A(g_{i}) \not\subset A'(g_{i}^{-1}) \}\right)$$
is nonempty and open.
Let $h \in S$ be an element of this intersection, and set
\begin{align*}
B_{i} &= (h\cdot A'(g_{i})) \vee (h\cdot A(g_{i})\cap A'(g_{i}^{-1})),\\
B'_{i} &= A'(g_{i}) \vee (A(g_{i}) \cap h \cdot A'(g_{i}^{-1})),
\end{align*}
and $U_{i} = \{x \in G \mid x\cdot A(g_{i}^{-1}) \not\subset B_{i}\text{ and } h\cdot
A(g_{i}^{-1}) \not\subset x\cdot B'_{i}\}$.
Because of the conditions set on $h$, we have $B_{i} \neq P_{i}$ and $B'_{i} \neq
P_{i}$, and it follows from Lemma~\ref{lem:tits} that $U = \bigcap_{i = 1}^{r} U_{i}$
is open and dense in $G$. Therefore $S \cap U$ is $k$-dense in $G(k)$. Let $u \in S \cap U$

Set $\pi_{i} = \mathrm{proj}(A'({g_{i}}), A(g_{i}))$ and
$\pi'_{i} = \mathrm{proj}(h\cdot A'({g_{i}}), h\cdot A(g_{i}))$. We have
$u \cdot A(g_{i}^{-1}) \not\subset h \cdot A'(g_{i})$ and $\pi'_{i}(u\cdot A(g_{i}^{-1})) \not\subset
A'(g_{i}^{-1})$. Similarly, $u^{-1}h\cdot A(g_{i}^{-1}) \not\subset A'(g_{i})$ and
$\pi_{i}(u^{-1}h\cdot A(g_{i}^{-1})) \not\subset h\cdot A'(g_{i})$. Let $Y_{i}$ (resp.\  $Y'_{i}$)
be a compact neighborhood of $A(g_{i}^{-1})$ (resp.\  $u^{-1}h\cdot A(g_{i}^{-1})$) such
that 
\begin{gather*}
u\cdot Y_{i} \cap h \cdot A'(g_{i}^{-1}) = \varnothing, \qquad
\pi'_{i}(u\cdot Y_{i}) \cap A'(g_{i}^{-1}) = \varnothing,\\
Y'_{i} \cap A'(g_{i})
= \varnothing, \qquad \pi_{i}(Y'_{i})\cap h\cdot A'(g_{i}^{-1}) = \varnothing.
\end{gather*}
Let $Z_{i}$ (resp.\  $Z'_{i}$) be a compact neighborhood of $\pi'_{i}(u\cdot Y_{i})$ (resp.\ 
$\pi_{i}(Y'_{i})$) in $P_{i}$ whose intersection with $A'(g_{i}^{-1})$ (resp.\  $h\cdot
A'(g_{i}^{-1})$) is empty. By Lemmas~\ref{lem:3.9}.(i) and Lemma~\ref{lem:3.5}, there
exists a number $q > 0$ such that
$$\norm{\widehat{\rho_{i}}(hg^{z}h^{-1}u)|_{Y_{i}}} < q\qquad \text{and}\qquad
\norm{\widehat{\rho_{i}}(g^{z})|_{Y'_{i}}} < q$$
for all $i = 1,\ldots, r$ and $z \in \mathbb{N}$.

Lemma~\ref{lem:rotationofeigs} implies that there exists an infinite subset $N \subset \mathbb{N}$ such
that for every pair $\lambda, \mu \in \Omega(g_{i})$, and every $i = 1,\ldots, r$ we
have $(\lambda\mu^{-1})^{z} \to 1$ as $z \to \infty$ along $N$.
The linear transformation $\rho_{i}(hgh^{-1})$ has the same eigenvalues as $g_{i}$.
Therefore, Lemma~\ref{lem:3.9}.(ii) implies that, for almost all $z \in N$
we have
\begin{equation}\label{eqn:tits3111}
hg^{z}h^{-1}u\cdot Y_{i} \subset Z_{i} \qquad \text{and}\qquad g^{z}\cdot Y'_{i} \subset
Z'_{i}.
\end{equation}
Upon replacing $N$ by a subset, we may assume in addition that the set $\{g^{z} \mid z \in N\}$
is $k$-connected.

By Lemma~\ref{lem:3.8}.(i), we have for almost all $z \in \mathbb{N}$ that
\begin{gather}
\norm{\widehat{g_{i}}^{-z}|_{Z_{i}}} < 1/q,\qquad g^{-z}\cdot Z'_{i} \subset u\cdot
\Int Y'_{i}\label{gat:tits3112}\\
\norm{\widehat{\rho_{i}}(hg^{-z}h^{-1})|_{Z'_{i}}} < 1/q\norm{u^{-1}},
\qquad hg^{-z}h^{-1}\cdot Z'_{i} \subset u\cdot \Int Y'_{i}.\label{gat:tits3113}
\end{gather}
Let $N'(u)$ denote the set of all $z \in N$ such that \eqref{eqn:tits3111},
\eqref{gat:tits3112}, and \eqref{gat:tits3113} hold simultaneously, all $i = 1,\ldots,
r$.

For all $z \in N'(u)$, we have
\begin{gather*}
g^{-z}hg^{z}h^{-1}u \cdot Y_{i} \subset \Int Y_{i}, \qquad
\norm{\widehat{\rho_{i}}(g^{-z}hg^{z}h^{-1}u)|_{Y_{i}}} < 1,\\
u^{-1}hg^{-z}h^{-1}g^{z} \cdot Y'_{i} \subset \Int Y'_{i}, \qquad
\norm{\widehat{\rho_{i}}(u^{-1}hg^{-z}h^{-1}g^{z})|_{Y'_{i}}} < 1;
\end{gather*}
hence, by Lemma~\ref{lem:3.8}.(ii),
$$g^{-z}hg^{z}h^{-1}u \in \bigcap_{i = 1}^{r} \Omega_{0}(\rho_{i}, S).$$
Since $N\setminus N'(u)$ is finite and the set $\{g^{z} \mid z \in N\}$
is $k$-connected, it follows that $$g^{-z}hg^{z}h^{-1}u \in
\mathrm{Cl}_{k}\left(\bigcap_{i = 1}^{r} \Omega_{0}(\rho_{i}, S)\right)$$
for all $z \in N$. That this is true for every $u$
in the $k$-dense set $S \cap U$ implies that $G \subset \mathrm{Cl}_{k}(\bigcap_{i =
1}^{r} \Omega_{0}(\rho_{i}, S))$. This completes the proof.
\end{proof}

\begin{proposition}[Compare to Lemma 9 of~\cite{MR613853}]\label{prop:coset}
Let $H$ be a finite-index subgroup of $S$ and $Hg$ a coset of $H$ in $S$. If 
$\bigcap_{i = 1}^{r} \Omega_{+}(\rho_{i}, S)$ is nonempty, then so is
$\bigcap_{i = 1}^{r} \Omega_{0}(\rho_{i}, Hg)$.
\end{proposition}
\begin{proof}
Since $G$ is $k$-connected, $H$ is Zariski-dense in $G$. The previous two
propositions tell us that we can find an element $h$ of $\bigcap_{i = 1}^{r}
\Omega_{0}(\rho_{i}, H)$. Set $h_{i} = \rho_{i}(h)$. Consider
the set
$$U_{i} = \{x \in G \mid xg \cdot A(h_{i}) \not\subset A'(h_{i})\}
\cap \{x \in G \mid A(h_{i}^{-1}) \not\subset xg \cdot A'(h_{i}^{-1})\}.$$
It follows from Lemma~\ref{lem:tits} that $U = \bigcap_{i = 1}^{r} U_{i}$
is open and dense in $G$. Therefore $H \cap U$ is $k$-dense in $G(k)$. Let $x \in H \cap
U$.

Set $h_{0} = xg \in Hg$. By construction, we have for all $i$ that $h_{0}\cdot A(h_{i})
\not\subset A'(h_{i})$ and $h_{0}^{-1}\cdot A(h_{i}^{-1}) \not\subset A'(h_{i}^{-1})$.
Since the sets $A(h_{i}^{\pm 1})$ are singletons, there exist compact sets $K_{i}$ and
$K_{i}^{-}$ in $P_{i}$ such that
\begin{gather*}
A(h_{i}) \subset \Int K_{i}, \qquad h\cdot K_{i} \subset P_{i}\setminus
A'(h_{i}),\\
A(h_{i}^{-1}) \subset \Int K_{i}^{-}, \qquad h^{-1}\cdot K_{i}^{-} \subset P_{i}\setminus
A'(h_{i}^{-1}),
\end{gather*}
all $i = 1,\ldots, r$. Write $M_{i} = h_{0}\cdot K_{i}$ and $M_{i}^{-} = h_{0}^{-1}\cdot
K_{i}^{-}$. By Lemma~\ref{lem:3.5}, there exists a number $q > 0$ such that
$$\max_{1\leqslant i \leqslant r}\{\norm{\widehat{\rho_{i}}(h_{0})|_{K_{i}}}, \norm{
\widehat{\rho_{i}}(h_{0}^{-1})|_{M_{i}}}\} < q.$$
Since $M_{i} \subset P_{i} \setminus A'(h_{i})$ and $M_{i}^{-} \subset P_{i} \setminus
A'(h_{i}^{-1})$, it follows from Lemma \ref{lem:3.8}.(i) that there exists a number
$N_{1}$ such that $\norm{\widehat{\rho_{i}}(h^{z})|_{M_{i}}} < 1/q$ and
$\norm{\widehat{\rho_{i}}(h^{-z})|_{M_{i}^{-}}} < 1/q$ for each $z > N_{1}$ and $i
= 1,\ldots,r$. Moreover, there exists a number $N_{2} > 0$ such that $h^{z}\cdot M_{i}
\subset \Int K_{i}$ and $h^{-z}\cdot M_{i}^{-}
\subset \Int K_{i}^{-}$ for each $z > N_{2}$ and $i = 1,\ldots, r$.

By choice of $q$, we have that
\begin{gather*}
h^{z}h_{0}\cdot K_{i} \subset \Int K_{i},\qquad
\norm{\widehat{\rho_{i}}(h^{z}h_{0})|_{K_{i}}} < 1,\\
h^{-z}h_{0}^{-1}\cdot K_{i}^{-} \subset \Int K_{i}^{-}, \qquad
\norm{\widehat{\rho_{i}}(h^{-z}h_{0}^{-1})|_{K_{i}^{-}}} < 1,
\end{gather*}
 for
each $z > \max(N_{1}, N_{2})$ and $i = 1,\ldots, r$. So from Lemma~\ref{lem:3.8}.(ii)
it follows that $h^{z}h_{0}, h^{-z}h_{0}^{-1} \in \bigcap_{i = 1}^{r} \Omega_{+}(\rho_{i},
S)$ for each $z > \max(N_{1}, N_{2})$. But $h^{-z}h_{0}^{-1}$ is conjugate to
$h_{0}^{-1}h^{-z} = (h^{z}h_{0})^{-1}$. It follows that for such $z$ we have
$h^{z}h_{0} \in \bigcap_{i = 1}^{r} \Omega_{0}(\rho_{i}, Hg)$.
\end{proof}

We conclude this section with a result on density, though not about proximality. It says
that for a subset of a finitely generated linear group, profinite density is better than
Zariski density, and almost as good as Zariski openness. 
\begin{lemma}\label{lem:YcapHg}
Let $\Gamma$ be a finitely generated, Zariski-dense subgroup of a reductive algebraic group
$G$; $Hg$ a coset in $\Gamma$ of a finite-index subgroup $H < \Gamma$; and $Y$ a Zariski-open subset
of $G$. Then $Y \cap Hg$ is Zariski dense in $G$. 
\end{lemma}
\begin{proof}
Upon replacing $Y$ by $Yg^{-1}$, we may reduce to the case $g = 1$. Since $H$ is of finite
index in $\Gamma$, it is Zariski dense, and we reduce to the case $\Gamma = H = Hg$.
The result now follows from the proof of Proposition~3 of~\cite{MR613853}, for example.
\end{proof}

\subsection{Existence}
As promised, in this section we will produce elements which act proximally
with respect to each of several given representations.

Fix an infinite field $k$ and assume
that $k_{1},\ldots, k_{r}$ is a family of local fields extending $k$.
As before, $G$ is a connected, reductive $k$-group, and $S$ a
Zariski-dense subgroup of $G(k)$.  For $i = 1,\ldots, r$, let
$E_{i}$ be a finite-dimensional vector space endowed with a $k_{i}$-structure.
Assume given a family $\rho_{i}: G \to
\mathrm{GL}(E_{i})$ of $k_{i}$-rational, $k_{i}$-irreducible representations.
\begin{lemma}\label{lem:simprox}
If $\Omega_{+}(\rho_{i}, S) \neq\varnothing$ for all $i = 1,\ldots, r$, then
also $$\bigcap_{i = 1}^{r} \Omega_{+}(\rho_{i}, S)\neq\varnothing.$$
\end{lemma}
\begin{proof}
By induction on $r$. For $r = 1$, there is nothing to prove. Now suppose that both
$\bigcap_{i = 1}^{r - 1} \Omega_{+}(\rho_{i}, S)$ and $\Omega_{+}(\rho_{r}, S)$ are nonempty.
By Proposition~\ref{prop:dense}, we can find semisimple elements $g$ and $h$
such that $g \in \bigcap_{i = 1}^{r - 1} \Omega_{+}(\rho_{i}, S)$ and
$h \in \Omega_{+}(\rho_{r}, S)$. For convenience, we shall write $h_{i} = \rho_{i}(h)$
and $g_{i} = \rho_{i}(g)$ for $i=1,\ldots, r$. Also, set
$\pi_{i} = \mathrm{proj}(A'(h_{i}), A(h_{i}))$ for $i=1,\ldots, r -1$, and
 $\pi_{r} = \mathrm{proj}(A'(g_{r}), A(g_{r}))$.
 
Since the $\rho_{i}$ are irreducible and $S$ is dense, we can find $v \in S$ be such
that $v\cdot A(g_{i})
\subset P_{i}\setminus A'(h_{i})$ for $i = 1,\ldots, r-1$.
Consider the set of elements $u\in G$ simultaneously satisfying the following $r+1$
conditions:
\begin{align}
u\cdot \pi_{i}(v\cdot A(g_{i}))\subset P_{i}\setminus A'(g_{i}),& \qquad i =1,\ldots,
r-1, \nonumber\\
u\cdot A(h_{r})\subset P_{r}\setminus A'(g_{r}).&\nonumber\\
u\cdot A(h_{r}) \not\subset (v^{-1}\cdot A'(h_{r})\cap A(g_{r}))\vee A'(g_{r}).
\label{eqn:2ndprojection}
\end{align} 
By Lemma~\ref{lem:tits}, the set of such $u$ is Zariski-open. Since $S$ is dense, we can choose such an element $u$ in $S$. By~\eqref{eqn:2ndprojection},
we have $v\cdot \pi_{r}(u\cdot A(h_{r})) \subset P_{r}\setminus A'(h_{r})$.

Choose compact sets $L_{i}$ such that
 \begin{gather*}
 L_{i} \subset P_{i}\setminus A'(g_{i})\\
 u\cdot \pi_{i}(v\cdot A(g_{i}))\subset \Int L_{i},
 \end{gather*}
 for $i = 1,\ldots, r-1$; and
 \begin{gather*}
 L_{r} \subset P_{r}\setminus A'(h_{r})\\
 v\cdot \pi_{r}(u\cdot A(h_{r}))\subset \Int L_{r}.
 \end{gather*}
 
 Let $N \subset \mathbb{N}$ be the infinite set whose existence is guaranteed by
 Lemma~\ref{lem:rotationofeigs}, where the set of given semisimple elements is 
 $h_{1}, h_{2},\ldots$, $h_{r-1}, g_{r}$. We choose compact sets $K_{i}$ such that
 \begin{gather*}
 A(g_{i}) \subset \Int K_{i}\\
 v\cdot K_{i} \subset P_{i}\setminus A'(h_{i}) \\
 u\cdot \pi_{i}(v\cdot K_{i})\subset \Int L_{i}
\end{gather*}
for $i = 1,\ldots, r-1$; and 
\begin{gather*}
 A(h_{r}) \subset \Int K_{r}\\
u\cdot K_{r} \subset P_{r}\setminus A'(g_{r}) \\
 v\cdot \pi_{r}(u\cdot K_{r})\subset \Int L_{r}.
\end{gather*}

By Lemma~\ref{lem:3.9}, we have that there exist numbers $N_{0}$ and $q$ such that,
if we set $N' = \{z \in N\mid z>N_{0}\}$, then
\begin{align*}
uh^{z}v\cdot K_{i} \subset \Int L_{i},&\qquad \text{ all $z \in N'$; and}\\
\norm{\widehat{\rho_{i}}(h^{z}) |_{v\cdot K_{i}}} < q &\qquad \text{ for all $z\in \mathbb{N}$}
\end{align*}
for $i = 1,\ldots, r-1$; and 
\begin{align*}
vg^{z}u\cdot K_{r} \subset \Int L_{r},&\qquad \text{ all $z \in N'$; and}\\
\norm{\widehat{\rho_{r}}(g^{z}) |_{u\cdot K_{r}}} < q &\qquad \text{ for all $z\in
\mathbb{N}$}.
\end{align*}

By Lemma~\ref{lem:3.8}.(i), there exists a number $N_{1}$ such that
\begin{gather*}
\norm{\widehat{\rho_{i}}(g^{z})|_{L_{i}}} < (q\norm{\widehat{\rho_{i}}(u)|_{u^{-1}\cdot L_{i}}}
\norm{\widehat{\rho_{i}}(v)|_{K_{i}}})^{-1}\\
g^{z}\cdot L_{i} \subset K_{i}
\end{gather*}
for all $z>N_{1}$ and $i = 1,\ldots, r-1$; and 
\begin{gather*}
\norm{\widehat{\rho_{r}}(h^{z})|_{L_{r}}} < (q\norm{\widehat{\rho_{r}}(v)|_{v^{-1}\cdot L_{r}}}
\norm{\widehat{\rho_{r}}(u)|_{K_{r}}})^{-1}\\
h^{z}\cdot L_{r} \subset K_{r}
\end{gather*}
for all $z>N_{1}$.

It follows from Lemma~\ref{lem:3.8}.(ii) that if
$z > N_{1}$ and $s \in N'$, that $\rho_{i}(g^{z}uh^{s}v)$ is proximal for
$i=1,\ldots,r-1$, as is $\rho_{r}(h^{z}vg^{s}u)$. But $h^{z}vg^{s}u$ is conjugate
to $g^{s}uh^{z}v$. It follows that if $s, z > \max(N_{0}, N_{1})$ are both taken in $N$,
then we have $g^{z}uh^{s}v \in \bigcap_{i=1}^{r}\Omega_{+}(\rho_{i}, S)$.
\end{proof}

Momentarily breaking with our notational convention, we denote by $\widehat{S}$ the
profinite completion of $S$.
\begin{corollary}\label{cor:existenceanddensity}
If $\Omega_{+}(\rho_{i}, S) \neq \varnothing$ for all $i = 1\ldots, r$, then the set $$\bigcap_{i = 1}^{r}
\Omega_{0}(\rho_{i}, S)$$ of simultaneously very proximal elements is Zariski-dense in
$G$, as well as dense in the profinite topology on $\widehat{S}$.
\end{corollary}
\section{Strong transversality for simply laced groups}\label{sec:StrongTransversality}
\subsection{Combinatorial data}\label{sec:CD}
Let $\Delta$ be a reduced root system, and $\Lambda$ the
associated weight lattice. Suppose given an affine hyperspace $K$ of $\Lambda \otimes
\mathbb{R}$, and denote by $\mathscr{D}K$ the parallel linear hyperspace. The goal of
this section is to define some combinatorial data, depending on (the linearization of) our
affine hyperspace $K$, which develop consistently as
we restrict our attention to smaller and smaller root subsystems. 

Fix for reference a system $\Pi$ of simple roots. As usual, we
identify $\Pi$ with the vertex set of the Dynkin diagram of $\Delta$.
Let us call a simple root $\omega$ \emph{extremal} if the subset $\Pi \setminus \{\omega\}$
corresponds to the vertex set of a connected subgraph of the Dynkin diagram of $\Delta$.
Given an extremal vertex $\omega$ and an element $w \in W(\Delta)$ such that $\Pi \setminus \{\omega\} \not \subset w^{-1} \cdot \mathscr{D}K$,
there exists an $\alpha \in \Pi \setminus \{\omega\}$
not parallel to $w^{-1} \cdot K$. If we insist that the combinatorial path
$$
\alpha = \beta_{1} \sim \cdots \sim \beta_{q} = \omega
$$
connecting $\alpha$ to $\omega$ in the Dynkin diagram of $\Delta$ passes entirely
through vertices $\beta_{2}, \ldots, \beta_{q - 1} \in \Pi$ lying parallel to $w^{-1} \cdot
K$, then $\alpha$ is nearly unique. Fix an extremal vertex $\omega$, and denote by $\alpha_{w} \in \Pi \setminus \{\omega\}$ one such 
root.

\begin{lemma}\label{lem:alpha-omega}
Let $\Delta$ be an irreducible, reduced root system of simply laced type and of higher
rank; $K$, $\mathscr{D}K$, $\Pi$, and $\omega$ as above. Then there exists an element $w
\in W$ such that
\begin{itemize}
\item[(i)] $w \cdot \omega$ is not parallel to $K$;
\item[(ii)] $\Pi \setminus \{\omega\} \not \subset w^{-1} \cdot \mathscr{D}K$;
and
\item[(iii)] $w \cdot \alpha_{w}$ is not orthogonal to $K$.
\end{itemize}
\end{lemma}
\begin{proof}
It follows from irreducibility of $\Delta$ that
the Weyl group cannot stabilize the union of a proper subspace and its
orthogonal complement. In particular, there exists $w \in W$ such that $w \cdot \omega$
lies outside $\mathscr{D}K \cup \mathscr{D}K^{\perp}$. For a nonzero root $\alpha \in
\Delta$, let $w_{\alpha} \in W$ be
the corresponding simple reflection. Suppose first that $w$ does not
satisfy condition~(ii). Then $w_{\omega}w_{\omega + \beta_{q - 1}}w$ verifies the
conclusion of the lemma.

Therefore, we suppose that $w$ already verifies condition~(ii).
If $\alpha_{w}$ is orthogonal to $\omega$, then we have that $w \cdot \alpha_{w}$ also
lies outside of $\mathscr{D}K \cup \mathscr{D}K^{\perp}$. We are left with the case
where $\alpha_{w} \sim \omega$, i.e.\ $\langle \alpha_{w}, \omega \rangle =
-1$ (using the normalization that $\langle \alpha, \alpha \rangle = 2$ for all roots
$\alpha
\neq 0$). In the case where $w$ fails to verify condition~(iii), then
$\alpha_{w_{\omega}w} = \alpha_{w}$, and $w_{\omega}w$ verifies the
conclusion of the lemma.
\end{proof}

Let $w$ be as in Lemma~\ref{lem:alpha-omega}. Upon replacing $\Pi$ by $w \cdot \Pi$, we may assume $w =
1$. That is, 
\begin{itemize}
\item[(CD1)] $\omega$ is not parallel to $K$;
\item[(CD2)] $\Pi \setminus \{\omega\} \not \subset \mathscr{D}K$;
\item[(CD3)] $\alpha = \alpha_{w}$ is not orthogonal to $K$; and
\item[(CD4)] the combinatorial path
$\alpha = \beta_{1} \sim \cdots \sim \beta_{q} = \omega$ 
connecting $\alpha$ to $\omega$ in the Dynkin diagram of $\Delta$ passes entirely
through vertices $\beta_{2}, \ldots, \beta_{q - 1} \in \Pi$ lying parallel to $K$.
\end{itemize}

Define $\Delta' = \Delta \cap
\operatorname{Span}(\Pi \setminus \{\omega\})$. By
construction, $\Delta'$ is again irreducible and of simply laced type.
By definition, $\alpha$ is not parallel to $K$. It follows that $K' = K \cap
\operatorname{Span} \Delta'$ has codimension one in $\operatorname{Span} \Delta'$.
Provided $\Delta'$ has higher rank, we can choose a set $\Pi'$ of simple roots for $\Delta'$
and a simple root $\omega' \in \Pi'$ such
that properties~(CD1)--(CD4) hold with $\omega$, $K$,
$\Pi$, and $\alpha$ replaced by $\omega'$, $K'$, etc. If $\operatorname{rank} \Delta' =
1$, then we simply set $\alpha' = \omega'$ to be a nonzero root, and $\Pi' =
\{\alpha'\}$.

\subsection{Weight spaces for simply laced groups}\label{sec:weightspacesforSLgroups}
Let $G$ be a connected, reductive complex algebraic group.
We fix a maximal torus $T$
and a system $\Pi$ of simple roots for $\Delta(T, G)$. Denote by $B$ the corresponding
Borel subgroup of $G$.

In this section, we will also assume that $G$ is a simply laced group, i.e.\ a
group of type $A$, $D$, or $E$ (modulo its unipotent radical). 
We begin essentially by formulating statements analogous to those in \S\ref{sec:rank1facts}, but in the higher rank
situation. We recall the
fact that the unipotent radical $B_{u}$ of $B$ is isomorphic as an algebraic variety
to affine space of dimension $\abs{\Delta^{+}}$. Explicitly, the map
$$\Vect{t} = (t_{\beta})_{\beta \in \Delta^{+}} \mapsto \prod_{\beta \in \Delta^{+}}
u_{\beta}(t_{\beta})$$
is an isomorphism of complex algebraic varieties, where the product is taken in any fixed
order on $\Delta^{+}$.

Fix an numbering $\langle \beta_{1}, \beta_{2}, \ldots \rangle$ of
$\Delta^{+}$ compatible with the poset structure: that is to say, with the
property that $\beta_{i} < \beta_{j}$ implies $i > j$. Fix an irreducible representation of $G$ on a finite-dimensional complex vector space $E$. Write
$u_{i} = u_{\beta_{i}}$ and $t_{i} = t_{\beta_{i}}$. To an
arbitrary but fixed nonzero vector $\Vect{e}
\in E$, we can associate an endomorphism
$\prod_{i = 1}^{\abs{\Delta^{+}}} (u_{i}(t_{i}) - I)^{z_{i}}$,
where the numbers $z_{i}$ are defined inductively as follows:
The vector-valued function
$\prod_{i = 1}^{j} (u_{i}(t_{i}) - I)^{z_{i}} \cdot \Vect{e}$
is not zero identically in $(t_{1}, \ldots, t_{j})$, but this vector itself lies
in the kernel of $u_{j}(t_{j}) - I$ for all $t_{j} \in \mathbb{C}$. If $\Vect{e}$ is a
weight vector, then so is $\prod_{i = 1}^{\abs{\Delta^{+}}} (u_{i}(t_{i}) - I)^{z_{i}}
\cdot \Vect{e}$.
\begin{lemma}\label{lem:raisingoperators}
The vector
$\prod_{i = 1}^{\abs{\Delta^{+}}} (u_{i}(t_{i}) - I)^{z_{i}} \cdot \Vect{e}$
is independent of $\Vect{t} = (t_{i})_{i = 1, \ldots, \abs{\Delta^{+}}}$, up to scalar. If
$E$ is
irreducible as a complex representation of $G$, then the line in this direction is
the unique $B$-invariant line in $E$.
\end{lemma}
\begin{proof}
The first assertion follows from the analogous phenomenon in rank one. For the second,
we need only show invariance under $B$.

Taking $\Vect{e} \in E^{\mu}$, we see that for almost all $\Vect{t}$, we have that $\mu
+ \sum_{i = j}^{\abs{\Delta^{+}}} z_{i}\beta_{i}$ is a lowest weight in
$\pounds_{E}\left(\prod_{i = j}^{\abs{\Delta^{+}}}(u_{i}(t_{i}) - I)^{z_{i}} \cdot
\Vect{e}\right)$, and likewise when $\Vect{e}$ is replaced by any element of $B \cdot
\Vect{e}$. Therefore
$$\pounds_{E}\left(b \cdot \prod_{i = j}^{\abs{\Delta^{+}}}(u_{i}(t_{i}) - I)^{z_{i}} \cdot
\Vect{e}\right) = \pounds_{E}\left(\prod_{i = j}^{\abs{\Delta^{+}}}(\mbox{}^{b}(u_{i}(t_{i})) -
I)^{z_{i}} \cdot
\Vect{e}\right)$$
for any $b \in B$. The result now follows from the compatibility of the numbering with
the poset structure on $\Delta^{+}$, together with the fact that $\prod_{i = 1}^{\abs{\Delta^{+}}} (u_{i}(t_{i}) - I)^{z_{i}} \cdot \Vect{e}$
is independent of $\Vect{t}$, up to scalar.
\end{proof}

Given a set of simple roots $S \subset \Pi$, we can identify $S$
with a subgraph of the Dynkin diagram of $G$. We shall be
interested in those $S$ for which $S$ and $S' = \Pi \setminus S$
correspond to connected (nonempty) subgraphs.

Denote by $G_{S}$ the group generated by the root subgroups $U_{\alpha}$, $U_{-\alpha}$
for $\alpha \in S$. Set $\Delta_{S} = \Delta \cap \operatorname{Span}
S$ and $\Delta^{+}_{S} =
\Delta^{+} \cap \Delta_{S}$. By symmetry, we can adopt similar notations for the situation
where $S$ and $S'$ exchange r\^oles.

Before turning to the next lemma, it is useful to recall the commutation relations for
unipotent subgroups $U_{\alpha'}$, $\alpha' \in \Delta^{+}$. If $\alpha$, $\beta \in
\Delta^{+}$ are such that
$\alpha + \beta$ is not a root, then $U_{\alpha}$ and $U_{\beta}$ commute. Otherwise,
$u_{\alpha}(t)$ and $u_{\beta}(s)$ commute to $u_{\alpha + \beta}(cst)$
for some nonzero constant $c \in \mathbb{C}$ depending only on $\alpha$ and $\beta$. Though not
strictly necessary for our purposes, we recall that, upon normalizing the unipotent
subgroups $U_{\alpha'}$, $\alpha' \in \Delta$, we can take $c = c(\alpha, \beta) = \pm
1$ for all positive roots $\alpha$ and $\beta$.

Suppose that $E' < E$ is a $G_{S}$-submodule of E. If we assume moreover that $E'$ is irreducible as a $G_{S}$-submodule, then there is a unique line in $E'$ which
is stable under all $U_{\beta}$, $\beta \in \Delta^{+}_{S}$. Let $\Vect{e}' \in E'$
denote a nonzero vector in the invariant direction.
Since $S'$ is the vertex set of a connected subgraph,
we can apply the preceding considerations to produce numbers
$y_{\alpha}$, $\alpha \in \Delta_{S'}^{+}$, such that the vector
$\prod_{\alpha \in \Delta_{S'}^{+}} (u_{\alpha}(t_{\alpha}) - I)^{y_{\alpha}} \cdot \Vect{e}'$
(in a fixed ordering on $\Delta_{S'}^{+}$) is independent of $\Vect{t} = (t_{\alpha})_{\alpha \in \Delta_{S'}^{+}}$, up to
scalar.
Let $\Vect{e}''$ be a nonzero vector lying in this direction. Applying
Lemma~\ref{lem:raisingoperators} to the group $G_{S'}$, we obtain that
$\Vect{e}''$ is invariant under those subgroups $U_{\alpha}$ with $\alpha \in
\Delta_{S'}^{+}$.

\begin{lemma}\label{lem:hwvLiesInLongChain}
Let $S$ be the complement in $\Pi$ of an extremal $\omega$.
Suppose $E'$, $\Vect{e}'$, $\Vect{e}''$ are as above. Then
$\Vect{e}''$ is invariant under $U_{\beta}$ for all $\beta \in
\Delta^{+} \setminus \Delta_{S}$.
\end{lemma}
\begin{proof}
We know that $\Vect{e}''$ is a weight vector, say $\Vect{e}'' \in
E^{\mu}$. Then $\mu$ is the highest weight in the set $\pounds_{E}
(S_{\omega} \cdot \Vect{e}') = \pounds_{E} (U_{\omega} \cdot
\Vect{e}') \cup \pounds_{E} (U_{-\omega} \cdot
\Vect{e}')$. Fix a simple root $\alpha \in S$. We compute
\begin{equation}\label{eqn:prod_j'}
\begin{split}
u_{\omega}(t_{\omega}) \cdot \Vect{e}'
 &= u_{\omega}(t_{\omega}) u_{\alpha}(t)
\cdot
\Vect{e}'\\
&= u_{\alpha + \omega}(t_{\omega} t) u_{\alpha}(t)
u_{\omega}(t_{\omega}) \cdot \Vect{e}'.
\end{split}
\end{equation}

Now consider a positive root $\beta \notin \Delta_{S}$ such that
$U_{\beta}$ stabilizes $\Vect{e}''$. If $\alpha + \beta$ is a
root, then $u_{\alpha}(-t) u_{\beta}(-s/t) u_{\alpha}(t) =
u_{\alpha + \beta}(s)u_{\beta}(-s/t)$. By hypothesis on $\beta$,
it follows that $u_{\alpha + \beta}(s) \cdot \Vect{e}'' =
u_{\alpha}(t) u_{\beta}(-s/t)u_{\alpha}(t) \cdot \Vect{e}''$.
Combining with equation~\eqref{eqn:prod_j'}, we have
\begin{equation*}
\begin{split}
\pounds_{E} (u_{\alpha + \beta}(s) \cdot \Vect{e}'') &=
\pounds_{E}(u_{\alpha}(-t) u_{\beta}(-s/t) u_{\alpha}(t) \cdot
\Vect{e}'')\\
&\subset \pounds_{E}(U_{\alpha} U_{\beta} u_{\omega}(t_{\omega})
u_{\alpha}(t) \cdot \Vect{e}'')\\
&= \pounds_{E}(U_{\alpha} U_{\beta} u_{\alpha + \omega}(t_{\omega}
t) u_{\alpha}(t)
u_{\omega}(t_{\omega}) \cdot \Vect{e}'')\\
&\subset \pounds_{E}(U_{\alpha} U_{\beta} u_{\alpha +
\omega}(t_{\omega} t) u_{\alpha}(t)
u_{\omega}(t_{\omega}) \cdot \Vect{e}')\\
&= \pounds_{E}(U_{\alpha} U_{\beta} u_{\omega}(t_{\omega}) \cdot
\Vect{e}')
\end{split}
\end{equation*}
for almost every $t_{\omega} \in \mathbb{C}$.

Assume contrary to fact that $\mu + \alpha + \beta$ lies in
$\pounds_{E} (U_{\alpha} U_{\beta} U_{\omega} \cdot \Vect{e}')$.
Then either $\mu + \beta$ lies in $\pounds_{E} (U_{\beta}
U_{\omega} \cdot \Vect{e}')$ or $\mu + \alpha - \ell \beta$ lies
in $\pounds_{E}(U_{\omega} \cdot \Vect{e}')$ for some $\ell > 0$.
The first option is impossible because $U_{\beta}$ fixes
$\Vect{e}''$. We proceed to show the second is impossible by
computing the scalar products of $\alpha - \ell \beta$ against the
fundamental weights.

Since $\mu$ is the highest weight of $\pounds_{E}(U_{\omega} \cdot
\Vect{e}')$, we see that $\ell \beta - \alpha$ must be a positive
integer multiple of $\omega$. In particular, $\ell \beta - \alpha$
must be orthogonal to $\varpi_{\alpha'}$ for $\alpha' \in S$.
Taking $\alpha = \alpha'$, we obtain
$$\ell \langle \beta, \varpi_{\alpha} \rangle= 1,$$
whence $\ell = 1$. For the $\alpha' \in S$ different from
$\alpha$, the orthogonality condition simply implies that $\langle
\beta, \varpi_{\alpha'} \rangle = 0$. It follows that $\beta =
\alpha + \ell' \omega$ for some positive integer $\ell'$. But
since $\beta \in \Delta^{+}$ and $\omega$ is extremal, we must
have $\ell' = 1$. But then it cannot happen that $\alpha + \beta$
is a root.

On the other hand, suppose that $\beta$ is as above, but $\alpha \in \Pi \setminus S = \{\omega\}$.
Then we have that $U_{\alpha}U_{\beta}U_{\alpha}$ fixes $\Vect{e}''$. But
$u_{\alpha}(t)u_{\beta}(s)u_{\alpha}(-t) = u_{\alpha + \beta}(st)u_{\beta}(s)$.
Combining these observations yields
$$\pounds_{E}(U_{\alpha + \beta} \cdot \Vect{e}'') = \pounds_{E}(U_{\alpha}U_{\beta}U_{\alpha} \cdot
\Vect{e}'') = \pounds_{E}(\Vect{e}'') = \{\mu\}.$$

So far, we have shown that if a positive root $\beta$ is such that $U_{\beta}$ leaves
$\Vect{e}''$ invariant, then the same is true for all roots for the form $\alpha +
\beta$, $\alpha \in \Pi$. By construction, $\Vect{e}''$
is invariant under $U_{\omega}$.
The result now follows.
\end{proof}

\begin{definition}
Let $k$ be a field. Given a positive integer $\ell$,
we shall say that a $k$-rational, absolutely irreducible representation of a $k$-split,
reductive algebraic $k$-group $G$ on a finite-dimensional vector space $E$ is
\emph{$\ell$-large}
if the set of its weights contains a root string of length $\ell$;
otherwise, we shall say that $E$ is \emph{$\ell$-small}.
\end{definition}
\noindent We remark that $E$ is $\ell$-large if and only
if there exist a weight $\mu$ of $E$ and a root $\alpha$ of $G$ with
$$\frac{2 \langle \mu, \alpha \rangle}{\langle \alpha, \alpha \rangle} \geqslant \ell - 1.$$
Using the normalization that $\langle \alpha, \alpha \rangle = 2$ for all roots $\alpha
\neq 0$, the previous inequality reads simply $\langle \mu, \alpha \rangle \geqslant \ell -
1$.

On the other hand, if $\mu$ is a weight of an $\ell$-small representation of $G$, then
we can write $\mu = \sum_{\alpha \in \Pi} n_{\alpha} \varpi_{\alpha}$, with
$\abs{n_{\alpha}} < \ell - 1$ for all $\alpha \in \Pi$.
It follows that there exists a constant $\kappa$, depending only on the type of the root
system $\Delta$, such that if $\mu$ is a weight of an arbitrary $\ell$-small
representation of a reductive group with absolute root system $\Delta$, then $\norm{\mu}
< \kappa(\ell - 1)$. Let $\kappa_{r}$ be a positive number no greater than the
minimum value of $\kappa(\Delta)$ as $\Delta$ ranges over all root systems of type $ADE$
and rank $r$.

Given a representation of $G$ on a finite-dimensional complex vector space $E$ and a set
$S$ of simple roots, we can decompose $E$ as an internal
direct sum of irreducible $G_{S}$-modules.
If $S$ is the vertex set of a connected subgraph, then it makes sense to speak of
$\ell$-small or $\ell$-large $G_{S}$-submodules of $E$.

Given an extremal vertex $\omega$, consider the set $D_{\omega} =
\{\beta \in \Delta^{+} \mid \langle \beta, \varpi_{\omega} \rangle
= 1\}$. We are interested in a certain convex cone in $\Lambda
\otimes \mathbb{R}$, where $\Lambda$ denotes the root lattice
associated to $\Delta$. Namely, define $C_{\omega}$ to be the cone
whose base is the convex hull of $D_{\omega} \setminus
\{\omega\}$, and whose vertex is $0$. Choose positive constants
$\theta_{i}$ such that for any root system of type $ADE$ and rank
$i$, we have that the angle between $\omega$ and $C_{\omega}$ is
no less than $\theta_{i}$ for any extremal root $\omega$.

Let $r$ denote the rank of $G$.
\begin{lemma}\label{lem:l-2l'}
Let $S$ be the complement in $\Pi$ of an extremal vertex $\omega$,
as above. Let $E$ be an irreducible complex representation of
$G$ whose set of weights contains a root string of length $\ell$.
Suppose that $\Vect{e}$ is a nonzero element of an $\ell'$-small
submodule $E' < E$ for the action of $G_{S}$. Then
$\pounds_{E}(S_{\omega} \cdot \Vect{e})$ contains a root string of
length exceeding $\ell - 4\kappa_{r - 1}(\ell' - 1)\tan \theta_{r}$.
\end{lemma}
\begin{proof}
We may assume without loss of generality that $\Vect{e}$ is a
weight vector, say $\Vect{e} \in E^{\nu}$. Define corresponding
vectors $\Vect{e}'$ and $\Vect{e}''$ as in
Lemma~\ref{lem:hwvLiesInLongChain}, and set $\nu'$ and $\mu$ to be
their respective weights. In this case, $\mu - \nu'$ is a
nonnegative multiple of $\omega$.

Denote by $\lambda$ the highest weight of the representation $E$.
By irreducibility we have $\lambda \in \pounds_{E}(B \cdot
\Vect{e}'')$. It follows from Lemma~\ref{lem:hwvLiesInLongChain}
that $\lambda - \mu$ lies in the span of $S$. Therefore, we have
the identity $\langle \mu, \varpi_{\omega} \rangle = \langle
\lambda, \varpi_{\omega} \rangle$.

Let $\mu^{+}$ be the highest weight in the string
$\pounds_{E}(U_{\omega} \cdot \Vect{e})$. Consider the trapezoid
with corners $\mu$, $\nu$, $\nu'$, and $\mu^{+}$. Denote by
$U_{S}$ the subgroup generated by the root subgroups associated to
roots in $S$. In other words, $U_{S}$ is the product of the root
subgroups associated to roots in $\Delta^{+}_{S}$. By definition,
$\mu$ lies in $\pounds_{E}(U_{\omega}U_{S} \cdot \Vect{e})$. The
commutation relations tell us in particular that
$\pounds_{E}(U_{\omega}U_{S} \cdot \Vect{e}) \subset
\pounds_{E}((\prod_{\beta \in D} U_{\beta}) U_{S}U_{\omega} \cdot
\Vect{e})$, where $D = D_{\omega} \setminus \{\omega\}$. By definition, we have that $\mu - \mu^{+}$ forms with
$\omega$ an angle no smaller than $\theta_{r}$.

Now since $E'$ is $\ell'$-small, it follows that $\norm{\nu -
\nu'} < 2\kappa_{r - 1}(\ell' - 1)$. Since $\langle \nu, \varpi_{\omega}
\rangle = \langle \nu', \varpi_{\omega} \rangle$, we conclude that
\begin{equation*}
\begin{split}
\langle \mu^{+}, \varpi_{\omega} \rangle &\geqslant \langle \mu,
\varpi_{\omega} \rangle - 2\kappa_{r - 1}(\ell' - 1)\tan \theta_{r}\\
&= \langle \lambda,
\varpi_{\omega} \rangle - 2\kappa_{r - 1}(\ell' - 1)\tan \theta_{r}
\end{split}
\end{equation*}
The same
argument used for the opposite system $-\Delta^{+}$ tells us that
the lowest weight $\mu^{-}$ of $\pounds_{E}(U_{-\omega}\cdot
\Vect{e})$ satisfies $\langle \mu^{-}, \varpi_{\omega} \rangle
\leqslant \langle w_{0} \cdot \lambda, \varpi_{\omega} \rangle +
2\kappa_{r - 1}(\ell' - 1)\tan \theta_{r}$. The result follows.
\end{proof}

For a subset $S \subset \Pi$ which is the vertex set of a connected subgraph of the
Dynkin diagram, we denote by $E_{S, \textup{ $\ell$-small}}$
the sum of all $\ell$-small $G_{S}$-submodules of $E$. In this case, we can canonically
identify $E_{S, \textup{ $\ell$-small}}$ as a $G_{S}$-submodule of $E$.
Set $\tau_{s} = 2 \kappa_{s - 1} \tan \theta_{s}$.

\begin{corollary}\label{cor:longOmegaChain}
Let $S$ and $\omega$ be as in Lemma~\ref{lem:l-2l'}. If $E$ is an
$(\ell' + 2 \ell \tau_{r})$-large representation of $G$, then
$$E_{S, \textup{ $\ell$-small}} \cap E_{\{\omega\}, \textup{ $\ell'$-small}} = 0.$$
\end{corollary}

Let us conclude with the observation that $\mu \in \pounds_{E}
(E_{\textup{$\ell$-small}})$ implies that $\langle \mu, \beta
\rangle < \ell$ for all $\beta \in \Delta_{S}$.

\subsection{Exterior powers}\label{sec:ExteriorPowers}
We begin with a very general lemma. Let $G$ be a connected
algebraic group and $E$ a $G$-module with basis $\mathscr{B}_{E}$.
Given a finite subset $I \subset \mathscr{B}_{E}$, we define
$\Vect{e}_{I} = \bigwedge_{\Vect{\upsilon} \in I}
\Vect{\upsilon}$. (We acknowledge yet ignore the subtlety that
$\Vect{e}_{I}$ is only defined modulo sign.) Then we have that the
set $\{\Vect{e}_{I} \mid I \subset \mathscr{B}_{E}, \abs{I} = m\}$
forms a basis for the $m$-fold exterior power
$\textstyle{\bigwedge^{\!m}} E$ of $E$.
\begin{lemma}\label{lem:indexsetdecomposition}
Let $E$ and $F$ be modules for an connected algebraic $k$-group
$G$, with bases $\mathscr{B}_{E}$ and $\mathscr{B}_{F}$,
respectively. Take finite subsets $I \subset \mathscr{B}_{E}$, $J \subset
\mathscr{B}_{F}$, $\Vect{e} \in \textstyle{\bigwedge^{\!\abs{I}}} E$, and
$\Vect{f} \in \textstyle{\bigwedge^{\!\abs{J}}} F$. Write $m = \abs{I} + \abs{J}$. Let $H$ be an irreducible
subvariety of $G$ defined over $k$.

Suppose that, for some $h \in H$, we have that $h \cdot \Vect{e}$
(resp.\ $h \cdot \Vect{f}$) has nonzero component in the
$\Vect{e}_{I}$-direction (resp.\ the $\Vect{e}_{J}$-direction).
Then for an open, dense set of $h \in H$, we have the following:
\begin{itemize}
\item[(i)] The vector
$h \cdot (\Vect{e} \wedge \Vect{f}) \in \textstyle{\bigwedge^{\!m}} (E \oplus F)$ has
nonzero component in the
$$\Vect{e}_{I \cup J} = \Vect{e}_{I} \wedge
\Vect{e}_{J}$$ direction with respect to the basis $\textstyle{\bigwedge^{\!m}} (\mathscr{B}_{E} 
\cup \mathscr{B}_{F})$; and \item[(ii)] likewise for $h \cdot
(\Vect{e} \wedge \Vect{e}_{J})$.
\end{itemize}
\end{lemma}
\begin{proof}
The set of $g$ such that $g \cdot \Vect{e}_{J}$ has nonzero component in the
$\Vect{e}_{J}$-direction is open and contains the identity. Thus, the second assertion is
a special case of the first.

We now treat assertion~(i). Consider the set $H_{\Vect{e}}$ (resp.\ $H_{\Vect{f}}$) of $g
\in H$ such that $g \cdot \Vect{e}$ (resp.\ $g \cdot \Vect{f}$) has nonzero component in
the $\Vect{e}_{I}$-direction (resp.\ the $\Vect{e}_{J}$-direction). By hypothesis, $H_{\Vect{e}}$
and $H_{\Vect{f}}$ are nonempty. But they are also open subsets of $H$. By irreducibility,
their intersection is nonempty, again open, and therefore dense. If $h$ lies in the
intersection $H_{\Vect{e}} \cap H_{\Vect{f}}$, then $h$ satisfies the conclusion of assertion~(i).
\end{proof}

Suppose now that $G$ is a reductive complex algebraic group, and that $E$ is a
finite-dimensional complex vector space on which $G$ acts.
Suppose furthermore that our basis $\mathscr{B}$ of $E$ consists of weight
vectors for the action of a fixed maximal torus $T$, say $\Vect{e}^{\mu}_{i}$, $1
\leqslant i \leqslant \dim E^{\mu}$. If $I$ is a subset of
$\{(\mu, i) \mid 1 \leqslant i \leqslant \dim E^{\mu}\}$, then
define
$$\Vect{e}_{I} = \!\!\!\!\bigwedge_{(\mu, i) \in I}\!\!\!\!\Vect{e}^{\mu}_{i}.$$
If $K$ is a set of weights of $E$ (or more generally, a subset of
$\Lambda \otimes \mathbb{R}$) and $m = \sum_{\mu \in K} \dim
E^{\mu}$, then let us denote by $\Vect{e}_{K}$ the distinguished
basis element $$\Vect{e}_{\{(\mu, i) \mid \mu \in K, 1 \leqslant i
\leqslant \dim E^{\mu}\}}$$ of the $m$-fold exterior power
$\textstyle{\bigwedge^{\!m}} E$. Also, let us write $[K]$ for the
weight $\sum_{\upsilon \in K} \pounds \upsilon$. In this notation,
$\Vect{e}_{K}$ is a weight vector of weight $[K]$ in the
representation $\textstyle{\bigwedge^{\!m}} E$.

Given a nonzero root $\beta \in
\Delta(T, G)$, we form the basis $\{\Vect{e}^{\mu}_{i}\}$ as
usual. Likewise, when $I$ is a subset of $\{(\mu, i) \mid 1 \leqslant i \leqslant \dim
E^{\mu}\}$ we write $\Vect{e}_{I}$ for the corresponding monomial basis element for
the exterior power. As we
have seen, given a $G$-module $E$, the index set $\{(\mu, i) \mid 1 \leqslant i \leqslant \dim E^{\mu}\}$
carries some additional structures, depending on the choice of $\beta$: Namely, the relation of $\beta$-linkedness and the notion of $\beta$-chains.
If $C$
is a $\beta$-chain of $E$, then let us denote by $\ell(C)$ the length (i.e.\ cardinality)
of such a chain.

We remark that, while the index set $\{(\mu, i) \mid 1 \leqslant i
\leqslant \dim
E^{\mu}\}$ is fixed, the very basis which it indexes depends on the choice of a nonzero
root. We will soon have occasion to compare these constructions for different $\beta$.
We will resolve eventual ambiguities by explicitly appending the root to the notation. For
example,
we shall speak of the ``component in the $\Vect{e}_{I}(\beta)$-direction''
to specify that component with respect to the basis of monomials in the
vectors $\Vect{e}^{\mu}_{i} = \Vect{e}^{\mu}_{i}(\beta)$ associated to $\beta$.

The
substance of the next lemma is that a certain property is invariant under
changing the basis corresponding to one nonzero root to the basis corresponding to
another.
\begin{definition}
Let us say that a subset $I \subset \{(\mu, i) \mid 1 \leqslant i \leqslant \dim
E^{\mu}\}$ is \emph{$(\beta, \ell)$-sparse} if $\abs{I \cap C} < \min \{\ell, \ell(C)\}$
for every maximal $\beta$-chain $C$ of $E$.
\end{definition}

\begin{lemma}\label{lem:beta<->beta'}
Let $G$ be a simply laced, connected, reductive complex algebraic group whose Dynkin
diagram
is a connected graph. Let $\beta$, $\beta'$ be nonzero roots. Let $I$ be a $(\beta,
\ell)$-sparse subset of
$\{(\mu, i) \mid 1 \leqslant i \leqslant \dim E^{\mu}\}$.
Then $g \cdot \Vect{e}_{K}$ has a
nonzero component in the $\Vect{e}_{I}(\beta)$-direction for some/almost every $g$ if and
only if there exists a $(\beta',
\ell)$-sparse subset $J \subset
\{(\mu, i) \mid 1 \leqslant i \leqslant \dim E^{\mu}\}$ such that $g \cdot \Vect{e}_{K}$ has a nonzero component in
the $\Vect{e}_{J}(\beta')$-direction for some/almost every $g$.
\end{lemma}
\begin{proof}
Let $\dot{w}$ be a representative in $G$ of an element $w \in W$ sending $\beta$ to
$\beta'$. Since $\dot{w}$ conjugates $U_{\beta}$ to $U_{\beta'}$, it follows that
$\dot{w} \cdot \Vect{e}^{\mu}_{i}(\beta) = \Vect{e}^{w \cdot \mu}_{i'}(\beta')$,
up
to scalar multiple, for some
$i'$, $1 \leqslant i' \leqslant \dim E^{w \cdot \mu} = \dim E^{\mu}$. Thus $w$ sends $\beta$-chains to $\beta'$-chains of the same
length. The result follows.
\end{proof}

We begin to specialize the situation. We take for $K$ an
affine hyperspace of $\Lambda\otimes \mathbb{R}$, and denote by
$\mathscr{D}K$ the parallel linear hyperspace. For a fixed nonzero root $\alpha
\in \Delta$, we denote by $\{\Vect{e}^{\mu}_{i}\}$ the basis
constructed in \S\ref{sec:rank1facts}. When a construction involving basis vectors
depends on the choice of $\alpha$, we shall use the notational convention of Lemma~\ref{lem:beta<->beta'} to
make the choice explicit.

Set $\ell_{s} = 2^{s}\prod_{i = 1}^{s} \lceil 1 + \tau_{i} \rceil$. Set $S = \Pi
\setminus
\{\omega\}$. We denote by $I_{\textup{$\ell$-small}}$
the set of indices $(\mu, i)$ with $\Vect{e}^{\mu}_{i} \in E_{S, \textup{
$\ell$-small}}$. Since $E_{S, \textup{
$\ell$-small}}$ is an $S_{\alpha}$-submodule of $E$, it follows
that the set of $\Vect{e}^{\mu}_{i}$ with $(\mu, i) \in
I_{\textup{$\ell$-small}}$ forms a basis for it. Denote by $I_{K}$ the set of indices $(\mu, i)$ with $\mu \in K$.

\begin{proposition}\label{prop:e_kVg.e_k}
Suppose $G$ is a simply laced group of rank $r$; $K$ an affine hyperspace of $\Lambda \otimes
\mathbb{R}$; and $\Pi$ and $\alpha = \beta_{1} \sim \cdots \sim \beta_{q} = \omega$
chosen as in properties~\textup{(CD)} of~\S\ref{sec:CD}.

If $E$ is an $\ell_{r}$-large $G$-module, then there exists an $(\alpha, 2\ell_{r -
1})$-sparse set $I$ such that the $\Vect{e}_{I}$-component of $g \cdot \Vect{e}_{K}$
is nonzero for some (equiv., almost every) $g \in G$.
\end{proposition}
\begin{proof}
Induction on $r$. The rank one case follows from
Corollary~\ref{cor:rank1exteriorpower}, for instance. Thus, assume $r > 1$.

For convenience, we shall denote by $w_{j}$ the simple reflection through the simple root
$\beta_{j}$ of the combinatorial path of property~(CD4) of~\S\ref{sec:CD}. We observe that the
subsystem
$\Delta'$ of~\S\ref{sec:CD} coincides with $\Delta_{S}$. Set $w = w_{2} \cdots w_{q - 1} \in W_{S}$. Denote by $K''$ the affine hyperspace
$(w ^{-1}\cdot \alpha^{\perp} \cap K) \oplus \mathbb{R}\omega$. Let $\alpha'$ and $\omega'$ be
as in~\S\ref{sec:CD}.

Let $E'$ be an irreducible, $\ell_{r - 1}$-large $G_{S}$-submodule of $E$. By the induction
hypothesis, there exists an $(\alpha', 2\ell_{r - 2})$-sparse set $I'$ (to
use the notation of Lemma~\ref{lem:beta<->beta'}) such
that $g \cdot \Vect{e}_{K \cap \pounds_{E} E'}$ has nonzero component in the $\Vect{e}_{I'}(\alpha')$-direction
for almost every $g \in G_{S}$.
Let $w' \in W_{S}$ transport $w^{-1} \cdot \alpha$ to $\alpha'$.
Since $w^{-1} \cdot \alpha$ is not parallel to $K''$, it follows from Corollary~\ref{cor:rank1exteriorpower}
that we can take $I'$ to satisfy $\pounds
I' \cap w' \cdot K'' = \varnothing$. Now Lemma~\ref{lem:beta<->beta'} tells us that there
exists a $(w^{-1} \cdot \alpha, 2\ell_{r - 2})$-sparse set $J'$ satisfying $\pounds
J' \cap K'' = \varnothing$, such
that $g \cdot \Vect{e}_{K \cap \pounds_{E} E'} = g \cdot \Vect{e}_{w^{-1} \cdot K \cap
\pounds_{E}
E'}$ has nonzero component in the $\Vect{e}_{J'}(w^{-1} \cdot \alpha)$-direction
for almost every $g \in G_{S}$.

Collecting the preceding facts using Lemma~\ref{lem:indexsetdecomposition},
we see that there exists a set $I$ satisfying
\begin{itemize}
\item[(I1)] $I \cap I_{\textup{$\ell_{r - 1}$-small}}$ coincides with $I_{w ^{-1}\cdot K} \cap I_{\textup{$\ell_{r -
1}$-small}} = I_{K} \cap I_{\textup{$\ell_{r - 1}$-small}}$
\item[(I2)] $I \cap I_{K''} \subset I_{w ^{-1}\cdot K} \cap I_{\textup{$\ell_{r - 1}$-small}} = I_{K} \cap I_{\textup{$\ell_{r - 1}$-small}}$
\item[(I3)] $I \setminus I_{\textup{$\ell_{r - 1}$-small}}$ is $(w^{-1} \cdot \alpha,
2\ell_{r - 2})$-sparse
\end{itemize}
such that $g \cdot \Vect{e}_{w^{-1} \cdot K} = g \cdot \Vect{e}_{K}$ has nonzero
component in the $\Vect{e}_{I}(w ^{-1}\cdot \alpha)$-direction for almost every $g \in G_{S}$.

Corollary~\ref{cor:longOmegaChain}
tells us that $\pounds_{E}(S_{\omega} \cdot \Vect{e}^{\mu}_{i}) =
\pounds_{E}(U_{\omega} \cdot \Vect{e}^{\mu}_{i}) \cup
\pounds_{E}(U_{-\omega} \cdot \Vect{e}^{\mu}_{i})$ is an
$\omega$-string of length at least $2 \ell_{r - 1}$ for each
index $(\mu, i) \in I_{\textup{$\ell_{r - 1}$-small}}$. In fact, combining
Corollary~\ref{cor:longOmegaChain} with Lemma~\ref{lem:unipotentorbits}
tells us even more: For $\mu \in
\pounds I_{\textup{$\ell_{r - 1}$-small}}$ \emph{fixed} and any integer $m$ satisfying
$\abs{\langle \mu, \omega \rangle + 2m} < 2 \ell_{r - 1}$,
the map
$$
\operatorname{Span} \langle\Vect{e}^{\mu}_{i} \rangle_{(\mu, i)
\in I_{\textup{$\ell_{r - 1}$-small}}} \to E^{\mu + m\omega}
$$
is injective. This map (defined up to scalar multiple) is simply the action of $u_{\pm
\omega}(t)$ followed by projection to $E^{\mu + m\omega}$.
Since any $\omega$-string in $K''$ meets $\pounds I$ at most once, it follows from Lemma~\ref{lem:tits} and Lemma~\ref{lem:indexsetdecomposition}.(ii)
that there exists a set $J$ satisfying
\begin{itemize}
\item[(J1)] $J \setminus I_{K''}$ coincides with $I \setminus I_{K''}$
\item[(J2)] If $\mu \in K'' \cap \pounds J$, then $\abs{\langle w \cdot \mu, \alpha\rangle}
\geqslant \ell_{r - 1} -1$
\end{itemize}
such that $sg \cdot \Vect{e}_{K}$ has nonzero component in the
$\Vect{e}_{J}(w ^{-1}\cdot \alpha)$-direction for almost every $s \in S_{\omega}$, $g \in G_{S}$.
In fact, we can say more:
\begin{itemize}
\item[(J2\mlprime)] If $\mu \in K'' \cap \pounds J$, then $\abs{\langle w \cdot \mu,
\alpha\rangle}$ equals $\ell_{r - 1} -1$ or $\ell_{r - 1}$.
\end{itemize}
In particular, every $w^{-1} \cdot \alpha$-chain $C$ meeting $I_{K''} \cap J$ has
length at least $\ell_{r - 1}$, and satisfies $\abs{J \cap C} < 2 \ell_{r -
2} + 2$. On the other hand, if $C$ does not meet $I_{K''} \cap J$, then 
$J \cap C = I \cap C \setminus I_{\textup{$\ell_{r - 1}$-small}}$.
The result now follows from $(w^{-1} \cdot \alpha)$-sparsity of $I \setminus I_{\textup{$\ell_{r -
1}$-small}}$ and a final application of Lemma~\ref{lem:beta<->beta'}.
\end{proof}

Proposition~\ref{prop:e_kVg.e_k}
provides a roundabout method to derive some
information about the direct sum decomposition of certain exterior power
representations of simple groups. Let $G$, $K$ and $E$ be as in that statement. Denote by $V_{K}$ the projectivization of the linear subspace $\bigoplus_{\mu
\in K} E^{\mu}$. A consequence of Proposition~\ref{prop:e_kVg.e_k}
is that $V_{K} \cap g \cdot V_{K} = \varnothing$ for almost all $g \in G$. In fact, it
follows from that proposition and the fact that $E$ is $\ell_{r}$-large that
$$
\dim E - 2 \sum_{\mu \in K} \dim E^{\mu} \geqslant \ell_{r} - 2 \ell_{r - 1}
$$
for any affine hyperspace $K$ of $\Lambda \otimes \mathbb{R}$.

\begin{proposition}\label{prop:h_u.eVe}
Take $G$, $K$, and $E$ as in Proposition~\ref{prop:e_kVg.e_k}, and take $h \in T$. 
Let $h = h_{s}h_{u} = h_{u}h_{s}$ be the Jordan decomposition of $h$.
If the semisimple component
$h_{s}$ is not torsion, then
the set of $u
\in G$ satisfying $$
h^{u} \cdot \Vect{e}_{K} \wedge \Vect{e}_{K} \neq 0$$
is nonempty and open.
\end{proposition}
\begin{proof}
Fix for the moment a basis $\mathscr{B}$ of $E$ consisting of weight vectors for the action
of $T$. Given a finite subset $I \subset \mathscr{B}$, we define
$\Vect{e}_{I} = \bigwedge_{\Vect{\upsilon} \in I} \Vect{\upsilon}$. Although only
defined up to sign, the vector $\Vect{e}_{I}$ is again a weight vector for the action of
$T$. We have already seen that the set $\mathscr{B}(m) = \{\Vect{e}_{I} \mid I \subset
\mathscr{B}, \abs{I} = m\}$ forms a basis for the $m$-fold exterior power $\textstyle{\bigwedge^{\!m}} E$ of $E$.

Take $m = 2 \sum_{\mu \in K} \dim E^{\mu} = 2 \dim V_{K} + 2$. Given elements
$u \in G$ and $X \in \Lie T$, it makes sense to consider the quantity
$\chi_{\mathscr{B}, I, u}(X)$, defined by
$$X \cdot u \cdot \Vect{e}_{K} \wedge u \cdot \Vect{e}_{K}
= \sum_{\substack{I \subset \mathscr{B}\\ \abs{I} = m}} \chi_{\mathscr{B}, I, u}(X)
\Vect{e}_{I}.$$
We note that if we regard elements of $\Lambda$ as functionals on $\Lie T$, then $\chi_{\mathscr{B}, I,
u}$ lies in $\Lambda \otimes \mathbb{C}$.

By Proposition~\ref{prop:e_kVg.e_k}, $g \cdot \Vect{e}_{K} \wedge
\Vect{e}_{K}$ is nonzero for almost all $g \in G$. By dominance of the morphism $G \times
T \to G$ given by $(u, t) \mapsto t^{u}$, there exist $I$ and $u$ such that $\chi_{\mathscr{B}, I, u}$ is not identically zero on $\Lie
T$.

If $\dot{w} \in G$ represents $w \in W(T, G)$, then we have the relation
$$w \cdot \chi_{\mathscr{B}, I, u} = \chi_{\dot{w} \cdot \mathscr{B}, \dot{w} \cdot I, \dot{w}u},$$
up to sign. We are now ready
to vary our basis $\mathscr{B}$. Since $\chi_{\mathscr{B}, I, u} \neq 0$ as a member of
$\Lambda \otimes \mathbb{C}$, it follows that $W \cdot \chi_{\mathscr{B}, I, u}$ spans 
$\Lambda \otimes \mathbb{C}$ over $\mathbb{C}$.

Upon conjugation, we may assume that $h_{s} \in T$.
We first prove the result in the special case of a semisimple element $h = h_{s}$.
By hypothesis, $h$ does not lie in the intersection of the kernels of all characters of
$T$. It follows from the preceding discussion that there exists a $w \in W$ such that
$h\dot{w}u \cdot \Vect{e}_{K} \wedge \dot{w}u \cdot \Vect{e}_{K}$ has nonzero component
in the $\Vect{e}_{\dot{w} \cdot I}$-direction with respect to the basis $(\dot{w} \cdot \mathscr{B})(m)$ of $\textstyle{\bigwedge^{\!m}}
E$. The result in the case $h = h_{s}$ now follows from connectedness of $G$.

Finally, we turn to the general case, $h_{u} \neq 1$. By the preceding argument, we see
that there exists a $u \in G$ satisfying the equation $h_{s}^{u} \cdot
\Vect{e}_{K} \wedge \Vect{e}_{K} \neq 0$. The computation
$$\left(\begin{array}{cc}
\lambda & 0\\
0 & \lambda^{-1}
\end{array}\right)
\left(\begin{array}{cc}
1 & \mu\\
0 & 1
\end{array}\right)
\left(\begin{array}{cc}
\lambda^{-1} & 0\\
0 & \lambda
\end{array}\right)
= \left(\begin{array}{cc}
1 & \lambda^{2}\mu\\
0 & 1
\end{array}\right)$$
for $\lambda$, $\mu \in \mathbb{C}^{\times}$
implies that there exists a sequence of matrices $t_{i} \in T$ such that $h_{u}^{t_{i}}
\to 1$ as $i \to \infty$. Thus
$$
h^{t_{i}u} \cdot \Vect{e}_{K} \wedge \Vect{e}_{K} = (h_{u}^{t_{i}u})h_{s}^{u} \cdot
\Vect{e}_{K} \wedge \Vect{e}_{K} 
\to h_{s}^{u} \cdot
\Vect{e}_{K} \wedge \Vect{e}_{K}
$$
since $h_{s}$ commutes with $T$.
\end{proof}
\subsection{Reduction modulo $k$}\label{sec:ReductionModk}
We wish to generalize the results of the previous section to the context of a simple
group $G$ over a local field $k$. To apply those results, we must first construct
universal objects corresponding to $G$, its representations, etc. Roughly speaking, we
can realize $G$ as the `extension of scalars' by $k$ of a Chevalley group defined over
$\mathbb{Z}$. Then we can extend scalars of the large representation of
Proposition~\ref{prop:e_kVg.e_k} to obtain a representation of $G$ with some favorable
properties. We will follow a somewhat different procedure.

To begin, let $k$ be an arbitrary field.
Every $k$-split, absolutely almost simple algebraic $k$-group $G$ comes from a group scheme
$H$ defined over
$\mathbb{Z}$ by change of base $\mathbb{Z} \to k$~\cite{MR1611814}. As usual, we identify $H$ with the
group of its points over $\mathbb{C}$. We can endow the
universal enveloping algebra
of $H$ with a $\mathbb{Z}$-structure such that any finite-dimensional, $k$-rational
representation $\rho$ of $G$ is obtained by extending scalars of the
corresponding representation of $H$ from $\mathbb{Q}$ to $k$, and then factoring out the
unique maximal invariant $k$-subspace~\cite{MR0207713}.

The following lemma is a special case of the corollary to Proposition~14, Chapter~II,
\S7.7 of~\cite{MR1727844}.
\begin{lemma}\label{lem:cor.toII.7.7.14}
Let $V$ be a rational vector space. If $W$ and $W'$ are linear subspaces of $V$, then $$(W
\cap W') \otimes k = (W
\otimes k) \cap(W'
\otimes k).$$
\end{lemma}

Let $T$ be a maximal $k$-split
torus, and write $\Delta$ for the root system $\Delta(T, G)$, which coincides with the
root system of the Chevalley scheme $H$. Denote by $\Lambda$ the
corresponding weight lattice. Given a $G$-module $E$ and a subset $K \subset \Lambda
\otimes \mathbb{R}$, denote by $V_{K}$ the projectivization of the linear subspace $\bigoplus_{\mu
\in K} E^{\mu}$.

Let $\ell_{r}$ be as in the statement of
Proposition~\ref{prop:e_kVg.e_k}.
\begin{corollary}\label{cor:localversionofproposition}
Let $G$ be a connected, simply laced, $k$-split, absolutely almost simple algebraic
$k$-group of rank $r$; T a maximal $k$-split torus of $G$; and
$\rho$ an $\ell_{r}$-large, $k$-rational
representation of $G$ on a vector space $E$. Suppose $h \in G(k)$ has
semisimple Jordan component of infinite order.

Denote by $X$ the set of elements $u \in G$ such that 
\begin{equation}\label{eq:VkVk-prime}
V_{K} \cap h^{u} \cdot V_{K'} = \varnothing
\end{equation}
for every pair of parallel affine hyperspaces $K$,
$K'$ of $\Lambda \otimes \mathbb{R}$.
Then $X$ is nonempty and open.
\end{corollary}
\begin{proof}
Since $G$ is $k$-split, the subspaces $V_{K}$ and $V_{K'}$ are defined over $k$. It follows
that the condition~\eqref{eq:VkVk-prime} is an open condition on $u \in G$. Therefore,
we need only show that this condition is nonvacuous for each pair $K$, $K'$ of parallel
affine hyperspaces.

Suppose first that $K \neq K'$.
Then there exists a choice of a system $\Delta^{+}$ of positive roots of $G$ with respect
to $T$ such that $\mu' - \mu > 0$ for any weights $\mu \in K$, $\mu' \in K'$. If $B$ is the
Borel subgroup corresponding to this choice of $\Delta^{+}$, then any $u$ which
conjugates $h$ into $B$ will satisfy condition~\eqref{eq:VkVk-prime}. We are left with
the case where $K = K'$.

Let us denote by $\rho_{\mathbb{C}}$, $E_{\mathbb{C}}$ the corresponding complex
representation of the associated Chevalley scheme $H$, which we identify with its
points over $\mathbb{C}$. Denote by $E_{\mathbb{Q}}$ the rational structure coming
from the $\mathbb{Z}$-structure on the universal enveloping algebra of $H$.
Write $$Y = \{v \in H \mid W_{K} \cap h^{v} \cdot W_{K} = \varnothing\},$$
where $W_{K} = \bigoplus_{\mu \in K} E_{\mathbb{Q}}^{\mu}$.

It follows from Lemma~\ref{lem:cor.toII.7.7.14} that 
\begin{gather*}
\mathbb{C}(Y) = \mathbb{Q}(Y) \otimes \mathbb{C}\\
k(X) = \mathbb{Q}(Y) \otimes k.
\end{gather*}
But $Y$ is nothing but $\{v \in H \mid \Vect{e}_{K} \wedge h^{v} \cdot \Vect{e}_{K} \neq
0\}$. Thus Proposition~\ref{prop:h_u.eVe} implies that $\mathbb{C}(Y) =
\mathbb{C}(H)$. The result now follows from dimension considerations.
\end{proof}

\section{Producing representations with contractive
dynamics}\label{sec:ContractiveDynamics}
\subsection{Producing representations with proximal behavior}\label{sec:Pi-im}
We begin this section by outlining some of the elementary representation theory of
semisimple
groups, which can be found in more detail in \S24.B of~\cite{MR1102012}. 
We restore $k$ to being a local field.
Let $G$ be a connected, semisimple,
$k$-split $k$-group, $T$ a maximal $k$-split torus.
We write $\Delta =
\Delta(T, G)$ and $W = W(T, G)$ for the root system and Weyl group associated to $T$.

We fix for the moment a system
$\Delta^{+}$ of positive roots for $\Delta$. Let $B$ denote the corresponding Borel
$k$-subgroup containing $T$, and $U$ its unipotent radical. If
$\lambda_{0}$ is a character of $T$, we can regard it as a function on $B = TU$ by
extending it trivially on $U$.

Fix a $k$-rational character $\lambda$ of $T$, dominant with respect to $\Delta^{+}$.
 Let $E(\lambda)$ denote the Weyl module
associated to $\lambda$. Namely,
$$E(\lambda) = \{f \in \Omega[G] \mid f(gb) = (-w_{0} \cdot \lambda)(b)\cdot f(g), g \in G, b
\in B\}.$$

First we suppose that $k$ has characteristic zero. Then $E(\lambda)$ is an absolutely irreducible $G$-module with highest
weight $\lambda$.
Because $G$ is $k$-split, we can endow $E(\lambda)$ with a $k$-structure with respect to
which the representation $\rho: G \to \mathrm{GL}(E(\lambda))$
is rational. That $E(\lambda)$ is finite-dimensional (and
nonzero) follows from dominance of $\lambda$.

On the other hand, if $k$ has positive characteristic, then the Weyl module $E(\lambda)$
need not be irreducible. Nevertheless, $E(\lambda)$ contains a unique absolutely
irreducible subrepresentation $F(\lambda)$. As before, we have that $F(\lambda)$ is
finite-dimensional and nonzero; and that $\lambda$ is still the highest weight of
$F(\lambda)$. In this case it is also true (cf.\ Theorem~2.5 of~\cite{MR0277536})
that $F(\lambda)$
can be endowed with a $k$-structure with respect to which $\rho:G \to \mathrm{GL}(F(\lambda))$
is rational.

In the remainder of this section, we consider the
problem of producing a representation wherein a fixed
semisimple element of $G(k)$ acts very proximally. Take $h$ in $T(k)$.

We wish to associate some root data with $h$. Choose a system $\Delta^{+}$ of positive
roots containing all roots $\phi \in \Delta$ such that $\abs{\phi(h)} > 1$, where
$\abs{\cdot}$ denotes the absolute value on $k$. Let $\Pi$ be the associated set of
simple
roots. Define the set $\Delta_{im} = \lbrace \phi \in \Delta \mid \abs{\phi(h)}
= 1\rbrace$. This is easily seen to be a root subsystem of $\Delta$, and
$\Delta_{im}^{+} = \Delta_{im} \cap \Delta^{+}$ is a valid choice of positive system.
It is also true that the simple roots associated to this choice of positivity is exactly
$\Pi_{im} = \Delta_{im} \cap \Pi$.

The main tool for constructing projective actions where a given nontorsion semi\-simple
element admits a ping-pong partner is the following observation.

\begin{proposition}\label{prop:key}
Suppose there exist an element
$w \in W$ and a dominant weight $\lambda$ satisfying the following conditions.
\begin{itemize}
\item[(i)] $w\cdot\lambda = -\lambda$
\item[(ii)] $w\cdot\Delta_{im}=\Delta_{im}$
\item[(iii)] $\lambda \perp \Pi_{im}$
\end{itemize}
Then there exists an absolutely irreducible, finite-dimensional
$G$-module $E$ and a
$k$-structure on $E$ with respect to which the map
$\rho:G \to \mathrm{GL}(E)$ is $k$-rational and 
$\rho(h) \in \mathrm{GL}(E_{k})$ is very proximal.
\end{proposition}
\begin{proof}
Let $E = E(\lambda)$ if $\operatorname{char} k = 0$; otherwise, let $E = F(\lambda)$. By the
discussion above, the representation $\rho: G \to \mathrm{GL}(E)$ is $k$-rational,
absolutely irreducible, and finite-dimensional.
We need only show that $\rho(h)$ is very proximal.

Any weight $\mu$ of $\rho$ must be of the form
\begin{equation}\label{eqn:weight}
\mu = \lambda - \sum_{\alpha \in \Pi} c_{\alpha} \alpha
\end{equation}
with $c_{\alpha} \geqslant 0$.
The definition of $\Delta^{+}$ immediately implies that $\lambda(h)$
is an eigenvalue of maximal norm, and that if $\abs{\mu(h)} =
\abs{\lambda(h)}$, then all coefficients $c_{\alpha}$ in \eqref{eqn:weight} vanish
for $\alpha \in \Pi \setminus \Pi_{im}$.

Denote by $E^{\mu}$ 
the weight space of $E$ corresponding to weight $\mu$. In this language, the
projectivization $\widehat{E^{\lambda}}$ of $E^{\lambda}$ lies in $A(\rho(h))$. We work to
establish the claim that if
\begin{equation}\label{eqn:maxeigs}
\mu = \lambda - \!\!\!\sum_{\alpha \in \Pi_{im}}\!\!\! c_{\alpha} \alpha
\end{equation}
is a weight for some integers $c_{\alpha}$, then all of the coefficients
$c_{a}$ are zero. It will follow
that $A(\rho(h))$ coincides with the projectivized weight space $\widehat{E^{\lambda}}$. Since $\lambda$
has multiplicity one as a weight of $\rho$, it will follow that the $\rho(h)$ is
proximal.

We have already observed that $\Delta_{im}$ is itself an abstract root system; that
$\Delta^{+}_{im} = \Delta_{im} \cap \Delta^{+}$ is a valid notion of positivity in $\Delta_{im}$;
and that $\Pi_{im} = \Pi \cap \Delta_{im}$ is a system of simple roots. Denote
by $W_{im}$ the Weyl group of this system. In particular, for any element $z = \sum_{\alpha \in \Pi_{im}} z_{\alpha}
\alpha$
of the root lattice of $\Delta_{im}$, there exists an element $w_{z} \in W_{im}$ such
that $w_{z}\cdot z$ lies in the closed fundamental Weyl chamber corresponding to
$\Pi_{im}$.

We note that the Weyl group $W_{im}$ of the subsystem $\Delta_{im}$ is generated by the
reflections $\sigma_{\alpha}$ through hyperspaces perpendicular to simple roots $\alpha
\in \Pi_{im}$. Although $W_{im}$ was defined abstractly in terms of a root system, we may
identify $W_{im}$ with the corresponding subgroup of $W$: namely,
the subgroup generated by the root reflections $\sigma_{\alpha}$ (a minor abuse of
notation) associated to simple roots in $\Pi_{im}$. Now for all $\alpha \in \Pi_{im}$,
condition (iii) implies that $\sigma_{\alpha}$ fixes $\lambda$. So if we take
$z = -\sum_{\alpha \in \Pi_{im}} c_{\alpha} \alpha$, then
$$w_{z}\cdot \mu = \lambda + w_{z} \cdot z \succcurlyeq \lambda,$$
with respect to the notion of positivity coming from $\Delta^{+}$. If $\mu$ is a weight
of $\rho$, then of course so is $w_{z} \cdot \mu$. However, $w_{z}
\cdot \mu$ cannot be a higher weight than $\lambda$. Therefore, $z = 0$ and $\mu =
\lambda$. This proves our claim, and establishes the proximality of $\rho(h)$.

Condition (i) implies in particular that $-\lambda$ is a weight of $\rho$. If there were
a lower weight $\mu$, then $w_{0}\cdot \mu$ would be higher than $\lambda$. 
By construction of $\Delta^{+}$, we see that $\abs{\mu(h)}$ is minimized over
all weights $\mu$ of $\rho$ when $\mu = - \lambda$, that is, when $\abs{\mu(h)} =
\abs{\lambda(h)}^{-1}$. Moreover, if $\mu$ is a weight realizing $\abs{\mu(h)} =
\abs{\lambda(h)}^{-1}$, then $\mu$ satisfies $\mu + \lambda = \sum_{\alpha \in
\Pi_{im}} c_{\alpha} \alpha$. Applying $w$ to both sides yields
$$w\cdot \mu = \lambda + w \cdot \!\!\!\!\sum_{\alpha \in
\Pi_{im}}\!\!\! c_{\alpha} \alpha.$$
But condition (ii) implies that $w$ preserves the root lattice corresponding to the
subsystem $\Delta_{im}$. Therefore, $w\cdot\sum_{\alpha \in
\Pi_{im}} c_{\alpha} \alpha = \sum_{\alpha \in
\Pi_{im}} c'_{\alpha} \alpha$. But we have already observed that the only weight of the
form~\eqref{eqn:maxeigs} is $\lambda$ itself. Thus $w\cdot \mu = \lambda$, or
equivalently, $\mu = -\lambda$. Thus $A(\rho(h^{-1}))$ coincides with the projectivization
of $E^{-\lambda}$. Since $\lambda$ has multiplicity one, so does $w\cdot \lambda =
-\lambda$, and the result follows.
\end{proof}
\begin{remark}
It often happens that a Weyl group contains the linear transformation $-1$ sending each
root to its negative. In this case, we note that the hypotheses of this theorem are much
easier to verify. Conditions (i)--(iii) can be replaced with the much simpler condition
that $\Pi_{im} \subsetneq \Pi$.
\end{remark}

Of course, this Proposition is not very useful if $\Pi_{im}$ coincides with $\Pi$. We
conclude this section by quoting a lemma which will eventually serve as an ``escape
hatch'' from this situation.

\begin{lemma}[Lemma 4.1 of~\cite{MR44:4105}]\label{lem:hatch}
Let $\ell$ be a finitely generated field, and $t \in \ell^{\times}$ an element of
infinite multiplicative
order. Then there exists a locally compact field $k$ endowed with an absolute value $\omega$, and
a homomorphism $\sigma: \ell \to k$ such that $\omega(\sigma(t)) \neq 1$.
\end{lemma}

\subsection[Quasi-projective transformations, quasi-proximality]{Quasi-projective
transformations and quasi-proximal unipotents}\label{sec:quasiProxUnipotent}
If $G$ is a connected, semisimple algebraic group over an archimedean local field $k$, $S
\subset G(k)$ a $k$-dense subgroup, and $u$ a unipotent
(and nontorsion) element of $S$, then our eventual goal is to produce a irreducible, 
finite-dimensional, $k$-rational linear representation $\rho$ of
$G$ such that $\rho(u)$ behaves ``as if it were proximal,'' and such that there exists
genuinely proximal elements in $\rho(S)$. By the phrase ``as if it were proximal,'' we
mean roughly that $\widehat{\rho}(u)$ acts on the projective space with an isolated fixed
point; that the basin of attraction for this fixed point is large; and that
convergence to this fixed point occurs uniformly on compact sets in the appropriate
sense.

We begin with a brief review of the notions of quasi-projective
transformations and contractions. The notion was introduced by
Furstenberg~\cite{MR0352328}, in the same paper where he introduced proximality. The
basic references for this section are the papers of 
Gol\mlprime dshe{\u\i}d and Margulis~\cite{MR1040268}, and Abels, Margulis, and
So{\u\i}fer~\cite{MR1348303}.

In this section, we let $E$ denote a finite-dimensional vector space over an ar\-chi\-medean
local field $k$, $P$ its projectivization. Endow $P$ with the angle metric
$$d(p, q) = \arccos\frac{\abs{\langle \Vect{p}, \Vect{q}\rangle}}{\norm{\Vect{p}}
\cdot\norm{\Vect{q}}}$$
where $\Vect{p}$ and $\Vect{q}$ are representatives in $E\setminus \{0\}$ of
$p$ and $q$ respectively;
 $\abs{\langle \cdot, \cdot\rangle}$ is a fixed Euclidean
(resp.\  Hermitian) scalar product on $E$ if $k = \mathbb{R}$ (resp., if $k =
\mathbb{C}$); and $\norm{\cdot}$ is the norm on $E$ associated to the scalar
product $\abs{\langle \cdot, \cdot\rangle}$.

\begin{definition}
A map $b: P \to P$ is called a \emph{quasi-projective transformation} if it is the pointwise
limit of a sequence of projective transformations of $P$.
\end{definition}

If $b$ is a quasi-projective transformation, we denote by $M_{0}(b)$ the image
under $b$ of the set of points of continuity for $b$. If $s = (s_{n})$ is a a sequence in $\mathrm{GL}(E)$
with $b = \lim \widehat{s_{n}}$ then we set $M_{0}(s) = M_{0}(b)$. Let $L_{1}(s)$ denote
the topmost space in the Lyapunov filtration associated to $s$: Namely, if $\beta$ is
the linear operator on $E$ given by
$$\beta(x) = \lim_{n \to \infty} \norm{s_{n}}^{-1} s_{n}$$
(which exists thanks to compactness of the unit ball in $E$), then $L_{1}(s)$ is by
definition the projectivization in $P$ of $\ker
\beta$.

\begin{definition}
We call $b$ a \emph{contraction} if $M_{0}(b)$ is a point.
\end{definition}

\begin{definition}
We call a sequence $s = (s_{n})_{n > 0}$ in $\mathrm{GL}(E)$ is \emph{contractive}
if the following conditions are satisfied:
\begin{itemize}
\item[(i)] The associated sequence $(\widehat{s_{n}})$ converges
pointwise on $P$; and
\item[(ii)] there is a sequence of scalars $(B_{n})$ in $k^{\times}$ such that
$(B_{n}s_{n}) \in \mathrm{GL}(E)$ converges pointwise to a linear map $\beta$ of rank
$1$.
\end{itemize} 
\end{definition}

\begin{remark}
If $s$ is a contractive sequence, then $b$ is a contraction, and moreover
agrees with $\widehat{\beta}$
on the complement of $L_{1}(s)$. Conversely, it is true that any
contraction is the limit of a contractive sequence.
\end{remark}

The following lemma is evident.

\begin{lemma}\label{lem:4basic}
Let $b$ be a quasi-projective transformation with
$b = \lim \widehat{s_{n}}$ for some sequence $s = (s_{n})$ in $\Sigma$; and take $h \in
\Sigma$.
\begin{itemize}
\item[(i)]  $M_{0}(hs) = h \cdot M_{0}(s)$
\end{itemize}

If in addition $s$ is a contractive sequence, then
\begin{itemize}
\item[(ii)] $L_{1}(hs) = L_{1}(s)$
\end{itemize}
\end{lemma}

We shall exploit the interplay between the notions of proximality and
contraction. For a subsemigroup $\Sigma$ of $\mathrm{GL}(E)$, let us denote by $\overline{\Sigma}$
the set of quasi-projective transformations arising as limits of sequences of projective
transformations induced by elements of $\Sigma$.
We quote the following result.
\begin{lemma}[Lemma 3.13 of~\cite{MR1348303}]\label{lem:3.13}
If $\Sigma$ is an irreducible subsemigroup of $\mathrm{GL}(E)$ and
$s = (s_{n})$ is a contractive sequence in $\Sigma$, then $\Sigma$ contains a proximal
element.
\end{lemma}
\begin{proof}
Upon replacing $s$ by $hs$ for some $h \in  \Sigma$, we may assume by
Lemma~\ref{lem:4basic} and Lemma~\ref{lem:tits} that
$M_{0}(s) \notin L_{1}(s)$. It now follows from Lemma~\ref{lem:3.8}.(ii) that $s_{n}$
is proximal for $n$ sufficiently large.
\end{proof}

Suppose that $u \in \mathrm{GL}(E)$ is a unipotent linear transformation. If $x \in P$ then $\widehat{u^{n}}x$
converges pointwise to an eigendirection for $u$; i.e., to a point in the
projectivization in $P$ of $\ker
u - I$.

\begin{definition} We shall say a unipotent element $u$ is \emph{quasi-proximal} if for
some integer $d$, the
endomorphism $u_{d} = (u - I)^{d}: E \to E$ has rank one.
\end{definition}

We shall assume for the remainder of this section that $u$ is quasi-proximal. Then if
$d$ is as in the definition,
$\binom{n}{d}^{-1}u^{n}$ converges pointwise to a nonzero scalar multiple of $u_{d}$.
Thus, the sequence $s = (u^{n})$ is contractive.
Denote by $B(u)$ and $B'(u)$ the projectivizations of $\operatorname{im} u_{d}$ and $\ker
u_{d}$, respectively. We have $M_{0}(s) = B(u)$ and $L_{1}(s) = B'(u)$. We note that we have the
following analogue for quasi-proximal unipotents of Tits' Lemma~\ref{lem:3.8}.(i):

\begin{lemma}\label{lem:quasiprox}
Suppose that $u \in \mathrm{GL}(E)$ is a quasi-proximal unipotent. Let $K \subset P
\setminus B'(u)$
be a compact set. For any $q > 0$ there exists an integer $N$ such that $\norm{
\widehat{u}^{z}|_{K}} < q$ for all $z > N$; and for every neighborhood $U$ of $B(u)$
there exists an integer $N'$ such that $u^{z}\cdot K \subset U$ for all $z > N'$.
\end{lemma}
\begin{proof}
The second assertion follows from the observation that the quasi-projective
transformation
$\lim_{z \to \infty}\widehat{u}^{z}$ agrees with $\widehat{u_{d}}$
on the complement in $P$ of $L_{1}(\{u^{n}\})$.

The first assertion requires more care. Let $p, q \in P$, and pick representatives $\Vect{p},
\Vect{q} \in E \setminus \{0\}$. We first observe
that if $\Vect{p}, \Vect{q} 
\in E \setminus \{0\}$, then we can regard
$\abs{\langle u^{z} \Vect{p}, u^{z} \Vect{q}\rangle}$ as a polynomial in $z$
with coefficients in $k$. Moreover, if $p, q \in P\setminus B'(u)$, then the degree of
this polynomial is precisely $2d$. It follows that the expression
$\cos^{2} d(\widehat{u}^{z}p, \widehat{u}^{z}q) \in k(z)$ is a ratio of polynomials
of degree $4d$ in $z$, with identical leading coefficients. (In a suitable basis the
coefficients
are homogeneous polynomials of degree four in $\Vect{p}$ and $\Vect{q}$.)

Let $K$ be as in the statement of the theorem. From the above description of $d(\widehat{u}^{z}p,
\widehat{u}^{z}q)$, we see that if $z$ is large relative to the distance $d(K, B'(u))$, 
then
$$\lim_{p \to q} \frac{d(\widehat{u}^{z}p, \widehat{u}^{z}q)}{d(p, q)} = 0$$
for $q \in K$.

Consider the functions $f_{z}: K \to \mathbb{R}_{+}$ given by
$$f_{z}(p) = \sup_{q \in K\setminus \{p\}}
\frac{d(\widehat{u}^{z}p, \widehat{u}^{z}q)}{d(p,q)}.$$
The function $f_{z}$ is continuous because if $p$ is near $q$, then either
$$\frac{d(\widehat{u}^{z}p, \widehat{u}^{z}p')}{d(p, p')}\text{ is close to }
\frac{d(\widehat{u}^{z}q, \widehat{u}^{z}p')}{d(q, p')}$$
or $q = p'$.

We observe that rational functions in $z$ have the property that they are eventually
monotone. From the remarks above, it follows that the sequence $f_{z}$
eventually decreases in $z$ monotonically 
to zero. It follows from Dini's theorem that this convergence is uniform. The result follows.
\end{proof}

Let us summarize the results of this section.
\begin{corollary}\label{cor:qpu}
If $\Sigma$ is an absolutely irreducible subsemigroup of $\mathrm{GL}(E)$ which contains
quasi-proximal unipotent elements, then it contains proximal elements.
\end{corollary}
\begin{proof}
Apply Lemma~\ref{lem:3.13}.
\end{proof}

It may be illuminating to supply a more direct proof of Corollary~\ref{cor:qpu}.
Suppose $u$ is a quasi-proximal unipotent. By irreducibility, there exists an element $h
\in \Sigma$ such that $h\cdot B(u) \not\subset B'(u)$. We will see that $hu^{z}$ is
proximal for all sufficiently large $z$.

Let $K \subset P$ be a compact set such that $B(u) \subset \Int K$ and
$h\cdot K \subset P\setminus B'(u)$. By Lemma~\ref{lem:quasiprox} there exists a number  
$N$ such that $u^{z}h\cdot K \subset \Int K$ and $\norm{
\widehat{u}^{z}|_{h\cdot K}} < 1/\norm{
\widehat{h}|_{K}}$ for all $z > N$. Then
$hu^{z}h\cdot K \subset h\cdot\Int K$ and $\norm{
\widehat{hu}^{z}|_{h\cdot K}} < 1$. We conclude by invoking Lemma~\ref{lem:3.8}.(ii).

In closing, we note that if a unipotent element is quasi-proximal, then so is its inverse.
Moreover, $B(u) = B(u^{-1})$ and $B'(u) = B'(u^{-1})$. We also have for any $x \in
\mathrm{GL}(E)$ and $z \in \mathbb{N}$ that $B(\mbox{}^{x}u) = x\cdot B(u)$,
$B(u^{z}) = B(u)$, and likewise for $B'$.

\section{Proof of the main theorem}\label{sec:proof}
We say a group is \emph{icc} if all of its nonidentity conjugacy classes are infinite.
Bekka and de la Harpe~\cite{MR2001j:46078} showed the following:
\begin{lemma}\label{lem:goingUp}
Suppose that $\Gamma_{0}$ 
is a finite-index subgroup of an icc group $\Gamma$. If $C_{r}^{*}(\Gamma_{0})$ is simple
(resp.\ has a unique trace up to normalization), then the same is true of $C_{r}^{*}(\Gamma)$.
\end{lemma}

\begin{lemma}\label{lem:icc}
Let $H$ denote the union of the finite conjugacy classes of a group $\Gamma$. Then $H$
is a characteristic, amenable subgroup of $\Gamma$.
\end{lemma}
\begin{proof}
The conjugacy class of a product lies in the product of conjugacy classes. Thus, it
remains only to prove that $H$ is amenable.

If $h \in \Gamma$ lies in a finite conjugacy class, then its centralizer $Z_{\Gamma}(h)$
is of finite index in $\Gamma$. If $S$ is a finite subset of $H$, then $\bigcap_{h \in
S} Z_{\Gamma}(h)$ is likewise of finite index in $\Gamma$. But this intersection
contains the center of the subgroup $\langle S \rangle$. It follows that every finitely
generated subgroup of $H$ is virtually abelian, and hence amenable. The result follows.
\end{proof}

We begin by deducing Theorem~\ref{thm:main} from the following result:
\begin{theorem}\label{thm:subgroup}
Let $G$ be a connected, semisimple algebraic group with trivial center, and let $\Gamma <
G$ be a finitely generated,
Zariski-dense subgroup. Then $C_{r}^{*}(\Gamma)$ is simple, and has a unique trace up
to normalization.
\end{theorem}
\begin{proof}[Proof of Theorem~\ref{thm:main}]
Paschke and Salinas~\cite{MR82c:22010} showed that if $\Gamma$ has a nontrivial normal
amenable subgroup, then $C_{r}^{*}(\Gamma)$ is not simple, nor is its trace unique up to
normalization.

Suppose $\Gamma$ has no nontrivial normal
amenable subgroup. Let $G$ be the Zariski closure of $\Gamma$. If we denote by $R$
the radical of $G$, then $\Gamma \cap R$ is a solvable normal subgroup of $\Gamma$,
hence trivial. Therefore, upon projecting modulo $R$, we may suppose that $G$ is
semisimple. Likewise, up to projecting to the adjoint group, we may assume that $G$ is
center-free.

For convenience, for a subgroup $H$ of $G$, let us denote by $\mathrm{Cl}(H)$ the Zariski closure
of $H$ and by $H^{\circ}$ and the connected component of $H$ containing the identity.
For any infinite subgroup $\Gamma_{0}$ of $\Gamma$, we have that the dimension of $\mathrm{Cl}(\Gamma_{0})$ is
positive. Since $\Gamma$ is nonamenable, it must contain an infinite, finitely generated
subgroup. Write
$$d = \max\{\dim \mathrm{Cl}(\Gamma_{0}) \mid \Gamma_{0}\text{ is finitely generated}\}.$$
We can write $\Gamma$ as the union of an increasing chain $\Gamma_{0} \subset \Gamma_{1} \subset
\cdots$ of finitely generated
subgroups. Write $G_{0} = \mathrm{Cl}(\Gamma_{0})^{\circ}$. Without loss of generality, we
may suppose that $\dim G_{0} = d$.

A theorem of Schur has it that if a linear torsion group is finitely generated, then it is
finite. In particular, $\Gamma_{0} \cap G_{0}$ contains an element of infinite order. It
follows that $\Gamma_{i} \cap G_{0}$ is nontorsion (i.e.\ not locally finite) for all $i
\geqslant 0$.

For $\gamma \in \Gamma$, we have $\langle \Gamma_{0}, \Gamma_{0}^{\gamma} \rangle \subset
\langle \Gamma_{0},\gamma \rangle$. Taking closures, we obtain
$$\langle G_{0}, G_{0}^{\gamma} \rangle \subset
\mathrm{Cl}(\langle \Gamma_{0},\gamma \rangle).$$
But the latter is an irreducible algebraic variety of dimension $d$. Since $G_{0}$ is
connected, we must have $G_{0} = G_{0}^{\gamma}$. Since $\gamma$ was arbitrary, we have
$\Gamma \subset N_{G}(G_{0})$. But this last is closed, while $\Gamma$ is dense. Hence
$G_{0}$ is a normal subgroup of $G$.

For an element $a \in G$, we have $\abs{a^{G_{0} \cap \Gamma_{i}}} \leqslant
\abs{a^{\Gamma_{i}}}$. But $a^{G_{0} \cap \Gamma_{i}}$ is Zariski dense in $a^{G_{0}}$.
In particular, if $a$ lies in a finite conjugacy class of $\Gamma_{i}$, then $a \in
Z_{G}(G_{0})$ = 1. Therefore, $\Gamma_{i}$ is icc.

The remainder of the proof consists in exhibiting an arbitrary Zariski-dense subgroup
$\Gamma$ of an arbitrary semisimple adjoint group $G$ as the increasing union of
subgroups whose reduced $C^{*}$-algebras are simple and have a unique tracial state. The
result will then follow by taking inductive limits.

We can write $G^{\circ}$ as an internal direct product of simple groups, say $G^{\circ}
= G_{1} \times \cdots \times G_{r}$. Induction on $r$: If $r = 1$, then $G_{0} =
G^{\circ}$. We can take $\Gamma_{i}$ Zariski dense for $i \geqslant \abs{G/G^{\circ}}$.
The result in this case follows from Theorem~\ref{thm:subgroup}.

Therefore, suppose that $r > 1$, and indeed that $G_{0} < G^{\circ}$. Then $\Gamma_{i}
\cap G_{0}$ is Zariski dense in $G_{0}$ and of finite index in $\Gamma_{i}$. By the
induction hypothesis, $C_{r}^{*}(\Gamma_{i} \cap G_{0})$ is simple and has a unique
trace, up to normalization. But then the algebra $C_{r}^{*}(\Gamma_{i})$ enjoys the same
properties by Lemma~\ref{lem:goingUp}.
\end{proof}

We shall prove that groups as in the statement of Theorem~\ref{thm:subgroup} `virtually' satisfy the
hypotheses of Lemma~\ref{lem:modifiedBCD}. Namely, if $\Gamma$ is such a group, then there exists a finite index subgroup $\Gamma_{0}$ such that for any finite subset $F \subset \Gamma_{0}$, there
exists an element $g \in \Gamma$ of infinite order such that etc. But, one may wish to know
more: For a fixed finite subset $F$, how large is the set of $g \in \Gamma$ verifying
the conditions of that lemma? We shall prove the following generalization of
Theorem~\ref{thm:subgroup}; it includes the statement that the set of such $g$ is very
large indeed. 

\begin{theorem}\label{thm:profinitedensity}
Let $G$ and $\Gamma$ be as in Theorem~\ref{thm:subgroup}, and let $F$ be a finite subset
of $G$. Suppose that for every simple direct factor $H$ of $G$, and every torsion
element $h \in F$, the centralizer $Z_{H}(h)^{\circ}$ is not a simple, proper subgroup of
full rank in $H$. Then the set of $g \in
\Gamma$ which verify conditions~(i) and (ii) of
Lemma~\ref{lem:modifiedBCD} with respect to some
decomposition $F = F_{p} \cup F_{r}$ is dense in the profinite topology on $\Gamma$.
\end{theorem}
\noindent Before turning to the proofs of Theorems~\ref{thm:subgroup} and~\ref{thm:profinitedensity},
we prove an important special case.
\begin{theorem}\label{thm:zariski-densePPPs}
Let $G$, $\Gamma$, and $F$ be as in Theorem~\ref{thm:profinitedensity}.
Write $G$ as the almost direct product $G_{-1}G_{+1}$,
where $G_{-1}$ is the largest normal subgroup of $G$ whose Weyl group contains the
linear transformation $-1$. Suppose that $G_{+1}$ contains no
nontrivial power of any element of $F$ whose semisimple part is nontorsion.

Denote by $Y_{F}$ the set of elements $\gamma \in \Gamma$ of
infinite order such that the subgroup
$\langle \gamma, h\rangle$ generated by $\gamma$ and $h$ is canonically isomorphic to the
free product
$\langle \gamma \rangle * \langle h\rangle$, all $h \in F$. Then $Y_{F}$ is dense in the
profinite topology on $\Gamma$.
\end{theorem}
\begin{proof}
Let $Hg$ be a coset of a finite-index subgroup $H \subset \Gamma$, $g\in \Gamma$. Upon
replacing $H$ be a finite-index subgroup, we may assume $H$ is normal in $\Gamma$; and we
note that $H$ is itself Zariski-dense in $G$, by virtue of connectedness of $G$.

Let $k$ be a finitely generated field over which $G$ is defined and such that $\Gamma
\subset G(k)$. Upon extending $k$, we may also assume that $k$ contains all eigenvalues
of all elements of $F$. Let $T$ be a maximal $k$-split torus of $G$. We write $F$ as a
disjoint union $F_{s} \cup F_{t} \cup F_{u}$, where
\begin{align*}
F_{s} &= \{h \in F \mid h\text{ has nontorsion semisimple part $h_{s}$}\},\\
F_{t} &= \{h \in F \mid h\text{ is torsion}\},\text{ and}\\
F_{u} &= \{h \in F \mid \text{$h$ is not torsion, $h_{s}$ is torsion}\}.
\end{align*}

Suppose first that $h \in F_{s}$. Upon conjugating, we may assume $T$
contains $h_{s}$. Then writing $h = h_{-1}h_{+1}$ with $h_{-1} \in G_{-1}$, then some
eigenvalue $\tau$ of $h_{-1}$ has infinite multiplicative order. By Lemma~\ref{lem:hatch},
there exists a local field $k_{h}$
extending $k$ such that $\abs{\tau} \neq 1$. If $\Pi_{im} =
\Pi_{im}(h_{s})$ is as in \S\ref{sec:Pi-im}, then $\rho_{h}(h)$ is proximal in the irreducible
representation $\rho_{h}$ of highest weight $\lambda$ for any positive integral sum
$\lambda$ of fundamental weights $\varpi_{\alpha}$ with $\alpha \in \Pi
\setminus \Pi_{im}$. By choice of $k_{h}$, there exists a connected component of the Dynkin
diagram, not contained in $\Pi_{im}$, whose Weyl group contains $-1$. In particular, by
constraining the coefficients
$\langle \lambda, \alpha \rangle$ in the sum to be nonzero only for those $\alpha$ lying in
this
connected component, we obtain that $\rho_{h}(h)$ is very proximal, and that $\rho_{h}$
factors through an absolutely simple factor of $G$.

Suppose $h \in F_{u}$. It follows that $h_{u}$ is nontorsion, whence $k$ has
characteristic zero. Let
$k_{h}$ be an archimedean local field extending $k$. Let $m_{h}$ be the order of $h_{s}$.
Let $\rho_{h}$ be any
 $k_{h}$-rational, absolutely irreducible, finite-dimensional representation of $G$.
Upon replacing $\rho_{h}$ by an exterior power and taking the correct irreducible
component, we can guarantee that $\rho_{h}(h^{m_{h}})$ is a quasi-proximal unipotent.
By Corollary~\ref{cor:qpu}, $\Omega_{+}(\rho_{h}, \Gamma)$ is nonempty. As before, we
may assume that $\rho_{h}$ factors
through an absolutely simple factor of $G$.

Finally, we turn to the torsion elements $F_{t}$. Denote by $\widetilde{F}_{t}$ the set
of all nontrivial projections of all elements of $F_{t}$ to all simple direct factors of
$G$. We remark that, in order to find an element of $\Gamma$ which satisfies no
nontrivial relation with any $h \in F_{t}$, it suffices to find an element of $\Gamma$ which satisfies no
nontrivial relation with any $h \in \widetilde{F}_{t}$. We also remark that $k$ already
contains all eigenvalues of all elements of $\widetilde{F}_{t}$. Therefore, we assume
that $\widetilde{F}_{t} = F_{t}$.

Fix an element $h \in F_{t}$. By the previous paragraph, $h$ lies in an absolutely
simple direct factor $G_{h}$ of $G$.
If $k_{h}$ is any local field extending $k$, and $\rho$
any $k_{h}$-rational, absolutely irreducible, finite-dimensional, nonzero representation of
$G$ factoring through $G_{h}$, then it follows from absolute simplicity that $G_{h} \cap
\ker \rho \subset Z(G_{h}) = 1$. Since a
finitely generated linear torsion group is finite, it follows that $\rho(\Gamma)$
contains nontorsion elements. Hence, the preceding paragraphs imply that there exists
such a representation $\rho_{h}$ with $\Omega_{+}(\rho_{h}, \Gamma)$ nonempty.
 
To summarize the proof so far, we have constructed a family $$\{\rho_{h}: G \to
\mathrm{GL}(E_{h})\mid h
\in F\}$$ of absolutely irreducible, finite-dimensional, nonzero linear representations of
$G$ rational over local extensions of $k$, factoring through absolutely simple direct
factors of $G$,
with the properties that
\begin{itemize}
\item[(i)] $\rho_{h}(h)$ is very proximal if $h \in F_{s}$;
\item[(ii)] $\langle\rho_{h}(h)\rangle \cong \langle h\rangle$ if
$h
\in F_{t}$;
\item[(iii)] $\rho_{h}(h)$ has a power which is a quasi-proximal unipotent if $h \in
F_{u}$; and
\item[(iv)] $\Omega_{+}(\rho_{h}, \Gamma) \neq \varnothing$ for all $h \in F$.
\end{itemize}
Let $P_{h}$ denote the compact Hausdorff topological space $P(E_{h})(k_{h})$.

Write $F_{i} = F_{s} \cup F_{u}$. For $h \in F_{i}$,
let us set the notation
\begin{equation}
\label{eqn:Krdef}
\mathrm{Kr}(h) = \begin{cases}
\mathrm{Cr}(\rho_{h}(h))& \text{if $h \in F_{s}$},\\
\bigcup_{\nu = 1}^{\infty}h^{\nu}\cdot B'(\rho_{h}(h_{u}))& \text{if $h \in F_{u}$}.
\end{cases}
\end{equation}
We note that the apparently infinite union on the right-hand side of~\eqref{eqn:Krdef}
is actually a finite union of projective hyperspaces by the hypothesis that $h_{s}$ is
torsion.

For convenience, we introduce the notation for $\gamma \in G$ that $\gamma_{h} =
\rho_{h}(\gamma)$. In analogy with the critical set defined by equation~\eqref{eqn:Krdef},
we set for $\gamma \in G$
\begin{equation*}
C(\gamma_{h}) = \begin{cases}
A(\gamma_{h}) \cup A(\gamma_{h}^{-1})& \text{if $A(\gamma_{h}) \neq E_{h}$},\\
\bigcup_{\nu = 1}^{\infty}\gamma^{\nu}\cdot B(\rho_{h}(\gamma_{u}))& \text{if $\gamma_{s}$ is
torsion and $\rho_{h}(\gamma_{u})$ quasi-proximal},
\end{cases}
\end{equation*}
noting again that the union on the right-hand side is a finite union of
singleton sets. When $h \in F_{i}$, we will write
$C(h) = C(h_{h})$ for simplicity.

Denote by $R$ the set of representations $\{\rho_{h}, \rho_{h}^{\vee}\mid h \in F\}$, where
we denote by $\rho^{\vee}$ the dual of the representation $\rho$. Since
$\rho^{\vee}(\gamma)$ is proximal if $\rho(\gamma^{-1})$ is, we see that
$\bigcap_{\rho \in R} \Omega_{0}(\rho, S)$ coincides with
$\bigcap_{h \in F} \Omega_{0}(\rho_{h}, S)$. Lemma~\ref{simulunnec}.(i) implies the
existence of a finite attractor family $\Phi \subset \bigcap_{h \in F} \Omega_{0}(\rho_{h},
S)$ for the family $R$. Upon conjugation, we may assume
\begin{align*}
C(h)&\subset P_{h}\setminus A'(\phi_{h}^{-1})\text{ and}\\
A(\phi_{h})&\subset P_{h}\setminus\mathrm{Kr}(h)
\end{align*}
all $h \in F_{i}$ and $\phi \in \Phi$.

For each $h \in F_{i}$, choose compact subsets $U_{h}$ and $K_{h} \subset
P_{h}$ satisfying
\begin{gather*}
K_{h} \subset P_{h}\setminus \bigcup_{\phi \in \Phi} A'(\phi_{h}^{-1})\qquad C(h)\subset
\Int K_{h}\\
U_{h} \subset P_{h} \setminus \mathrm{Kr}(h)\qquad \bigcup_{\phi \in \Phi} A(\phi_{h}) \subset \Int
U_{h}.
\end{gather*}
By Lemma~\ref{lem:3.8}.(i) there exist numbers $N_{h}$ such that $h^{z} \cdot U_{h}
\subset \Int K_{h}$, all $h \in F_{i}$ and $\abs{z} \geqslant N_{h}$. If $h \in F_{t}$,
denote by $N_{h}$ the order of $h$.

Upon replacing $\Phi$ by a high (elementwise) power of itself,
Corollary~\ref{attractfamily} implies the following: For any $\gamma \in \bigcap_{h \in F}
\Omega_{0}(\rho_{h}, H)$ there exists $\phi \in \Phi$ such that for all $h \in
F_{i}$ we have
\begin{gather}
\phi\cdot A(\gamma_{h}) \subset \Int U_{h}\label{attractor1}\\
\phi\cdot A'(\gamma_{h}^{-1}) \cap K_{h} = \varnothing. \label{attractor2}
\end{gather}

Now fix $\gamma \in \bigcap_{h \in F} \Omega_{0}(\rho_{h}, H)$. By virtue of
Lemma~\ref{lem:tits} and
Proposition~\ref{prop:dyslexictits}, there exists $y \in S$ satisfying
\begin{gather*}
y^{-1}\phi^{-1}\cdot C(h) \subset P_{h}\setminus A'(\gamma_{h}^{-1})\\
y\cdot C(\gamma_{h}) \subset P_{h}\setminus (\phi^{-1}\cdot \mathrm{Kr}(h))
\end{gather*}
for all $h \in F_{i}$ and all $\phi \in \Phi$, and
$$h^{\phi yz}\cdot A(\gamma_{h}) \subset P_{h}\setminus A'(\gamma_{h}^{-1})$$
for all $h \in F$, $\phi \in \Phi$, and $0 < \abs{z} < N_{h}$. Since $H$ is
normal in $S$, we may as well replace $\gamma$ by $\mbox{}^{y}\gamma$. Thus we
obtain an element $\gamma \in \bigcap_{h \in F} \Omega_{0}(\rho_{h}, H)$ satisfying
\begin{gather*}
\phi^{-1}\cdot C(h) \subset P_{h}\setminus A'(\gamma_{h}^{-1})\\
C(\gamma_{h}) \subset P_{h}\setminus (\phi^{-1}\cdot \mathrm{Kr}(h))
\end{gather*}
for all $h \in F_{i}$ and all $\phi \in \Phi$, and
$$h^{\phi z}\cdot A(\gamma_{h}) \subset P_{h}\setminus A'(\gamma_{h}^{-1})$$
for all $h \in F$, $\phi \in \Phi$, and $0 < \abs{z} < N_{h}$. Let $\phi \in \Phi$ be that element
satisfying conditions~\eqref{attractor1} and~\eqref{attractor2} for this choice of
$\gamma$.

On the other hand, if $h \in F_{i}$ and $\abs{z} \geqslant N_{h}$, we have $h^{z}\phi
\cdot A(\gamma_{h})\subset \Int h^{z} \cdot U_{h} \subset \Int K_{h}$. Then it follows
in this case by~\eqref{attractor2} that $(h^{z}\phi \cdot A(\gamma_{h}))\cap(\phi\cdot
A'(\gamma_{h}^{-1})) = \varnothing$. Thus if we write $\beta = \mbox{}^{\phi}\gamma$, then
$\beta\in
\bigcap_{h\in F} \Omega_{0}(\rho_{h}, H)$ satisfies the following conditions:
\begin{align}
C(h)&\subset P_{h}\setminus A'(\beta_{h}^{-1})\text{ and}\nonumber\\
C(\beta_{h})&\subset P_{h}\setminus\mathrm{Kr}(h)\nonumber
\intertext{for all $h \in F_{i}$; and}
h^{z}\cdot A(\beta_{h}) &\subset P_{h}\setminus A'(\beta_{h}^{-1})\label{supertrans}
\end{align}
for all $h \in F$ and all nontrivial powers $h^{z}$ of $h$.

By Lemma~\ref{lem:3.8}.(i) and Lemma~\ref{lem:quasiprox},
we have that $h^{z}\cdot U_{h} \subset \Int K_{h}$
for all $h \in F_{i}$ and for almost all $z \in \mathbb{N}$.
Equations~\eqref{attractor1}
and~\eqref{attractor2} imply that $A(\beta_{h}) \subset \Int U_{h}$ and $K_{h} \subset
P_{h}\setminus A'(\beta_{h}^{-1})$, $h \in F_{i}$.
Upon shrinking $U_{h}$,
enlarging $K_{h}$, and using~\eqref{supertrans} if necessary, we obtain compact sets $K_{h}^{+}, M_{h} \subset P_{h}$, $h \in
F_{i}$, with
\begin{gather*}
K_{h}^{+} \subset P_{h}\setminus \mathrm{Kr}(h)\qquad A(\beta_{h})\subset \Int K_{h}^{+}\\
M_{h } \subset P_{h} \setminus A'(\beta_{h}^{-1})\qquad C(h) \subset \Int M_{h}.
\end{gather*}
and  $h^{z}\cdot K_{h}^{+} \subset \Int M_{h}$ for all $h \in F_{i}$ and all $z \neq 0$.
For $h \in F_{t}$, by virtue of~\eqref{supertrans}, we can choose compact sets
$K_{h}^{+} \subset P_{h}$ such that 
\begin{align*}
A(\beta_{h}) &\subset \Int K_{h}^{+}\\
h^{z}\cdot K_{h}^{+} &\subset P_{h} \setminus A'(\beta_{h}^{-1})
\end{align*}
for all nontrivial powers $h^{z}$ of $h$. For the sake of consistency, when $h \in F_{t}$,
we shall write $M_{h} = \bigcup_{h^{z} \neq 1} h^{z}\cdot K_{h}^{+}$. With this
definition,
we have $M_{h}\subset P_{h} \setminus A'(\beta_{h}^{-1})$ for all $h \in F$.

Now by Proposition~\ref{prop:coset}, there exists $\delta \in \bigcap_{h \in F}
\Omega_{0}(\rho_{h}, Hg)$. By Lemma~\ref{lem:tits}
and the fact that $H$ is Zariski-dense, there
exists an element $x \in H$ such that for all $h \in F$ we have
\begin{gather*}
x\cdot A(\beta_{h}^{-1}) \subset P_{h} \setminus \mathrm{Cr}(\delta_{h})\\
x^{-1}\cdot C(\delta_{h}) \subset P_{h} \setminus A'(\beta_{h}).
\end{gather*}
Since $x \in H$ and $H$ is normal,
we have $\delta^{x} \in Hg$. Therefore, we may as well suppose that  $\delta = 
\delta^{x} \in \bigcap_{h \in F}
\Omega_{0}(\rho_{h}, Hg)$ satisfies the following conditions
\begin{gather}
A(\beta_{h}^{-1}) \subset P_{h} \setminus \mathrm{Cr}(\delta_{h})\label{eqn:delta2}\\
C(\delta_{h}) \subset P_{h} \setminus A'(\beta_{h}),\label{eqn:delta3}
\end{gather}
all $h \in F$.

By
equations~\eqref{eqn:delta2} and~ \eqref{eqn:delta3}, for each $h \in F$ and we can choose compact sets
$K_{h}^{-}$ and $L_{h}$ with
\begin{gather*}
K_{h}^{-} \subset P_{h}\setminus \mathrm{Cr}(\delta_{h})\qquad A(\beta_{h}^{-1})\subset \Int K_{h}^{-}\\
L_{h } \subset P_{h} \setminus A'(\beta_{h})\qquad C(\delta_{h}) \subset \Int L_{h}.
\end{gather*}
By Lemma~\ref{lem:3.8}.(i) there exists a number $\nu$ such that
\begin{align}
\beta^{-\nu}\cdot M_{h} &\subset \Int K_{h}^{-}\label{eqn:cdelta1}\\
\beta^{\nu}\cdot L_{h}&\subset \Int K_{h}^{+}\label{eqn:cdelta2}
\end{align}
for all $h \in F$. Finally, by Lemma~\ref{lem:3.8}.(i), there exists a number $\mu_{0}$ such that
$\delta^{\mu}\cdot K_{h}^{-} \subset \Int L_{h}$ for all $h \in F$, whenever $\abs{\mu} >
\mu_{0}$.

Consider the expression $\beta^{-\nu}h^{z}\beta^{\nu}\cdot L_{h}$. Using
equations~\eqref{eqn:cdelta1} and~\eqref{eqn:cdelta2}, we have
\begin{align*}
\beta^{-\nu}h^{z}\beta^{\nu}\cdot L_{h}
&\subset \Int \beta^{-\nu}h^{z}\cdot K_{h}^{+}\\
&\subset \Int \beta^{-\nu}\cdot M_{h}\\
&\subset \Int K_{h}^{-}.
\end{align*}
In particular, we have for $\mu > \mu_{0}$ and $j \neq 0$ that
$$\delta^{j\mu}\beta^{-\nu}h^{z}\beta^{\nu}\cdot L_{h} \subset \Int L_{h}.$$
Thus whenever $h \in F$ and $z \in \mathbb{Z}$ is such that $h^{z} \neq 1$, we have
$$(\delta^{\mu})^{j}(\beta^{-\nu}h\beta^{\nu})^{z}\cdot L_{h}
\subset \Int L_{h}.$$
Also, taking $z = 1$ in the above sequence of inclusions yields $$\beta^{-\nu}h\beta^{\nu}\cdot L_{h}
\subset \Int K_{h}^{-} \subset P_{h}\setminus \mathrm{Cr}(\delta_{h}),$$
whence $\beta^{-\nu}h\beta^{\nu} \cdot L_{h} \not\supset L_{h}$.
By Lemma~\ref{lem:pingpong}, we have for all $h \in F$ and $\mu > \mu_{0}$ that
$$\langle \delta^{\mu}, \beta^{-\nu}h\beta^{\nu}\rangle
\cong \langle\delta^{\mu}\rangle * \langle\beta^{-\nu}h\beta^{\nu}\rangle.$$
Conjugating both sides yields
$$\langle \beta^{\nu}\delta^{\mu}\beta^{-\nu}, h\rangle
\cong \langle\beta^{\nu}\delta^{\mu}\beta^{-\nu}\rangle * \langle h\rangle.$$
If we take $\mu >\mu_{0}$ such that $\mu \equiv 1 \mod{\abs{H\backslash\Gamma}}$
then $H\beta^{\nu}\delta^{\mu}\beta^{-\nu} = H\delta^{\mu} = Hg$, because $H$ was assumed
normal, $\beta \in H$, and $\delta \in Hg$.
\end{proof}
Finally, we return to the general case.
\begin{lemma}\label{lem:awayFromClass}
Let $C$ be a nonunipotent conjugacy class in a connected, reductive group $G$. If $\Gamma$
is a
finitely generated subgroup of $G$, then there exists a finite index subgroup $H < \Gamma$ which
does not meet $C$.
\end{lemma}
\begin{proof}
Since $\Gamma$ lies in the points of $G$ over a finitely generated field, the identity
map $\Gamma \to \Gamma$ extends to a continuous map from the profinite completion of $\Gamma$
into the points of $G$ over a local field. As the (Zariski) closure of $C$
does not contain the identity, it follows that the identity does not lie in the
profinite closure of $C \cap \Gamma$.
\end{proof}
\begin{proof}[Proof of Theorem~\ref{thm:subgroup}]
Let $H$ and simple direct factor of $G$. Denote by $\pi_{H}$ the projection $G
\to H$. If $H$ is neither of type $B_{n}$, $n \geqslant 2$, nor of type $G_{2}$, then
Remark~\ref{rem:conjtobigcell} implies that there exists no class of elements $h \in H$
with the property that $Z_{H}(h)^{\circ}$ is a simple, proper subgroup of full rank
in $H$. If $H$ is of of type $B_{n}$, $n \geqslant 2$, or of type $G_{2}$, then there
exists a unique such class, say $C_{H}$. In this case, $C_{H}$ is semisimple torsion.

It follows from Lemma~\ref{lem:awayFromClass} that there exists a subgroup $\Gamma_{0}$
of finite index in $\Gamma$ with the property that $\pi_{H}(\Gamma_{0}) \cap C_{H} =
\varnothing$ for every $H$ of type $B_{n}$, $n \geqslant 2$, or $G_{2}$. The result now
follows from Theorem~\ref{thm:profinitedensity} and Lemma~\ref{lem:goingUp}.
\end{proof}
\begin{proof}[Proof of Theorem~\ref{thm:profinitedensity}]
Let $G_{\pm 1}$ be as in the statement of Theorem~\ref{thm:zariski-densePPPs}. Note that
$G_{+1}$ lies in the product of the simply laced, absolutely simple direct factors of
$G$. Let $Hg$ be a coset of a finite-index subgroup $H \subset \Gamma$, $g\in \Gamma$.

Let $F_{r}$ be the set of all elements $h \in F$ whose semisimple parts $h_{s}$
satisfy the condition
$$h_{s}^{j} \in G_{+1} \text{ implies } j = 0.$$
Fix for the moment an element $h \in F_{r}$. Upon conjugating, we may assume $T$
contains $h_{s}$. Then writing $h = h_{-1}h_{+1}$ with $h_{-1} \in G_{-1}$, then some
eigenvalue $\tau$ of $h_{-1}$ has infinite multiplicative order. By Lemma~\ref{lem:hatch},
there exists a local field $k_{h}$
extending $k$ such that $\abs{\tau} \neq 1$. If $\Pi_{im} =
\Pi_{im}(h_{s})$ is as in \S\ref{sec:Pi-im}, then there exists a simply laced
component of the Dynkin diagram of $G$ which is not contained in $\Pi_{im}$. If $\rho$
is any $k_{h}$-rational, absolutely irreducible, finite-dimensional, nonzero representation of
$G$ which factors through the
corresponding absolutely simple
direct factor of $G$, then $\rho(h)$ lies outside every compact subgroup of $\rho G(k_{h})$.

Choose such a representation $\rho_{h}$ which is $\ell_{r(h)}$-large, where $r(h)$ is the rank of the simple factor of $G$ through which $\rho_{h}$ factors. 
Given an affine hyperspace $K$ of $\Lambda \otimes \mathbb{R}$, denote by $V_{h, K}$ the
projectivization of $\bigoplus_{\mu \in K} E^{\mu}$. By
Corollary~\ref{cor:localversionofproposition}, we have that the set $X$ of $u \in G$
satisfying $V_{h, K} \cap h^{u} \cdot V_{h, K'} = \varnothing$ for every pair $K$, $K'$ of
parallel affine hyperspaces of $\Lambda \otimes \mathbb{R}$ and every $h \in F_{r}$ is
nonempty and open.

Set $F_{p} = F \setminus F_{r}$. Denote by $Y$ the set of elements $\gamma \in \Gamma$
of infinite order such that the subgroup
$\langle \gamma, h\rangle$ generated by $\gamma$ and $h$ is canonically isomorphic to the
free product
$\langle \gamma \rangle * \langle h\rangle$, all $h \in F_{p}$. By Theorem~\ref{thm:zariski-densePPPs}, $Y$ is dense in the profinite topology on $\Gamma$. Since $\Gamma$ is finitely generated it follows from
Lemma~\ref{lem:YcapHg} that $Y \cap Hg$ is Zariski-dense. For convenience, denote by
$\phi$ the morphism $T \times G \to G$ given by $(t, u) \mapsto \mbox{}^{u}t$. Then
there exists an element $g \in Y \cap Hg$ in the image of $T(k) \times X$ under the
dominant morphism $\phi$.

Fix $h \in F_{r}$, and positive a number $c$. Writing $g = \mbox{}^{u}t$ for $u
\in X$ and $t \in T(k)$ and $K(h, c)$ for the affine hyperspace $\{\sum_{\mu \in \Lambda} \mu \otimes q_{\mu} \in
\Lambda
\otimes
\mathbb{R} \mid \sum_{\mu \in \Lambda} q_{\mu} \log \abs{\mu(h_{s})} = \log c\}$,
we observe that $A^{c}(\rho_{h}(t))$ coincides with $V_{h, K(h, c)}$. In particular,
$A^{c}(\rho_{h}(t)) \cap h^{u} \cdot A^{d}(\rho_{h}(t)) = \varnothing$ for every pair of
positive numbers $c$ and $d$ and every $h \in F_{r}$. Finally, we compute
$$A^{c}(\rho_{h}(g)) \cap h \cdot
A^{d}(\rho_{h}(g)) = u \cdot A^{c}(\rho_{h}(t)) \cap hu \cdot A^{d}(\rho_{h}(t))
=\varnothing$$
for every pair of positive numbers $c$ and $d$ and every $h \in F_{r}$.
\end{proof}

\providecommand{\bysame}{\leavevmode\hbox to3em{\hrulefill}\thinspace}


\begin{thebibliography}{10}

\bibitem{MR1348303}
H.~Abels, G.~A.~Margulis, and G.~A.~So{\u\i}fer, \emph{Semigroups containing
  proximal linear maps}, Israel J.\ Math.\ \textbf{91} (1995), no.~1--3, 1--30.

\bibitem{MR96a:22020}
M.~Bekka, M.~Cowling, and P.~de~la Harpe, \emph{Some groups whose reduced
  ${C}\sp *$-algebra is simple}, Inst.\ Hautes \'Etudes Sci.\ Publ.\ Math.\ (1994),
  no.~80, 117--134 (1995).

\bibitem{MR2001j:46078}
M.~Bekka and P.~de~la Harpe, \emph{Groups with simple reduced ${C}\sp
  *$-algebras}, Expo.\ Math.\ \textbf{18} (2000), no.~3, 215--230.

\bibitem{MR0372054}
A.~Borel, \emph{Linear representations of semi-simple algebraic groups},
  Algebraic geometry ({P}roc.\ {S}ympos.\ {P}ure {M}ath., {V}ol.\ 29, {H}umboldt
  {S}tate {U}niv., {A}rcata, {C}alif., 1974), Amer.\ Math.\ Soc., Providence,
  R.I., 1975, pp.~421--440.

\bibitem{MR1102012}
\bysame, \emph{Linear algebraic groups}, second ed., Graduate Texts in
  Mathematics, vol.\ 126, Springer-Verlag, New York, 1991.

\bibitem{MR0207712}
A.~Borel and J.~Tits, \emph{Groupes r\'eductifs}, Inst.\ Hautes \'Etudes Sci.\ 
  Publ.\ Math.\ (1965), no.~27, 55--150.

\bibitem{MR0316587}
\bysame, \emph{Homomorphismes ``abstraits'' de groupes alg\'ebriques simples},
  Ann.\ of Math.\ (2) \textbf{97} (1973), 499--571.

\bibitem{MR1727844}
N.~Bourbaki, \emph{Algebra {I}. {C}hapters 1--3}, Elements of Mathematics
  (Berlin), Springer-Verlag, Berlin, 1998, Translated from the French, Reprint
  of the 1989 English translation [ MR0979982 (90d:00002)].

\bibitem{MR1890629}
\bysame, \emph{Lie groups and {L}ie algebras. {C}hapters 4--6}, Elements of
  Mathematics (Berlin), Springer-Verlag, Berlin, 2002, Translated from the 1968
  French original by Andrew Pressley.

\bibitem{MR1611814}
C.~Chevalley, \emph{Certains sch\'emas de groupes semi-simples}, S\'eminaire
  {B}ourbaki, {V}ol.\ 6, Soc.\ Math.\ France, Paris, 1995, pp.~Exp.\ No.\ 219,
  219--234.

\bibitem{MR87b:22007}
P.~de~la~Harpe, \emph{Reduced ${C}\sp \ast$-algebras of discrete groups which
  are simple with a unique trace}, Operator algebras and their connections with
  topology and ergodic theory (Bu\c steni, 1983), Springer, Berlin, 1985,
  pp.~230--253.

\bibitem{1123.22004}
\bysame, \emph{On simplicity of reduced ${C}\sp \ast$-algebras of
  groups}, Bull.\ London Math.\ Soc.\ \textbf{39} (2007), no.~1, 1--26.

\bibitem{MR0352328}
H.~Furstenberg, \emph{Boundary theory and stochastic processes on homogeneous
  spaces}, Harmonic analysis on homogeneous spaces (Proc.\ Sympos.\ Pure Math.,
  Vol.\ XXVI, Williams Coll., Williamstown, Mass., 1972), Amer.\ Math.\ Soc.,
  Providence, R.I., 1973, pp.~193--229.

\bibitem{MR1040268}
I.~Ya.~Gol{\mlprime}dshe{\u\i}d and G.~A.~Margulis, \emph{Lyapunov exponents of
  a product of random matrices}, Uspekhi Mat.\ Nauk \textbf{44} (1989),
  no.~5(269), 13--60.

\bibitem{MR0207713}
B.~Kostant, \emph{Groups over {$Z$}}, Algebraic {G}roups and {D}iscontinuous
  {S}ubgroups ({P}roc.\ {S}ympos.\ {P}ure {M}ath., {B}oulder, {C}olo., 1965),
  Amer.\ Math.\ Soc., Providence, R.I., 1966, pp.~90--98.

\bibitem{MR613853}
G.~A.~Margulis and G.~A.~So{\u\i}fer, \emph{Maximal subgroups of infinite index
  in finitely generated linear groups}, J.\ Algebra \textbf{69} (1981), no.~1,
  1--23.

\bibitem{MR82c:22010}
W.~L.~Paschke and N.~Salinas, \emph{${C}\sp{\ast} $-algebras associated with
  free products of groups}, Pacific J.\ Math.\ \textbf{82} (1979), no.~1,
  211--221.

\bibitem{MR1642713}
T.~A.~Springer, \emph{Linear algebraic groups}, second ed., Progress in
  Mathematics, vol.~9, Birkh\"auser Boston Inc., Boston, MA, 1998.

\bibitem{MR0180554}
R.~Steinberg, \emph{Regular elements of semisimple algebraic groups}, Inst.\ 
  Hautes \'Etudes Sci.\ Publ.\ Math.\ (1965), no.~25, 49--80.

\bibitem{MR0430094}
\bysame, \emph{On the desingularization of the unipotent variety}, Invent.\ 
  Math.\ \textbf{36} (1976), 209--224.

\bibitem{MR0277536}
J.~Tits, \emph{Repr\'esentations lin\'eaires irr\'eductibles d'un groupe
  r\'eductif sur un corps quelconque}, J.\ Reine Angew.\ Math.\ \textbf{247}
  (1971), 196--220.

\bibitem{MR44:4105}
\bysame, \emph{Free subgroups in linear groups}, J.\ Algebra \textbf{20} (1972),
  250--270.

\end{thebibliography}
\end{document}